\documentclass[preprint,11pt]{elsarticle}
\usepackage{subcaption}
\usepackage{amsmath}
\usepackage{multirow}

\DeclareMathOperator*{\argmin}{arg\,min}
\usepackage{amsfonts}
\usepackage{graphics}
\usepackage[%nohead,%nofoot
letterpaper, top=1.2in, bottom=1.2in, left=1.2in, right=1.2in]{geometry}
\usepackage{xcolor}
\usepackage{amssymb}
\usepackage{psfrag}
\usepackage{graphicx}
\usepackage[utf8]{inputenc}
\usepackage{bm}
\usepackage{mathtools}
\usepackage{epstopdf}
\usepackage[colorlinks,linkcolor=black,anchorcolor=black,citecolor=black]{hyperref}
\usepackage{cleveref}
\usepackage{setspace}
\usepackage{amsmath,amssymb,bbm}
\usepackage[utf8]{inputenc}
\usepackage{url}
\usepackage{color}
\usepackage{booktabs}
\usepackage[linesnumbered]{algorithm2e}
\usepackage{placeins}
\newtheorem{remark}{Remark}
\usepackage[title]{appendix}

\DeclareMathAlphabet\mathbfcal{OMS}{cmsy}{b}{n}

\makeatletter
\def\ps@pprintTitle{%
  \let\@oddhead\@empty
  \let\@evenhead\@empty
  \def\@oddfoot{\reset@font\hfil\thepage\hfil}
  \let\@evenfoot\@oddfoot
}
\makeatother

\begin{document}

\begin{frontmatter}

% \title{Data-driven discovery of highly noisy dynamical systems via Mixed Integer Optimization}

\title{Compressive-sensing-assisted mixed integer optimization for 
dynamical system discovery with highly noisy data}

% Use letters for affiliations, numbers to show equal authorship (if applicable) and to indicate the corresponding author
\author[a]{Zhongshun Shi}
\author[a]{Hang Ma}
\author[b]{Hoang Tran}
\author[b]{Guannan Zhang}

\address[a]{Department of Industrial and Systems Engineering, University of Tennessee Knoxville, Knoxville, TN 37996}
\address[b]{Computer Science and Mathematics Division, Oak Ridge National Laboratory, Oak Ridge, TN 37831}

%

% Please give the surname of the lead author for the running footer
%\leadauthor{Shi}

%% Please add a significance statement to explain the relevance of your work
%\significancestatement{Governing equations are critical for shaping our understanding of the physical word. There are ubiquitous dynamical systems in which the governing equations are unknown by human, yet large datasets can be collected through the evolution. In the context of science and engineering, a model is worth a thousand datasets. This work proposes an effective approach for uncovering the governing equations from data. Via a modern discrete optimization lens, the proposed method addresses a key challenge faced by prior research, that is, exact recovery of governing equations under large noise. It also addresses the high dimensional limitation. Our work sketches a path for exact recovery of governing equations under large noise.
%}

%% Please include corresponding author, author contribution and author declaration information
% \authorcontributions{Z.S., H.M., H.T. and G.Z. designed research; Z.S. and H.M. performed research;
% Z.S., H.M., H.T., and G.Z. analyzed data; and Z.S., H.M., H.T., and G.Z. wrote the paper.}
%\authordeclaration{The authors declare no conflict of interest.}
%\correspondingauthor{\textsuperscript{1}To whom correspondence should be addressed. E-mail: tzshi@utk.edu and zhangg@ornl.gov}

% At least three keywords are required at submission. Please provide three to five keywords, separated by the pipe symbol.

%
\begin{abstract}
The identification of governing equations for dynamical systems is everlasting challenges for the fundamental research in science and engineering. Machine learning has exhibited great success to learn and predict dynamical systems from data.
However, the fundamental challenges still exist: discovering the exact governing equations from highly noisy data.
In present work, we propose a compressive sensing-assisted mixed integer optimization (CS-MIO) method to make a step forward from a modern discrete optimization lens. 
In particular, we first formulate the problem into a mixed integer optimization model.
The discrete optimization nature of the model leads to exact variable selection by means of cardinality constraint, and hereby powerful capability of exact discovery of governing equations from noisy data.
Such capability is further enhanced by incorporating compressive sensing and regularization techniques for highly noisy data and high-dimensional problems.
%variable selection and capability for denoise
%The formulated MIO is NP-hard, for which there is a belief for decades in many scientific communities that such discrete optimization problems were intractable. 
%Nevertheless, recent advancement in algorithms and hardware significantly increase the capability of discrete optimization and make it probable to be computationally tractable for such hard problems. successfully demonstrate the effectiveness and great potential to thrust system dynamics identification with discrete optimization.
%We fully leverage the advances of modern discrete optimization techniques to solve the MIO model to demonstrate the applicability and numerical tractability in system dynamics identification.
The case studies on classical dynamical systems have shown that CS-MIO can discover the exact governing equations from large-noise data, with up to two orders of magnitude larger noise comparing with state-of-the-art method. We also show its effectiveness for high-dimensional dynamical system identification through the chaotic Lorenz 96 system.
\end{abstract}


%\dates{This manuscript was compiled on \today}
%\doi{\url{www.pnas.org/cgi/doi/10.1073/pnas.XXXXXXXXXX}}

%\maketitle
% \tnotetext[fn]{dfad}
\tnotetext[fn]{{\bf Notice}: This manuscript has been authored by UT-Battelle, LLC, under contract DE-AC05-00OR22725 with the US Department of Energy (DOE). The US government retains and the publisher, by accepting the article for publication, acknowledges that the US government retains a nonexclusive, paid-up, irrevocable, worldwide license to publish or reproduce the published form of this manuscript, or allow others to do so, for US government purposes. DOE will provide public access to these results of federally sponsored research in accordance with the DOE Public Access Plan.}

\end{frontmatter}
%\thispagestyle{firststyle}
%\ifthenelse{\boolean{shortarticle}}{\ifthenelse{\boolean{singlecolumn}}{\abscontentformatted}{\abscontent}}{}
% If your first paragraph (i.e. with the \dropcap) contains a list environment (quote, quotation, theorem, definition, enumerate, itemize...), the line after the list may have some extra indentation. If this is the case, add \parshape=0 to the end of the list environment.

%\dropcap{T}his PNAS journal template is provided to help you write your work in the correct journal format. Instructions for use are provided below.

%Note: please start your introduction without including the word ``Introduction'' as a section heading (except for math articles in the Physical Sciences section); this heading is implied in the first paragraphs.

%\dropcap{T}his PNAS journal template is provided to help you write your work in the correct journal format. 

\section{Introduction}\label{sec:intro}

Governing equations of the ubiquitous dynamical systems are of critical significance to shape our comprehension of the physical world. Traditional regime of obtaining the equations respects to the mathematical or physical derivations following the first principles, including conservation laws, mathematical symmetry and invariants. 
This paradigm, however, might be intractable for dealing with many complex phenomena. 
With the availability of large dataset due to the advances of sensors and technology, a new paradigm of discovering governing equations purely from data has been evolved. 
Machine learning plays the pivotal role under this paradigm with a wide scope of methods including 
symbolic regression \cite{bongard2007automated, schmidt2009distilling}, Gaussian processes \cite{raissi2017machine}, deep neural network \cite{doi:10.1073/pnas.1814058116,RUDY2019483,doi:10.1137/20M1342859,CHEN2021110362,raissi2018multistep}, Bayesian inference \cite{doi:10.1098/rspa.2018.0305}, etc.

% Anther prominent research focuses on symbolic regression mainly using genetic programming \cite{koza1994genetic}, and has been successfully applied to learn free form nonlinear dynamics \cite{bongard2007automated, schmidt2009distilling,ly2012learning}. Although its flexible representations can express many highly nonlinear dynamics, the symbolic regression may suffer from the time expensiveness, overfitting and inability of identifying the coefficients.
% Many other research explores from statistical and probabilistic learning \cite{yuan2019data,raissi2017machine,raissi2018hidden}, Spline network \cite{sun2021physics}, compressive sensing and neural network\cite{lusch2018deep,mardt2018vampnets}.
% The spatial-temporal correlation and large noise in system dynamics identification, however, pose huge challenges to these MIO based algorithms for which novel studies are still of need. 

% Other techniques: 
% Gaussian processes \cite{raissi2017machine}, deep neural network \cite{RUDY2019483}, Bayesian inference \cite{doi:10.1098/rspa.2018.0305}, classical numerical method \cite{Kang_IDENT21}. 

Even though neural networks have been proved to be an effective tool in learning and predicting trajectories dynamical systems, it is often challenging to extract new physical laws out of neural network models. Thus, this work focuses on another thrust of data-driven discovery of governing equations exploit sparse regression approaches \cite{brunton2016discovering, rudy2017data, brunton2016sparse, loiseau2018sparse, champion2019data, mangan2019model}. Studies along this path typically construct a large library of candidate terms and eventually transform into a sparse regression problem, grounded on the realistic assumption that only parsimonious terms are active in the governing equations. 
The breakthrough work by \cite{brunton2016discovering} introduced a novel architecture called Sparse Identification of Nonlinear Dynamical Systems (SINDy), which used a sequential threshold least squares (or ridge regression \cite{de2020pysindy}) to advocate sparsity.
The SINDy framework is impressive for its succinct but useful rationale, that is, the sparsity is essentially incurred by the penalty on coefficients. 
On this regard, studies have been conducted from many perspectives, including Lasso-based approach with a dictionary of partial derivatives \cite{doi:10.1098/rspa.2016.0446}, $\ell_{2,1}$ norm for data with highly corrupted segments \cite{doi:10.1137/16M1086637}, weak SINDy and discretization accounting for white noise \cite{doi:10.1137/20M1343166}, integral formulation of the differential equation\cite{PhysRevE.96.023302}, weak formulation with the orthogonal matching pursuit \cite{Pantazis2019AUA} and compressed sensing technique \cite{Wang2011PredictingCI}, 
to name a few.
However, these methods perform terms selection essentially via imposing penalty on coefficients, subject to which they are usually sensitive to noise, and unable to control the exact level of sparsity for differential equations.

%which shares certain similarity with the regularization techniques, including Lasso \cite{tibshirani1996regression}, Elastic-net \cite{zou2005regularization} and non-convex regularization \cite{fan2001variable,mazumder2011sparsenet}

From the perspective of discrete optimization, this sparse regression problem can be formulated as a Mixed Integer Optimization (MIO) model which is to identify a combination of $k$ terms from a pool of $p$ candidates and simultaneously regress the coefficients.
This $\ell_0$ norm constrained MIO problem is non-convex and $\mathcal{NP}$-hard \cite{natarajan1995sparse}, corresponding to the best subset selection in the larger statistics community \cite{miller2002subset, bertsimas2016best}.
The $\mathcal{NP}$-hardness of the problem has contributed to the belief that discrete optimization problems were intractable \cite{bertsimas2020rejoinder}.
For this reason, plenty of impressive sparsity-promoting techniques have focused on computationally feasible algorithms for solving the approximations, including Lasso \cite{tibshirani1996regression}, Elastic-net \cite{zou2005regularization}, non-convex regularization \cite{fan2001variable,mazumder2011sparsenet} and stepwise regression \cite{draper1998applied}. 
These approximations induce obscure sparsity via regularization that often includes a large set of active terms (many are correlated terms and the coefficients are shrunken to zero to avoid overfitting) in order to deliver good prediction. That is, regularization is used for both variable selection and shrinkage. In contrast, the MIO based exact method allows to control the exact level of sparsity via setting the value of $k$. When MIO based exact method decides to select a term, it purely takes it in without any shrinkage on the coefficients thereby draining the effect of its correlated terms \cite{bertsimas2016best}. 
Indeed, there is nothing more important than correct terms selection in the identification of governing equations. Although existing methods lean heavily on the sparsity-promoting parameters to achieve indirect terms selection, domain-educated researchers and practitioners actually might have an intuition for the ground truth $k$.
This motivates us in present work to enable independent and direct terms selection and coefficients shrinkage for solving the sparse regression problem in governing equations identification.

We propose a discrete optimization based method for exact recovery of differential equations under large noise. Our method takes advantages of the nature of discrete optimization in the means of cardinality constraints for terms selection, and is able to separately control the exact sparsity of the governing equations and estimate the associated coefficients. The powerful capability of terms selection is the cornerstone for exact recovery under large noise in the data, and is further enhanced by combining compressive sensing and regularization techniques for large noise and high dimensional problems.
We demonstrate the capability of our method with a wide variety of examples from \cite{brunton2016discovering}, including the chaotic Lorenz 3 system, the fluid dynamics of vortex shedding behind a cylinder, and two dynamical systems with bifurcations. In addition, we test on the famous high-dimensional Lorenz 96 system. 
Our results show the proposed method can recover exact governing equations with up to two orders of magnitude larger noise comparing with state-of-the-art method. This shows the modern discrete optimization is significantly effective for identifying governing equations from noisy and high-dimensional data.

\section{Problem setting}

In this section, we describe the problem setting of data-driven discovery of dynamical system. In addition, we introduce the highly noisy data setting for the problem.

\subsection{Data-driven discovery of dynamical system}

%Refer to 2019. \note{Technical intruduction of SINDy, and then give the gap - MIO.}
We introduce the data-driven dynamical system discovery problem from the perspective of sparse recovery \cite{brunton2016discovering}. We define $\mathcal{J}=\{1,2,\cdots,J\}$ for any $J\in \mathbb{Z}_+$ throughout this paper. 
We consider the following dynamical system consisting of $J$ state variables, i.e.,
\begin{align}\label{eq:ode}
\frac{d}{dt}\mathbf{x}(t)=\mathbf{f}(\mathbf{x}(t)),
\end{align}
where the vector $\mathbf{x}(t)=[x_1(t),\cdots,x_J(t)] \in \mathbb{R}^{1\times J}$ denotes the state of a system at time $t$, and $\mathbf{f}(\mathbf{x}(t)) = [{f}_1(\mathbf{x}(t)), \ldots, {f}_J(\mathbf{x}(t)) ]\in \mathbb{R}^{1\times J}$ represents the forcing term with ${f}_j(\mathbf{x}(t))$ being the forcing term of the $j$-th state variable $x_j$ for $j\in \mathcal{J}$. 
% For general aspects of variables $\mathbf{x}(t)$, we omit $t$ in the following. 
%
A dictionary $\bm \theta (\mathbf{x})$ consisting a total of $P$ terms, denoted by
\begin{equation}\label{dict}
    \bm \theta (\mathbf{x}) := [\theta_1(\mathbf{x}), \theta_2(\mathbf{x}),  \ldots,\theta_P(\mathbf{x})],
\end{equation}
which consists of nonlinear combinations of state $\mathbf{x}$ that can be candidate terms in $\mathbf{f}$. 
For example, $\bm \theta (\mathbf{x})$ may consist of polynomial, and trigonometric terms of $\mathbf{x}$. 
% That is $
%     \bm \theta (\mathbf{x})=[1,~\mathbf{x}, ~\mathbf{x}^{P_2}, ~\mathbf{x}^{P_3}, ~\cdots, ~ \sin{(\mathbf{x})}, ~\cos{(\mathbf{x})}],
% $
% where $\mathbf{x}^{P_2}, ~\mathbf{x}^{P_3}$ denote the higher polynomial terms with corresponding order.
%For instance, $\mathbf{x}^{P_2}=[x_1^2, ~ x_1 x_2, ~ \cdots, ~ x_2^2 ,~ \cdots ,~ x_m^2]$. 
Each term of $\bm \theta (\mathbf{x})$ represents a candidate term for right-hand side of Eq.~\eqref{eq:ode}. 

This work is based upon two assumptions that were also used in  \cite{brunton2016discovering}. The first is that we assume the right-hand-side $\mathbf{f}$ in Eq.~\eqref{eq:ode} lives in the function space expanded by the dictionary $\bm \theta (\mathbf{x})$ in Eq.~\eqref{dict}.
In other words, there exists a coefficient matrix $\bm \Xi := [\bm \xi_1, \ldots, \bm \xi_J] \in \mathbb{R}^{P\times J}$ such that
\begin{equation}\label{eq:sparse}
f_j(\mathbf x) =  \bm \theta(\mathbf{x}) \cdot \bm \xi_j, \;\; \text{ for } j = 1, \ldots, J.
\end{equation}
We remark that the definition of the dictionary $\bm \theta (\mathbf{x})$ would require domain knowledge about the specific scientific problem, in order to ensure that all the terms in $\mathbf{f}(\mathbf x)$ are included in $\bm \theta(\mathbf{x})$. Since this work is to study sparse recovery of $\mathbf{f}(\mathbf{x})$, how to properly choosing $\bm{\theta}(\mathbf{x})$ to ensure Eq.~\eqref{eq:sparse} is out of the scope of this paper. The second assumption is that the forcing term $\mathbf{f}$ consists of only a few terms, i.e., very sparse in the function space expanded by the dictionary $\bm{\theta}(\mathbf{x})$, regardless of the dimensionality $J$. 
%
% We develop an MIO based method, denoted by Discrete Optimization for Nonlinear Dynamics Identification (DONDI), to determine the underlying nonlinear terms. 
% Although the choice of active terms in $\mathbf{f}$ from $\bm \theta (\mathbf{x})$ is combinatorially many, we believe that only a few of these terms are active in each equation.
Specifically, to indicate the presence of each term of $\bm{\theta}(\mathbf{x})$ in the right hand side $\mathbf{f}$, we introduce the following indicator matrix
\begin{equation}
   \bm \Gamma :=[\bm \gamma_1, \cdots, \bm \gamma_J] =
   \begin{bmatrix}
   \gamma_{11}& \cdots& \gamma_{1J}\\
   \vdots & \ddots & \vdots\\
   \gamma_{P1}& \cdots& \gamma_{PJ}\\
   \end{bmatrix},\qquad
   \gamma_{pj} := \left\{
   \begin{aligned}
   & 1, \quad\text{ if $f_j(\mathbf{x})$ includes $\theta_p(\mathbf{x})$},\\
   & 0, \quad\text{ otherwise, }
   \end{aligned}
   \right.
\end{equation}
where $\bm \gamma_j = (\gamma_{1j}, \ldots, \gamma_{Pj})^{T}\in \mathbb{B}^{P\times 1}$ and $\mathbb{B}$ is the Boolean domain $\mathbb{B}=\{0,1\}$. 
% If the $p$-th term of $\bm \theta (\mathbf{x})$ is active in the $i$-th equation, $\bm \gamma_i^j=1$; otherwise $\bm \gamma_i^j=0$.
Moreover, 
we denote the number of active terms in $\mathbf{f} = [f_1, \ldots, f_J]$ by a vector 
\begin{equation}
    \bm k=[k_1,~ k_2,~ \cdots,~ k_J],
\end{equation}
where $k_j$ is the number of non-zeros in $\bm \gamma_j$. 
%Provided that $\mathbf{f}_i$ is parsimonious, $k_i$ is realistically chosen to be small number either by a prior domain knowledge or the expected sparsity for the governing equations in the dynamical systems of interests.
When the dynamical system satisfies the above two assumptions, Eq.~\eqref{eq:ode} can be written as
\begin{equation}\label{eq:ode1}
\dot{\mathbf{x}}(t)= \bm \theta(\mathbf{x}(t)) (\bm{\Gamma} \circ \bm \Xi)  ,
\end{equation}
where $\bm{\Gamma} \circ \bm \Xi$ is the element-wise product (Hadamard product) of $\bm \Gamma$ and $\bm \Xi$.

\subsection{The noisy data}

The state $\mathbf x$ and its time derivative $\dot{\mathbf x}$ can be measured and collected at a series of time instants $t_1,t_2,\dots,t_N$. 
With the measurements of $\mathbf{x}(t)$ and $\dot{\mathbf{x}}(t)$, we will be given two data matrices, denoted by $\mathbf{X}\in \mathbb{R}^{N\times J}$ and 
$\dot{\mathbf{X}}\in \mathbb{R}^{N\times J}$, of the following forms, 
\begin{equation}
\small
    \mathbf{X} =
    % \left[
    % \begin{array}{c}
    %      \mathbf{x}^{T}(t_1)  \\
    %      \mathbf{x}^{T}(t_2)  \\
    %      \vdots \\
    %      \mathbf{x}^{T}(t_n)  \\
    % \end{array}
    % \right]
    \left[
    \begin{array}{cccc}
         x_1(t_1)  & \cdots & x_J(t_1)\\
         x_1(t_2)  & \cdots & x_J(t_2)\\
         \vdots      & \ddots & \vdots  \\
         x_1(t_N)  & \cdots & x_J(t_N)\\
    \end{array}
    \right],
   \quad\text{ and }\quad
     \dot{\mathbf{X}} 
    %  =\left[
    % \begin{array}{c}
    %      \dot{\mathbf{x}}^{T}(t_1)  \\
    %      \dot{\mathbf{x}}^{T}(t_2)  \\
    %      \vdots \\
    %      \dot{\mathbf{x}}^{T}(t_N)  \\
    % \end{array}
    % \right]
    =\left[
    \begin{array}{cccc}
         \dot{x}_1(t_1)  & \cdots & \dot{x}_J(t_1)\\
         \dot{x}_1(t_2)  & \cdots & \dot{x}_J(t_2)\\
         \vdots     & \ddots & \vdots  \\
         \dot{x}_1(t_N)  & \cdots & \dot{x}_J(t_N)\\
    \end{array}
    \right],
\end{equation}
where the measurements of $\dot{\mathbf{x}}(t)$ can be numerically approximated 
using the data $\mathbf{X}$ if $\dot{\mathbf{x}}(t)$ is not directly measurable. %
% A data matrix $\mathbf{X}\in \mathbb{R}^{N\times J}$ is formed as $\mathbf{X} =[\mathbf{x}(t_1) ~ \mathbf{x}(t_2) ~ \cdots ~ \mathbf{x}(t_n)]^T$, and the time derivative $\dot{\mathbf{x}}(t)$ is either measured or approximated numerically from $\mathbf{x}(t)$, that is, $\mathbf{\dot{X}}=[\mathbf{\dot{x}}(t_1) ~ \mathbf{\dot{x}}(t_2) ~ \cdots ~ \mathbf{\dot{x}}(t_n)]^T$. 
In practice, the measurements $\mathbf{X}$ and $\dot{\mathbf{X}}$ are usually corrupted with random noises, so that the matrices $\mathbf{\Gamma}$ and $\mathbf{\Xi}$ in Eq.~\eqref{eq:ode1} need to be recovered with noisy data, denoted by
\begin{equation}\label{eq:noisy_data}
    \mathbf{X}^{\rm noisy} := \mathbf{X} + \mathbfcal{U}
    \;\;\text{ and }\;\; \dot{\mathbf{X}}^{\rm noisy} := \dot{\mathbf{X}} + \mathbfcal{V},
\end{equation}
where $\mathbfcal{U} \in \mathbb{R}^{N\times J}$ and $\mathbfcal{V} \in \mathbb{R}^{N\times J}$ are additive noise.

Evaluating the library $\bm \theta (\mathbf{x})$ at each data point in $\mathbf{X}^{\rm noisy}$, we can construct an augmented data matrix, denoted by $\bm \Theta (\mathbf{X}^{\rm noisy})$, consisting of candidate nonlinear functions of the columns of $\mathbf{X}^{\rm noisy}$.
For ease of notation, we use $\bm \Theta^{\rm noisy}$ instead of $\bm \Theta (\mathbf{X}^{\rm noisy})$ in the following.
Since there are $P$ terms in $\bm \theta (\mathbf{x})$, the matrix $\bm \Theta ^{\rm noisy} \in \mathbb{R}^{N\times P}$ is represented by
\begin{equation}\label{eq:dict}
    \bm \Theta^{\rm noisy} := [
         \theta_1(\mathbf{X}^{\rm noisy}),   \cdots, \theta_P(\mathbf{X}^{\rm noisy})].
\end{equation}
%
% The goal of dynamical system discovery is to identify the correct terms in $\bm \theta (\mathbf x)$ for $\mathbf{f}(\mathbf{x})$ by minimizing the error 
% \begin{equation}\label{eq:sparse_reg_noise} 
%     \dot{\mathbf{X}}^{\rm noisy} - \bm \Theta(\mathbf{X}^{\rm noisy})\, (\bm{\Gamma} \circ \bm \Xi),
% \end{equation}
% where $\bm \Xi=[\bm \xi_1 ~~ \bm \xi_2 ~~ \cdots \bm \xi_J]$, $\bm \xi_j\in \mathbb{R}^{P\times 1}$ is the vector of coefficients to be determined for $\mathbf{f}_j$, and $\mathbf{Z}\in\mathbb{R}^{N\times J}$ is the noise matrix.
%
Similar to the standard SINDy method in \cite{brunton2016discovering}, we assume the entries of the noise matrices $\mathbfcal{U}$ and $\mathbfcal{V}$ in Eq.~\eqref{eq:noisy_data} are independent and identically distributed (i.i.d.) Gaussian random variables with zero mean and standard deviation $\sigma$. In this work, we are particularly interested in the scenario with relatively large standard deviation of the noises, i.e., low signal-to-noise ratio. Details about the definition of the noises are give in Section \ref{sec:ex}.

The goal of sparse recovery of the dynamical system in Eq.~\eqref{eq:ode1} is to correctly identify $\bm \Gamma$ and calculate the non-zero elements of $\bm \Xi$ from measurement data of $\mathbf{x}$ and/or $\dot{\mathbf{x}}$. As discussed in Section \ref{sec:intro}, existing work on sparse recovery of dynamical systems, e.g., \cite{brunton2016discovering, rudy2017data, brunton2016sparse, loiseau2018sparse, champion2019data, mangan2019model,Wang2011PredictingCI}, perform term selection and promote sparsity by imposing penalties on the coefficients. In other words, these methods try to recover the product $\bm \Gamma \circ \bm \Xi$ as a whole. Despite the success of these methods, they are usually very sensitive to the noise in the measurement data. When the signal-to-noise ratio is low, the method like SINDy may fail to identify the correct terms of $\mathbf{f}$ in the dictionary $\bm \theta(\mathbf x)$. The motivation of this work is to recover $\bm \Gamma$ and $\bm \Xi$ separately, where $\bm \Gamma$ is recovered by solving a compressive-sensing-assisted mixed integer optimization, in order to identify the correct terms in the case of having data with low signal-to-noise ratio.

\section{The compressive sensing-assisted mixed integer optimization method}\label{sec:method}
This section describes the details of the proposed method. Specifically, a linear regression model subject to sparsity constraints for Eq.~\eqref{eq:ode1} can be set up as follows,
\begin{equation}\label{eq:sparse_reg_noise_k} 
    \min_{\bm \Gamma,\bm \Xi}\;\big\|\dot{\mathbf{X}}^{\rm noisy} -  \bm \Theta^{\rm noisy} (\bm \Gamma \circ \bm \Xi)\big\|_2^2,~~ s.t.~ \bm \Gamma^T \bm e \le k^{\max},
\end{equation}
where $\bm e \in \mathbb{R}^{P\times 1}$ is a vector with all the entries to be one, such that the product $\bm \Gamma^T \bm e$ are exactly the cardinality constraints to indicate the active terms in each equation, and $\bm k^{\max} = [k^{\max}_1, \ldots ,k^{\max}_J]$ consists of the maximum allowable sparsity for the $J$ components. 
The main idea of the CS-MIO method, is to separately identify the physical terms (i.e., $\bm \Gamma$) the and the corresponding coefficients (i.e., $\bm \Xi$) in a two stage manner. The indicator matrix $\bm \Gamma$ is determined by mixed integer optimization. Once $\bm \Gamma$ is chosen, we can estimate the corresponding components of $\bm \Xi$ using the standard least-squares method. Nevertheless, when the size of the original dictionary, i.e., the number of columns of $\bm \Gamma$, is large, it is computationally intractable for the state-of-the-art MIO algorithms. To resolve this issue, we propose to use compressive sensing, i.e., $\ell_1$ minimization, to reduce the size of the dictionary to the extent that can be handled by MIO algorithms.

% In particular, we first identify the physical terms in the first stage using the standardized data which removes the impacts from the original scale and unit. In the second stage, we identify the coefficients using original data.

% In this section, we identify $\mathbf{\Gamma}$ by first using the standardized data and then $\mathbf{\Xi}$ by the original data.

% \begin{equation}\label{eq:sparse_reg_noise_k} 
%     \widetilde{\mathbf{\dot{X}}}^{\rm noisy} =  \widetilde{\bm \Theta}^{\rm noisy} (\bm \Gamma\circ\bm \Xi}),~~ s.t.~ \bm \Gamma^T\bm e=\bm k,
% \end{equation}

In the rest of this section, we take the $j$-th component of $\bm x$ in  Eq.~\eqref{eq:ode} as an example in the following derivation, which means we intend to  use the $j$-th column of the data matrices  $\mathbf{{X}}^{\rm noisy}$ and  $\dot{\mathbf{X}}^{\rm noisy}$ to infer the $j$-th columns of $\bm \Gamma $ and $ \bm \Xi$. For notational simplicity, we omit the subscript $j$ and use $\dot{\bm x}^{\rm noisy}$, $\bm \gamma$, $\bm \xi$ to represent the $j$-column of $\dot{\mathbf{X}}^{\rm noisy}$, $\bm \Gamma $ and $ \bm \Xi$, respectively.

\subsection{Compressive sensing for reducing the size of the dictionary $\mathbf{\Theta}^{\rm noisy}$}\label{sec:CS}
The goal of this subsection is to reduce the size of the original dictionary $\bm \Theta^{\rm noisy}$ in Eq.~\eqref{eq:dict}, so that the modern integer optimization solvers, e.g., \texttt{CPLEX} or \texttt{GUROBI}, can be used to determine the indicator vector $\bm \gamma$.
%
% When the size $P$ of the dictionary is large, it is computationally intractable to 
% Formulation \eqref{MIO2} might have some computational issues when $P$ is large. 
% In this case, it can be time consuming for solving the combinatorial optimization problem using modern integer optimization solvers like \texttt{CPLEX} or \texttt{GUROBI}. 
To this end, we first solve the following $\ell_1$ minimization problem:
\begin{equation}\label{eq:cs}
    \bm \xi^{\rm CS} = \argmin_{\bm \xi}  ||{\dot{\bm x}}^{\rm noisy} - {\bm \Theta}^{\rm noisy} \bm \xi||_2^2  + \lambda_{1} ||\bm \xi||_1,
\end{equation}
where $\|\cdot\|_1$ is the $\ell_1$ norm and $\bm \xi^{\rm CS}$ is the recovered coefficient by the $\ell_1$ minimization. In this paper, we used LARS algorithm in \citep{efron2004least} for $\ell_1$ minimization. Then, we define a subset, denoted by $\mathcal{S}$, of $\mathcal{P} = [1,2,\ldots,P]$ based on the magnitude of the components of $\bm \xi^{\rm CS}$, i.e.,  
\begin{equation}\label{eq:l1norm}
    \mathcal{S} :=\Big\{\;i \in \mathcal{P}\; \big |\; |\xi_i^{\rm CS}| \ge \varepsilon,\;  \xi_i^{\rm CS} \in \bm \xi^{\rm CS} \Big\},
\end{equation}
where the threshold $\varepsilon > 0$ is chosen such that the reduced dictionary can be handled by the state-of-the-art MIO algorithm. We denote the reduced dictionary, the reduced indicator vector, and the reduced coefficient vector by 
\begin{equation}\label{eq:reduce_dic}
\begin{aligned}
   & \bm \Theta^{\rm noisy}_{\mathcal{S}}:=  
   \{ \theta_i(\mathbf X^{\rm noisy}) \in \bm \Theta^{\rm noisy}\; | \; i \in \mathcal{S}\},\\[2pt]
   & \bm \gamma_{\mathcal{S}} := \{ \gamma_i \in \bm \gamma\; |\; i \in \mathcal{S}\},\\[2pt]
   & \bm \xi_{\mathcal{S}} := \{ \xi_i \in \bm \xi, \; |\; i \in \mathcal{S}\},
\end{aligned}
\end{equation}
respectively. 
We emphasize that the recovered coefficients $\bm \xi^{\rm CS}$ in solving the $\ell_1$ minimization problem is not used to determine the final estimation of the coefficients. Instead, it is only used to help screening and narrowing down the range of candidate terms for high-dimensional problems with large $P$. 
% Hence, we can always have $k < S < P$ by setting a small value of $\lambda_1$.

\subsection{Mixed-integer optimization for determining the indicator $\mathbf{\gamma}_\mathcal{S}$}\label{sec:MIO}
%
% \begin{itemize}
%     \item give a specific reason why l2 regularization is needed.
%     \item the MIO solver can only take k sparsity or <=k sparsity?
%     \item 
% \end{itemize}
% %
We start from converting the problem in Eq.~\eqref{eq:sparse_reg_noise_k} into an MIO problem.  Then, the 
%
% For the $j$-th equation, we omit the subscript $j$ in $\bm \gamma_j$, $\bm \xi_j$ and $k_j$ for the simplicity of notation.
MIO problem constrained by a given sparsity can be written as
\begin{align}
\label{MIO0} \tag{$P_0$} \min_{\bm \xi_{\mathcal{S}}, \bm \gamma_{\mathcal{S}}}~ & \big\|\dot{\bm{x}}^{\rm noisy} - {\bm \Theta}^{\rm noisy}_{\mathcal{S}} (\bm \gamma_{\mathcal{S}} \circ \bm \xi_{\mathcal{S}})\big\|_2^2 + \lambda_{2} ||\bm \xi_{\mathcal{S}}||_2^2 \\[0pt]
% \label{c1}\textrm{s.t.}~  &|\bm \xi|\leq M \bm \gamma, & &\\
% \label{c2}  &\bm \gamma^T \bm e = k,\\
% \label{c1} \textrm{s.t.}~ & ||\bm \xi||_0 = k, & &\\[2pt]
% \label{c3} &\; \bm \xi \in \mathbb{R}^{P\times 1},
\label{c1}\textrm{s.t.}~ & \|\bm \xi_{\mathcal{S}}\|_{\infty} \leq B, & \\[3pt]
\label{c3}  &\;\bm \gamma_{\mathcal{S}}^T\bm e = k,
\end{align}
where $k$ denotes the sparsity of $\gamma_\mathcal{S}$ and $B$ is the upper bound of the coefficient $\bm \xi_{\mathcal{S}}$, and the $L^2$ regularization term $\lambda_{2} ||\bm \xi_{\mathcal{S}}||_2^2$ is commonly added to help alleviate the influence of the measurement noises on the MIO optimization. 

\begin{remark}[Using normalized data for MIO]
The scales of different components of the dynamical system could be
significantly different, which can affect the performance of the MIO solver in determining the optimal $\bm \gamma_{\mathcal{S}}$. To resolve this issue, 
we standardize the data $\bm \Theta^{\rm noisy}$ and $\dot{\bm{x}}^{\rm noisy}$, and use the standardized data in MIO. 
\end{remark}

\begin{remark}
We emphasize that splitting the coefficient of ${\bm \Theta}^{\rm noisy}_{\mathcal{S}}$ into $\bm \gamma_{\mathcal{S}}$ and $\bm \xi_\mathcal{S}$ is to indicate that the goal of solving the MIO is only to determine $\bm \gamma_{\mathcal{S}}$, i.e., identify the correct terms in the reduced dictionary ${\bm \Theta}^{\rm noisy}_{\mathcal{S}}$. Even though an MIO algorithm, e.g., the \texttt{CPLEX}, will also provide an estimate of $\bm \xi_{\mathcal{S}}$, we will not use the estimate as our final solution.  
\end{remark}

The goal of this subsection is to determine $\bm \gamma_{\mathcal{S}}$ by solving the MIO problem in Eq.~\eqref{MIO0}. However, there are two hyperparameters, i.e., the sparsity $k$ and the $\ell_2$-norm weight $\lambda_2$, that could significantly affect the outcome of the MIO solver. To address this issue, we perform a grid search with a cross validation metric to tune the two hyperparameters and obtain $\bm \gamma_{\mathcal{S}}$.

% Their values significantly affect the performance of CS-MIO for correct physical term identification. Hereby, we introduce effective strategy for the given dataset to select the best $k$ and $\lambda_2$, which are further used for accurate identification of physical terms.

% We tune to select the best parameter pair $(k, \lambda_2)$ by cross validation. %
We first define a tensor grid of $(k, \lambda_2)$. The grid for $k$ is easily defined as $\{1, 2, \ldots, k^{\max}\}$ based on the maximum allowable sparsity $k^{\max}$. The upper bound for $\lambda_2$, denoted by $\lambda_2^{\max}$, is defined by the norm
\[
    \lambda_2^{\max} := \|({\bm \Theta}^{\rm noisy})^{\top} {\dot{\bm x}}^{\rm noisy}\|_{\infty}.
\]
% It is a reasonable upper bound because {\bf needs a justification here}. 
%
% we set a maximum possible value $k_{max}$ for the selection of $k$.
% This is reasonable per our assumption that the forcing term $\mathbf{f}$ consists of only a few terms of $\bm \theta(\mathbf{x})$, regardless of the dimensionality $J$.
% For each $k=1,2,\dots,k_{max}$, we further create a set of $\lambda_2$ values following a data-driven manner.
% We examine all the candidate models associated with different parameter pairs $(k,\lambda_2)$ to select the best among them. 
% To do this, we find the maximum possible value $\lambda_2^{max}$ for $\lambda_2$ by computing the infinite norm defined by:
% \begin{equation}
%     \lambda_2^{max} := ||(\widetilde{\bm \Theta}^{\rm noisy})^T \widetilde{\dot{X}}^{\rm noisy}||_{\infty}.
% \end{equation}
This is followed by setting a small ratio $r$ of $\lambda^{\max}_{2}$ to set the minimum allowed value $\lambda^{\min}_{2}$, i.e., $\lambda^{\min}_{2}=r\lambda^{\max}_{2}$. 
Empirically, if $N>P$, we set $r=0.0001$; otherwise $r=0.01$. 
%\cite{glmnet}. 
Afterwards, we uniformly sample $m$ values from interval $[\log{(\lambda^{\min}_{2})}, \log{(\lambda^{\max}_{2})}]$ 
, where $m$ is practically set to 50 or 100. Then by taking exponentials of the sampled values, we obtain a set of $\ell_2$-norm weight, denoted by $\bm \Lambda = \{\lambda_2^{1},\dots,\lambda_2^{m}\}$.

We next perform the cross validation to choose the best hyperparameters from the tensor grid of $(k, \lambda_2)$.
% Thus, for each pair of $(k, \lambda_2)$ where $k\in \{1,\dots, k_{max}\}$ and $\lambda_2 \in \bm \Lambda$, we train the model by iteratively partitioning $\widetilde{\bm \Theta}^{\rm noisy}$ into training set and validation set following the procedure of cross validation. 
Specifically, we evenly partition the data set $\{\dot{\bm x}^{\rm noisy}, \bm \Theta^{\rm noisy}_{\mathcal{S}}\}$ with a total of $N$ measurements into 
%
%$\mathcal{I} =\{1,\dots,N\}$ into 
$T$ disjoint subsets, denoted by $\{\dot{\bm x}^{\rm noisy}_{1}, \bm \Theta^{\rm noisy}_{\mathcal{S},1}\}, \ldots,  \{\dot{\bm x}^{\rm noisy}_{T}, \bm \Theta^{\rm noisy}_{\mathcal{S},T}\}$, respectively. For one subset $\{\dot{\bm x}^{\rm noisy}_{t}, \bm \Theta^{\rm noisy}_{\mathcal{S},t}\}$ and one pair of $(k, \lambda_2)$, the error of the MIO solution is defined by
\begin{equation}\label{cv_err}
    e_t(k, \lambda_2)  := \big\|{\dot{\bm x}}^{\rm noisy}_{t} - {\bm \Theta}^{\rm noisy}_{\mathcal{S}, t}\,  {(\bm \gamma_{\mathcal{S},t}\circ \bm \xi_{\mathcal{S},t}})\big\|_2^2,
\end{equation}
where $\bm \gamma_{\mathcal{S},t}$ and $ \bm \xi_{\mathcal{S},t}$ are obtained by solving the MIO problem in Eq.~\eqref{MIO0} using the complementary data set $\{\dot{\bm x}^{\rm noisy} \backslash \dot{\bm x}^{\rm noisy}_{t}, \bm \Theta^{\rm noisy}_{\mathcal{S}}\backslash\bm \Theta^{\rm noisy}_{\mathcal{S},t}\}$. The errors for other subsets and choices of $k, \lambda_2$ can be obtained similarly. 
%
% $\mathcal{I}_1,\dots, \mathcal{I}_T$. We will assume $T$ divides $N$ for convenience. 
% Consider fold $\mathcal{I}_t$, let $\widehat{\bm \gamma}^{(-t)}$ and $ \widehat{\bm\xi}^{(-t)}$ be the coefficients obtained by solving Formulation \eqref{MIO2} with only those paired data not in fold $\mathcal{I}_t$, that is $(\widetilde{\bm \Theta}^{\rm noisy}_{\mathcal{I} \setminus \mathcal{I}_t}, \widetilde{\dot{X}}^{\rm noisy}_{\mathcal{I} \setminus \mathcal{I}_t})$. 
% Then, we compute sum of squared error (SSE) loss of fold $\mathcal{I}_j$ as
% \begin{equation}
%     e_t =  ||\widetilde{\dot{X}}^{\rm noisy}_{\mathcal{I}_t} - \widetilde{\bm \Theta}^{\rm noisy}_{\mathcal{I}_t}  \widehat{\bm \xi}^{(-t)}||_2^2.
% \end{equation}
% The error loss for other folds are defined analogously. 
Then the total error for the pair $(k, \lambda_2)$ is defined by
\begin{equation}\label{totalE}
    {\mathcal{E}}(k, \lambda_2) := \sum_{t=1}^{T} e_t(k, \lambda_2),
\end{equation}
% which represents the total loss for parameter pair $(k,\lambda_2)$.
and the best hyperparameters are obtained by 
\begin{equation}
    (k^*,\lambda_2^*) = \argmin_{(k,\lambda_2)} {\mathcal{E}}{(k, \lambda_2)}.
\end{equation}
The final step in this subsection is to solve the MIO problem with the best hyperparameters $(k^{*},\lambda_2^*)$ to obtain the optimized indicator vector, denoted by 
$\bm \gamma_\mathcal{S}^*$.

% \subsection{Least-squares approach for determining the coefficient $ \xi_{\mathcal{S}}$}\label{LS}

% The physical terms are finally identified by applying the selected best parameters $(k^*,\lambda_2^*)$ in solving Formulation \eqref{MIO2} with the whole standardized dataset $(\widetilde{\bm \Theta}^{\rm noisy}, \widetilde{\dot{X}}^{\rm noisy})$. Note here we introduce cross validation for given limited dataset. Other tuning procedures can also be applied for selecting $(k^*,\lambda_2^*)$. For example, with unlimited dataset from simulation, we can use additional dataset as validation dataset, as shown in Section \ref{sec:lorenz3}. 

%\subsection{Coefficient identification}
After the optimal $\bm \gamma_{\mathcal{S}}^*$ is determined using the MIO method, we use the standard least-squares approach to estimate the coefficient $\bm \xi_{\mathcal{S}}$. In this case, we use the original data, not the standarized data used in the MIO method, to solve the following least-squares problem
\begin{equation}\label{LS}
     \bm \xi^*_{\mathcal{S}} = \argmin_{\bm \xi_{\mathcal{S}}} \|\dot{\bm x}^{\rm noisy} - (\bm \Theta^{\rm noisy} \bm \gamma_{\mathcal{S}}^*)\; \bm \xi_{\mathcal{S}} \|_2^2,
\end{equation}
where the matrix $\bm \Theta^{\rm noisy} \bm \gamma_{\mathcal{S}}^*$ only contains the columns of $\bm \Theta^{\rm noisy}$ identified by $\bm \gamma_{\mathcal{S}}^*$.

\subsection{Summary of the CS-MIO algorithm}
We summarize the proposed CS-MIO method in Algorithm \ref{alg:cs_mio}. 
The CS-MIO algorithm is general by combining the capability of expected sparsity control with physical term selection and coefficient estimation.
The key of the CS-MIO algorithm is on solving the MIO formulation.
%Due to the apparent difficulties with the cardinality constraints of MIO, Problem \eqref{MIO2} is nonconvex and $\mathcal{NP}$-hard
%, and has attracted little attention in the past decades \cite{bertsimas2016best,bertsimas2020sparse}. 
%Until recent years, studies for solving this MIO formulation start to emerge \cite{bertsimas2016best,bertsimas2020sparse}.
In present study, we fully take advantages of state-of-the-art algorithm in modern optimization solver \texttt{CPLEX} for solving the MIO problems.
With appropriate settings for the time limit and optimality gap, the solver returns the optimal solution. Even if we terminate the algorithm early, it still provides a solution with suboptimality guaranteed. We will discuss the details of parameter settings for the optimization solver in the following experiment studies.

\RestyleAlgo{ruled}
\SetKwComment{Comment}{/* }{ */}

% \begin{algorithm}
% \setstretch{1.1}
%         \caption{The CS-MIO algorithm}\label{alg:cs_mio}
%         \begin{algorithmic}[1]
%             \REQUIRE Time series noisy data $\mathbf{X}^{\rm noisy}$, $\dot{\mathbf{X}}^{\rm noisy}$
%             \STATE Construct $\bm \Theta^{\rm noisy}$
%     by evaluating $\bm \theta(\mathbf{x})$ at each data point in $\mathbf{X}^{\rm noisy}$;
%             \STATE Standardize the columns of $\bm \Theta^{\rm noisy}$ and ${\dot{\mathbf{X}}}^{\rm noisy}$ to have zero means and unit variance;
%             %
%             \FOR{$j\in \mathcal{J}$}
%             \IF{$P>S$}
            
%             \ENDIF
%             \ENDFOR
%             \RETURN d
%         \end{algorithmic}
%     \end{algorithm}

\begin{algorithm}[h!]
\caption{The CS-MIO algorithm}\label{alg:cs_mio}
\setstretch{1.2}
\SetKwInOut{Input}{Input}
\SetKwInOut{Output}{Output}
\Input{The noisy data $\mathbf{X}^{\rm noisy}$, $\dot{\mathbf{X}}^{\rm noisy}$}
% \textbf{Initialization}:
%     $\mathbf{\Xi}^* \gets \textbf{0} \in \mathbb{R}^{P\times J}$, $\mathbf{\Gamma}^* \gets \textbf{0} \in \mathbb{R}^{P\times J}$\;
    Construct matrix $\bm \Theta^{\rm noisy}$
    by evaluating $\bm \theta(\mathbf{x})$ at the data points in $\mathbf{X}^{\rm noisy}$\; 
    Standardize the columns of $\bm \Theta^{\rm noisy}$ and ${\dot{\mathbf{X}}}^{\rm noisy}$ to have zero means and unit variance\;
\For{$j\in \mathcal{J}$}{
    % Let $\widetilde{\dot{X}}^{\rm noisy}$ and $\dot{X}^{\rm noisy}$ be the $j$-th column of $\widetilde{\mathbf{\dot{X}}}^{\rm noisy}$ and $\mathbf{\dot{X}}^{\rm noisy}$, respectively\;
    %Err$_{min} \gets \infty$\;
    
    \If(\tcc*[f]{Compressive sensing-base dictionary reduction}){$P > P_{\max}$}{
    %$\mathbf{\Theta(X}^{\rm noisy}_1) \gets $Apply $\ell_1$ regularization to select $T$ columns of $\mathbf{\Theta(X}^{\rm noisy}_1)$\;
    % $\mathcal{S} =\Big\{i: \xi_i^{\ell_1} \neq 0,\; i \in \mathcal{P}, \; \bm \xi^{\ell_1} = \argmin_{\bm \xi}  ||\widetilde{\dot{X}}^{\rm noisy}- \widetilde{\bm \Theta}^{\rm noisy} \bm \xi||_2^2 + \lambda_{1} ||\bm \xi||_1 \Big\}$\;
    Construct the reduced dictionary $\mathcal{S}$ in Eq.~\eqref{eq:reduce_dic} by solving Eq.~\eqref{eq:l1norm}
    % $\mathbf{\Theta}^{\rm noisy}, \widetilde{\mathbf{\Theta}}^{\rm noisy} \gets$ select columns indexed by $\mathcal{S}$ to form new $\mathbf{\Theta}^{\rm noisy}$ and $\widetilde{\mathbf{\Theta}}^{\rm noisy}$\; 
    }
 
    \For(\tcc*[f]{MIO for determining $\bm\gamma_{\mathcal{S}}$ in Eq.~\eqref{MIO0}} ){$k =1,2,\cdots, k^{\max}$}
    {  
        Construct a grid $\bm \Lambda$ of $\lambda_2$ in Eq.~\eqref{MIO0}\;
        % Compute $\lambda_2^{max} = ||(\widetilde{\bm \Theta}^{\rm noisy})^T \widetilde{\dot{X}}^{\rm noisy}||_{\infty}$; let $\lambda_2^{min} = r\lambda_2^{max}$\; 
        % $\bm \Lambda = \{\lambda_2^{[1]},\dots,\lambda_2^{[m]}\} \gets$ Uniformly sample $m$ values from $[\log (\lambda_2^{min}), \log( \lambda_2^{max})]$ and then take the exponentials to obtain $\lambda_2^{[1]},\dots,\lambda_2^{[m]}$\;
        
        Divide $\{\dot{\bm x}^{\rm noisy}, \bm \Theta^{\rm noisy}_{\mathcal{S}}\}$ into $T$ disjoint subsets\;
        
        \For{$\lambda_2 \in \bm \Lambda$}
        {
            % Let $\mathcal{I}=\{1,\dots,N\}$ be the row index set of the $N$ measurement data; evenly partition $\mathcal{I}$ into $T$ disjoint subsets (folds) $\mathcal{I}_1,\dots, \mathcal{I}_T$\;
            
            %(\tcc*[f]{Cross validation})
            \For{$t = 1, \ldots, T$}{
                
                Solve the MIO problem in Eq.~\eqref{MIO0} using $k$ and $\lambda_2$\;
                
                Compute the error $e_t(k, \lambda_2)$ in Eq.~\eqref{cv_err}\;
        
                % Solve $( \widehat{\bm \gamma}^{(-t)}, \widehat{\bm\xi}^{(-t)}) = \argmin_{(\bm \gamma, \bm \xi)\in \Delta} ||\widetilde{\dot{X}}^{\rm noisy}_{\mathcal{I} \setminus \mathcal{I}_t} - \widetilde{\bm \Theta}^{\rm noisy}_{\mathcal{I} \setminus \mathcal{I}_t} \bm \xi||_2^2 + \lambda_{2} ||\bm \xi||_2^2$ where $\Delta=\{(\bm \gamma, \bm \xi): \gamma_i = 0 \Rightarrow \xi_i = 0, i \in \mathcal{P}, ||\bm \xi||_{\infty} \leq M, \bm \gamma^T \bm e = k, \bm \gamma \in \mathbb{B}^{P\times 1}, \bm \xi \in \mathbb{R}^{P\times 1}\}$\;
        
                % Compute squared error loss of $\mathcal{I}_t$: $e_t =  ||\widetilde{\dot{X}}^{\rm noisy}_{\mathcal{I}_t} - \widetilde{\bm \Theta}^{\rm noisy}_{\mathcal{I}_t}  \widehat{\bm\xi}^{(-t)}||_2^2$\;
                
            }
            Compute total error ${\mathcal{E}}{(k, \lambda_2)}$ in Eq.~\eqref{totalE}\;
            % \If{${\rm Err}^{(k, \lambda_2)} < {\rm Err}_{min}$}{
            %         ${\rm Err}_{min} = {\rm Err}^{(k, \lambda_2)}$\;
            %         %$\bm \xi^* = \bm \xi_k$\;
            % }
        }
    }
    Find the best hyperparameters $(k^*,\lambda_2^*) = \argmin_{(k,\lambda_2)} {\mathcal{E}}{(k, \lambda_2)}$
    
    Identify the optimal indicator $\bm \gamma^*_{\mathcal{S}}$ by solving the problem \eqref{MIO0} using $(k^*,\lambda_2^*)$\;
    
    %  Identify the terms: $(\bm \gamma^*, \bm \xi^*) = \argmin_{(\bm \gamma, \bm \xi)\in \Delta} ||\widetilde{\dot{X}}^{\rm noisy} - \widetilde{\bm \Theta}^{\rm noisy} \bm \xi||_2^2 + \lambda_{2}^* ||\bm \xi||_2^2$
    %   where $\Delta=\{(\bm \gamma, \bm \xi): \gamma_i = 0 \Rightarrow \xi_i = 0, i \in \mathcal{P}, ||\bm \xi||_{\infty} \leq M, \bm \gamma^T \bm e = k^*, \bm \gamma \in \mathbb{B}^{P\times 1}, \bm \xi \in \mathbb{R}^{P\times 1}\}$\Comment*[r]{Physical term identification}
       
      %Set $\bm \gamma^* = \widehat{\bm\gamma}^*$ 

    Determine the optimal coefficient $\bm \xi_{\mathcal{S}}^{*}$ by solving the problem in Eq.~\eqref{LS}\;
    % Identify the coefficients with original data: $\bm \xi^* = \argmin_{\bm \xi \in \Delta} ||\dot{X}^{\rm noisy} - \bm \Theta^{\rm noisy} \bm \xi ||_2^2~$where$~\Delta =\{\bm \xi: \gamma_i = 0 \Rightarrow \xi_i = 0, i\in\mathcal{P}, ||\bm \xi||_{\infty} \leq M, \bm \gamma = \bm \gamma^*, \bm \xi \in \mathbb{R}^{P\times 1}\}$\Comment*[r]{Coefficient identification}
    
    Set $\bm \gamma^*_{\mathcal{S}}$ and $\bm \xi^*_{\mathcal{S}}$ as the $j$-th column of $\bm \Gamma^*$ and $\bm \Xi^*$, respectively\;
}
\textbf{Return} $\dot{\mathbf{x}} = \bm \theta (\mathbf{x})(\bm \Gamma^* \circ \bm \Xi^*)$
\end{algorithm}

\section{Numerical experiments}\label{sec:ex}
We demonstrate the effectiveness of the proposed CS-MIO method for recovery of governing equations from large noise data.
We use several classical dynamical systems in \cite{brunton2016discovering} as the testing problems, including chaotic Lorenz 3 system, vortex shedding after a cylinder, bifurcation dynamical systems like Hopf normal form and logistic map. In addition, we also study the high-dimensional Lorenz 96 system. We compare CS-MIO with state-of-the-art method SINDy, specifically, the Python version solver PySINDy \cite{desilva2020,kaptanoglu2021pysindy}. 
For all the example systems, the experiments are deployed on a mobile workstation with Intel(R) Xeon(R) W-10885M CPU @ 2.40GHz, 128 GB memory, 64 bit Windows 10 Pro operating system for workstations. 

\begin{remark}[Reproducibility]
The algorithm of CS-MIO is implemented in Python. The code is publicly available at \url{https://github.com/utk-ideas-lab/CS-MIO}. All the numerical results presented in this section can be exactly reproduced using the code on Github. 
\end{remark}

\subsection{Experimental settings}
We first give the experimental settings throughout the case studies. To better measure the noise level and the anti-noise capability of the method, we consider the signal-to-noise ratio (SNR).
In this work, we consider the averaged SNR of the dynamical system consisting of a set of $J$ governing equations,
\begin{equation}\label{SNR}
    \textrm{SNR}:=\frac{1}{J}\sum_{j=1}^J \frac{\textrm{Var}(S_j)}{\textrm{Var}(\mathcal{Z}_j)},
\end{equation}
where $S_j \in \{X_j, \dot{X}_j\}$ is the signal data, i.e., the $j$-th column of matrices $\mathbf{X}$ or $\dot{\mathbf{X}}$, and $\mathcal{Z}_j \in \{\mathcal{U}_j,\mathcal{V}_j\}$ are the additive Gaussian noise, i.e., the $j$-th column of matrices $\mathbfcal{U}$ or $\mathbfcal{V}$ in Eq.~\eqref{eq:noisy_data}. The 
SNR gives a good indicator to assess the ability of methods to withstand noise in the data. Smaller SNR indicates a system with larger noise. We examine the anti-noise capability of the methods over a wide range of SNRs for the studied examples.
%Per the practice from statistical community, an SNR of 0.25 is typical for noisy experimental data, and an SNR around 0.02 is considered as huge noise in financial returns data \cite{hastie2020best}. 
%Under a wide variety of examples and SNRs, our results show the proposed CS-MIO method can recover exact governing equations with the SNR as low as 0.05, exhibiting the strong capability of denoising. 
In the following cases, we impose the below two types of noise by considering the signal that can be measured.

\begin{itemize}
    \item \emph{Type 1 Noise}: Both the state variables $\mathbf{x}$ and time derivatives $\dot{\mathbf{x}}$ can be measured; Gaussian noise is added to $\dot{\mathbf{x}}$.
    \item \emph{Type 2 Noise}: Only state variables $\mathbf{x}$ can be measured. Gaussian noise is added to $\mathbf{x}$. The time derivatives $\dot{\mathbf{x}}$ are computed by total variation derivative (TVD) \cite{chartrand2011numerical}.
\end{itemize}

%Other experimental settings are detailed in \emph{SI Appendix, Section 4.B}.
%The setting of $k$ is based on our assumption that there are only a few terms in the governing equations, and also domain-educated researchers and practitioners could have an intuition for the ground truth $k$.
%This is reasonable since our aim is focused on the exact recovery capability of CS-MIO for a given $k$. Notably, we could always enhance the choice of $k$ through cross validation for model selection as described in Algorithm . We will discuss this for the Lorenz 3 system in Sec. \ref{sec:lorenz3}.

% The settings for CS-MIO are presented as below.
We use state-of-the-art algorithm in modern optimization solver \texttt{CPLEX} (Python package \texttt{docplex}) for solving the MIO problems.
Unless specifically mentioned, we use up to fifth order total-degree polynomials throughout the examples to define the initial dictionary. 
The choice of the upper bound $B$ in Eq.~\eqref{c1} impacts the strength of the MIO formulation, especially when looking for good lower bounds.
$B\in\mathbb{R}$ is a sufficiently large constant such that $B \geq \|\boldsymbol{\xi}^*\|_{\infty}$. This setting is, however, not applicable because the $\boldsymbol{\xi}^*$ is not known a prior. 
%In general, a small value of $M$ could possibly result in a solution that is different with the optima of the primal best subset selection problem with the same $k$. If on the contrary, $M$ is too large, it could lead to loose lower bound of $\boldsymbol{\xi}^*$ and thus deteriorate the computing efficiency.
Some methods have been studied to set $B$ values by finding the upper bound of $\boldsymbol{\xi}^*$ using data-driven manners such as cumulative coherence function and solving convex optimization methods \cite{bertsimas2016best}. In this paper, we use a loosing upper bound $B=1000$ for all the examples. Besides, we set the \texttt{timelimit} to be 600 seconds, and the \texttt{mipgap} to be 0 for the invoked branch-and-cut algorithm in \texttt{docplex}. This refers to that if the branch-and-cut finds a solution within 600 seconds, it will be the optimal solution with zero gap; otherwise, the provided solution will be suboptimal and its gap to the lower bound, and thus to the optima, will be clearly quantified.

%We next discuss the setting of $\ell_2$-norm regularizations.
%the $\lambda_{1}$ is to induce a very loose pre-selection of screening of significant terms when dealing with high-dimensional problems. 
% In present study, we use the in \texttt{LassoLars} method from Python package \texttt{scikit-learn} with $\lambda_{1,i}=$1e-6, a very small regularization weight but capable of narrowing down thousands of terms to hundreds. Other settings remain as default. In addition, we set $S=100$.
% %the maximum number of terms kept in $\boldsymbol{\theta}_s(\mathbf{x})$. 
% It is seen in most times, the subset of terms resulting from LARS has a larger size than $S$. In this case we order the nonzero terms in decreasing order of the absolute values of their coefficients and select the top $S$ terms to form $\mathcal{S}$.

{\em The metrics for performance comparison.} We evaluate the performance of the identification of $\bm \Gamma$ in Eq.~\eqref{eq:sparse_reg_noise_k} by {\em the number of exactly recovered equations of the target dynamical system}, defined by
\begin{equation}\label{metric}
    A (\bm \Gamma):=\sum_{j=1}^J \mathbf{1}_{\boldsymbol{\gamma}_j = \boldsymbol{\gamma}_j^{\dagger}}, \;\text{ with }\; \mathbf{1}_{\boldsymbol{\gamma}_j = \boldsymbol{\gamma}_j^{\dagger}}
    = \left\{\begin{aligned}
       1, & \text{ if } \boldsymbol{\gamma}_j := \boldsymbol{\gamma}_j^{\dagger},\\
       0, & \text{ if } \boldsymbol{\gamma}_j \neq \boldsymbol{\gamma}_j^{\dagger},\\
    \end{aligned}\right.
\end{equation}
where $\boldsymbol{\gamma}_j^{\dagger}$ is the ground truth and  $\boldsymbol{\gamma}_j$ is recovered by a method.
% The equal supports indicate $\boldsymbol{\gamma}_j$ tells the whole truth and nothing but the truth, thereby is an exact recovery of the $j$-th equation. 
% Function $1(\cdot)$ equals to 1 if the equality condition holds, otherwise 0.
The exact recovery for the entire dynamical system occurs when $A(\bm \Gamma)=J$. 
When the exact $\bm \Gamma$ can be recovered, we evaluate the accuracy of the approximation of the coefficients in $\bm \Xi$ in Eq.~\eqref{eq:sparse_reg_noise_k} by the differences between the approximate and the exact coefficients and trajectories.

\subsection{The chaotic Lorenz 3 system}\label{sec:lorenz3}

Consider the 3-dimensional chaotic Lorenz system governed by the following equations: 
\begin{align}
     \dot{x} & = \alpha (y - x), \\
     \dot{y} & = x (\rho - z) - y, \\
     \dot{z} & = x y - \beta z.
\end{align}
With $\sigma = 10$, $\beta = 8/3$ and $\rho = 28$, the Lorenz 3 system performs chaotically. We generate the data using initial point $(x,y,z)=(-8,8,27)$ with time step $\Delta t=0.001$ in $t\in [0,60]$.
A set of noise standard deviation $\sigma$ is used to better quantify the spectrum of anti-noise capability of the methods. 
In particular, under Type 1 noise, Gaussian noise is added to $\mathbf{\dot{x}}$ with $\sigma$ ranging from 1 to 3000. SNR is computed by the added noise and $\mathbf{\dot{x}}$.
When under Type 2 noise, the Gaussian noise is added to $\mathbf{x}$ with $\sigma$ ranging from 0.01 to 20,
and the SNR is computed by the added noise and $\mathbf{x}$.
In this case, $\mathbf{\dot{x}}$ is smoothed using total variation derivative (TVD) of \cite{chartrand2011numerical}. 
The comparison results of CS-MIO and PySINDy for both cases are presented in Tables \ref{tab:lorenz3_gaussian} and \ref{tab:lorenz3_tvd}, respectively.

\begin{table}[h!]
\footnotesize
\caption{Comparison of the number of exactly recovered equations, i.e., the metric ${A}(\mathbf{\Gamma})$ in Eq.~\eqref{metric}, for the Lorenz 3 system. Compared with PySINDy, our CS-MIO method can correctly recover all the three equations, i.e., identifying the correcting $\mathbf{\Gamma}$ in Eq.~\eqref{eq:ode1}, under smaller SNR values.}
\begin{subtable}{0.43\textwidth}
  \centering
  \caption{Results under Type 1 noise.}
  \label{tab:lorenz3_gaussian}
  \begin{tabular}{cc|cc}
      \toprule
             &	 & \multicolumn{2}{c}{The metric ${A}(\mathbf{\Gamma})$} \\[2pt]
  Noise std  &	SNR & PySINDy  &   CS-MIO	\\
   \midrule
%0.01	&	41914317.1292	&	3	&	3	\\
%0.1	&	419143.1713	&	3	&	3	\\
1	&	4191.621	&	\textbf{3}	&	\textbf{3}	\\
10	&	41.921	&	2	&	\textbf{3}	\\
50	&	1.677	&	1	&	\textbf{3}	\\
100	&	0.419	&	0	&	\textbf{3}	\\
%150	&	0.186	&	0	&	\textbf{3}	\\
200	&	0.105	&	0	&	\textbf{3}	\\
%250	&	0.067	&	0	&	\textbf{3}	\\
300	&	0.047	&	0	&	\textbf{3}	\\
500	&	0.017	&	0	&	2	\\
1000	&	0.004	&	0	&	1	\\
3000	&	0.001	&	0	&	0	\\
\bottomrule
\end{tabular}
\end{subtable}
\hspace{0.5cm}
\begin{subtable}{0.5\textwidth}
   \centering
  \caption{Results under Type 2 noise.}
  \label{tab:lorenz3_tvd}
  \begin{tabular}{cc|cc}
      \toprule
&	 & \multicolumn{2}{c}{The metric ${A}(\mathbf{\Gamma})$} \\[2pt]
  Noise std  &	SNR & PySINDy  &   CS-MIO	\\
   \midrule
0.01	&	729427.159	&	\textbf{3}	&	\textbf{3}	\\
0.05	&	29178.677	&	2	&	\textbf{3}	\\
0.1	&	7295.616	&	2	&	\textbf{3}	\\
0.5	&	292.848	&	1	&	\textbf{3}	\\
1	&	73.981	&	0	&	\textbf{3}	\\
2	&	19.255	&	0	&	\textbf{3}	\\
5	&	3.926	&	0	&	2	\\
10	&	1.733	&	0	&	1	\\
20	&	1.184	&	0	&	0	\\
\bottomrule
\end{tabular}
\end{subtable}
\end{table}

% \begin{table}
% \footnotesize
%   \centering
%   \caption{Number of exactly recovered equations for Lorenz 3: only $\mathbf{x}$ can be measured; noise is added to $\mathbf{x}$; $\mathbf{\dot{x}}$ is computed by TVD \cite{chartrand2011numerical}.}
%   \label{tab:lorenz3_tvd}
%   \begin{tabular}{rr|rr}
%       \toprule
%   Noise: $\sigma$  &	SNR &   PySINDy &   CS-MIO	\\
%   \midrule
% 0.01	&	729427.159	&	\textbf{3}	&	\textbf{3}	\\
% 0.05	&	29178.677	&	2	&	\textbf{3}	\\
% 0.1	&	7295.616	&	2	&	\textbf{3}	\\
% 0.5	&	292.848	&	1	&	\textbf{3}	\\
% 1	&	73.981	&	0	&	\textbf{3}	\\
% 2	&	19.255	&	0	&	\textbf{3}	\\
% 5	&	3.926	&	0	&	2	\\
% 10	&	1.733	&	0	&	1	\\
% 20	&	1.184	&	0	&	0	\\
% \bottomrule
% \end{tabular}
% \end{table}

Table \ref{tab:lorenz3_gaussian} and \ref{tab:lorenz3_tvd} show that CS-MIO significantly outperforms PySINDy in terms of the number of exactly recovered equations. 
Under Type 1 noise as shown in Table \ref{tab:lorenz3_gaussian}, CS-MIO is able to exactly recover the differential equations with SNR as low as 0.047. 
Comparing to the SNR value of 4191.621 by PySINDy, this results in a tremendous difference of almost 100,000 times. 
Similar conclusions can be made under Type 2 noise as shown in Table \ref{tab:lorenz3_tvd}. 
It is noted in this case the white noise added to $\mathbf{x}$ is no longer Gaussian after using numerical differentiation and is difficult to handle. Thus, the performance of both CS-MIO and PySINDy is downgraded at smaller SNRs.

\begin{table}[h!]
\footnotesize
\caption{\footnotesize Comparison of discovered equations by PySINDy and CS-MIO under (a) Type 1 noise at $\sigma$=300; and (b) Type 2 noise at $\sigma$=2. The CS-MIO method correctly identified all the terms in the Lorenz 3 system, but
PySINDy picked up incorrect terms.}
\begin{subtable}[c]{0.6\textwidth} 
  \centering
  %\captionof{figure}{\footnotesize Number of exactly recovered equations.}
  %\resizebox{\columnwidth}{!}{
  \renewcommand{\arraystretch}{1}
  \begin{tabular}{p{1.2cm}|p{0.12cm}p{6.4cm}}
    \toprule
    \multirow{3}{1.2cm}{Ground Truth}  &   $\dot{x}=$    & $-10x+10y$ \\
                    &   $\dot{y} =$    & $28x-y-xz$\\
                    &   $\dot{z} =$    & $-\frac{8}{3}z +xy$\\
    \midrule
    \multirow{3}{1.2cm}{PySINDy}       &   $\dot{x} =$   & $-6.62 - 13.62x + 11.80y + 0.38z - 0.14x^2 + 0.23xy + 0.11xz - 0.11y^2 - 0.05yz$\\
                    &   $\dot{y} =$    & $3.46 + 29.32x - 1.32y - 0.05z - 1.03xz$\\
                    &   $\dot{z} =$    &   $-7.14 - 0.13x + 0.15y - 2.28z - 0.08x^2 + 1.05xy$\\
    \midrule
    \multirow{3}{1.2cm}{CS-MIO}        &   $\dot{x} =$    & $- 9.72x + 9.70y$\\
                    &   $\dot{y} =$    & $29.28x  - 1.31y - 1.03xz$\\
                    &   $\dot{z} =$    & $- 2.64z + 1.00xy$\\
    \bottomrule
    \end{tabular}
    %}
        \caption{Results under Type 1 noise with $\sigma$=300.}  
    \label{tab:lorenz3_Gaussian_eqs}
\end{subtable}
\begin{subtable}[c]{0.34\textwidth} 
  \centering
  %\captionof{figure}{\footnotesize Number of exactly recovered equations.}
  %\resizebox{\columnwidth}{!}{
  \renewcommand{\arraystretch}{1}
  \begin{tabular}{p{1.2cm}|p{0.12cm}p{3.2cm}}
    \toprule
    \multirow{3}{1.2cm}{Ground Truth}  &   $\dot{x}=$    & $-10x+10y$ \\
                    &   $\dot{y} =$    & $28x-y-xz$\\
                    &   $\dot{z} =$    & $-\frac{8}{3}z +xy$\\
    \midrule
    \vspace{0.2cm}
    \multirow{3}{1.2cm}{PySINDy}       &   $\dot{x} =$   & $-0.21 - 9.87x + 9.89y$\\
    \vspace{0.1cm}
                    &   $\dot{y} =$    & $0.10 + 27.23x - 0.73y$\\
    \vspace{0.1cm}
                    &   $\dot{z} =$    &   $-1.05 - 2.62z + 1.00xy$\\
    \midrule
    \multirow{3}{1.2cm}{CS-MIO}        &   $\dot{x} =$    & $- 9.87x + 9.89y$\\
                    &   $\dot{y} =$    & $27.23x  - 0.73y - 0.98xz$\\
                    &   $\dot{z} =$    & $- 2.66z + 1.00xy$\\
    \bottomrule
    \end{tabular}
    %}
    \caption{Results under Type 2 noise with $\sigma$=2.} 
    \label{tab:lorenz3_TVD_eqs}
\end{subtable}
\label{tab:lorenz3_300_example}
\end{table}

\begin{figure}[h!]
    \begin{subfigure}[c]{0.24\textwidth}
    \centering
    \includegraphics[width=1.52in]{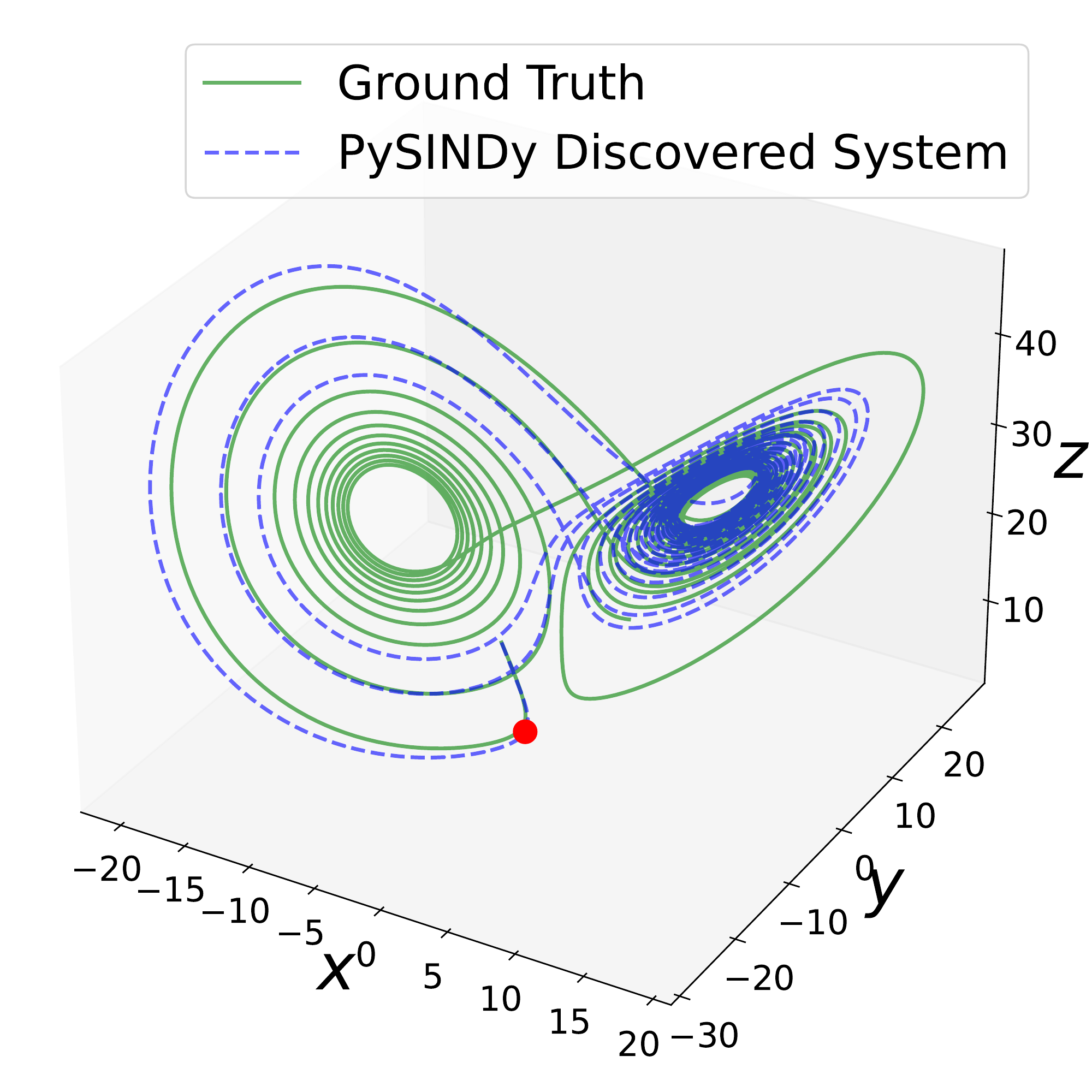}    
    \caption{}\label{fig:Lorenz3_sindy_gaussian_300_traj}
\end{subfigure}
\begin{subfigure}[c]{0.24\textwidth}
    \centering
    \includegraphics[width=1.52in]{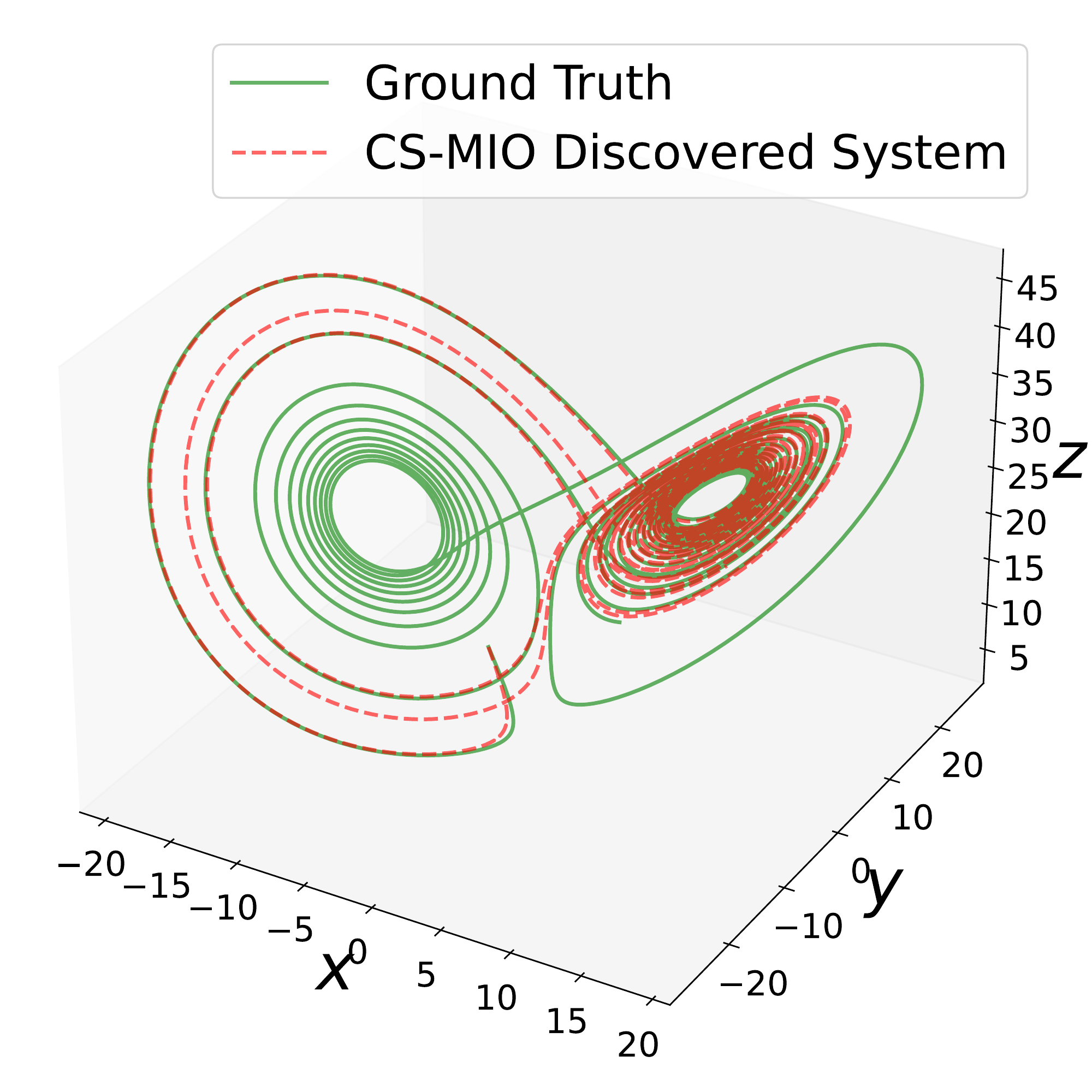}    
    \caption{}\label{fig:Lorenz3_cs_mio_gaussian_300_traj}
\end{subfigure}
\begin{subfigure}[c]{0.24\textwidth}
    \centering
    \includegraphics[width=1.52in]{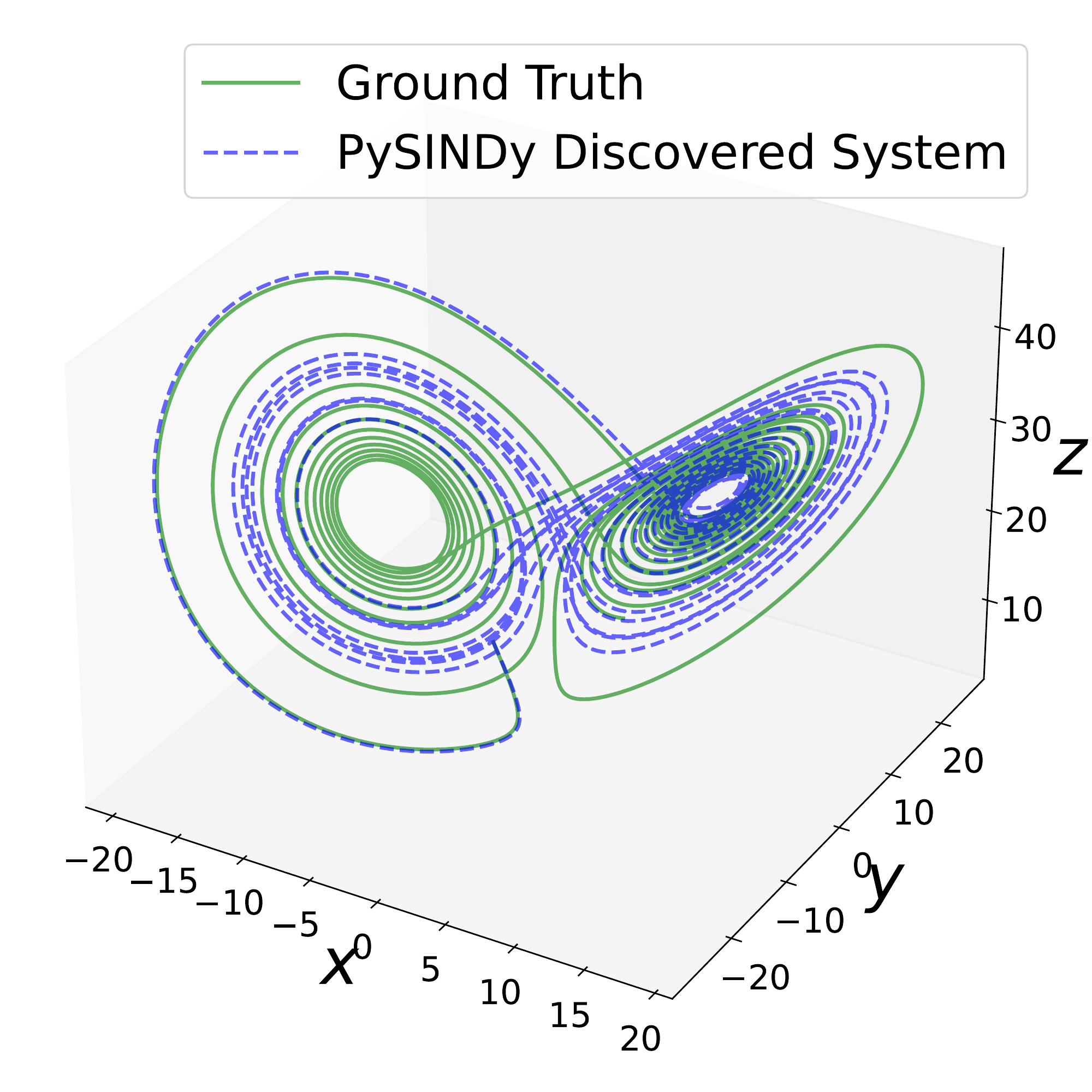}    
    \caption{}\label{fig:Lorenz3_sindy_tvd_300_traj}
\end{subfigure}
\begin{subfigure}[c]{0.24\textwidth}
    \centering
    \includegraphics[width=1.52in]{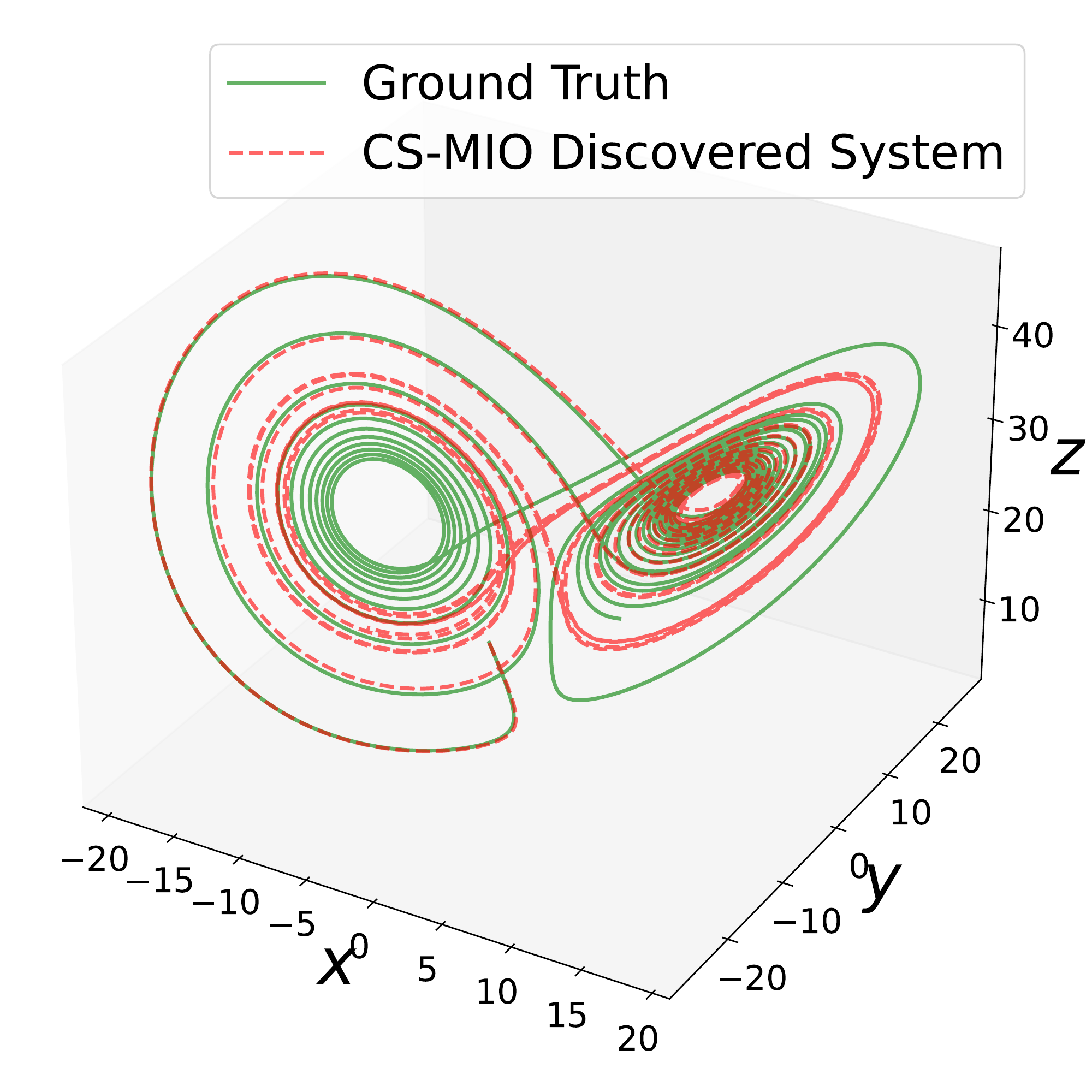}    
    \caption{}\label{fig:Lorenz3_cs_mio_tvd_300_traj}
\end{subfigure}
    \caption{(a) and (b) show the trajectories of PySINDy and CS-MIO discovered equations in Table \ref{tab:lorenz3_Gaussian_eqs}, respectively. The trajectory of PySINDy identified system deviates from the ground truth right at the beginning (the red dot) of the trajectory. (c) and (d) show the trajectories of PySINDy and CS-MIO identified equations in Table \ref{tab:lorenz3_TVD_eqs}, respectively.}
    \label{fig:lorenz3_traj_comparison}
\end{figure}

Tables \ref{tab:lorenz3_Gaussian_eqs} and \ref{tab:lorenz3_TVD_eqs} show the discovered equations by CS-MIO and PySINDy under Type 1 noise (at noise magnitude 300) and Type 2 noise (at noise magnitude 2), respectively.
Obviously from these tables, PySINDy includes redundant false terms.
On the contrary, CS-MIO identifies all and only the ground truth terms, while remains small deviation of the identified parameters from the ground truth.
This can be seen from the trajectories of the discovered equations by both PySINDy and CS-MIO in Figure \ref{fig:lorenz3_traj_comparison}. Figure \ref{fig:Lorenz3_sindy_gaussian_300_traj} shows the trajectory of PySINDy identified system under Type 1 noise at $\sigma$=300. It is seen the trajectory starts to deviate the ground truth right at the beginning, shown by the red dot. In contrast, the trajectory of CS-MIO identified system can coincide for longer time well with the ground truth, as shown in Figure \ref{fig:Lorenz3_cs_mio_gaussian_300_traj}. Appendix \ref{sec:appendix_lorenz3} gives more details of identified models of Lorenz 3 system by CS-MIO.

\subsection{The chaotic Lorenz 96 system}\label{ref:lorenz96}

%High-dimensional dynamic systems are commonly seen especially in fluid dynamics with data on a spatial-temporal resolved snapshots, resulting in up to billions of variables to represent the system dynamics. Traditional dimensional reduction techniques, such as proper orthogonal decomposition \cite{holmes2012turbulence} or dynamic mode decomposition \cite{taira2017modal}, are often used to extract dominant features (modes) from the high dimensional space for construcing reduced-order models \cite{loiseau2018sparse}. Instead of using coordinate transformation, the CS-MIO method considers a collection of regularization techniques for dimension reduction and dynamics identification. In the algorithm, CS-MIO begins its application with the $\ell_1$ conservative screening and then proceeds the $\ell_2$ coefficients shrinkage for withstanding the noise, and terms selection by the $\ell_0$ constraint (\emph{SI Section 3.B}).

We demonstrate the effectiveness of the proposed CS-MIO method for high dimensional problems using Lorenz 96 dynamic system, which is defined as follows. For $j=1,\cdots,J$,
\begin{equation}
    \dot{x}_j = (x_{j+1}-x_{j-2})x_{j-1} - x_j + F,
\end{equation}
where $x_{j}$ is the state variable and $F$ is a forcing constant. Here it is assumed that $x_{-1}=x_{J-1}$, $x_{0}=x_{J}$, $x_{J+1}=x_{1}$ and $J\geq 4$. In this study, we set $J=96$. $F=8$ is a common value known to cause chaotic behavior. We use the initial condition $\mathbf{x}(0)=\boldsymbol{1}$ with a small perturbation 0.01 added to $x_1(0)$ to generate the dataset with time step $\Delta t=0.01$ in $t\in[0,600]$. 
We use second order polynomials in CS-MIO for the Lorenz 96 system with 96 variables, which results in 4752 polynomial terms. This leads to huge difficulties to deal with the high dimension.

For this high-dimensional Lorenz 96 system, we use compressive sensing approach as described in Eq.~\eqref{eq:l1norm} for pre-selecting a subset of at most $S=100$ significant terms from 4752 candidate terms. We use the \texttt{LassoLars} algorithm from Python package \texttt{scikit-learn}. We set $\lambda_{1}=10^{-6}$, a very small regularization weight but capable of narrowing down thousands of terms to hundreds. Other settings remain as default. 
It is seen in most times, the subset of terms resulting from  has a larger size than $S$. In this case, we order the nonzero terms in decreasing order of the absolute values of their coefficients and select the top $S$ terms to form the preselected subset $\mathcal{S}$.

\begin{table}[h!]
\footnotesize
\caption{Comparison of the number of exactly recovered equations, i.e., the metric ${A}(\mathbf{\Gamma})$ in Eq.~\eqref{metric}, for the Lorenz 96 system. Compared with PySINDy, our CS-MIO method can correctly recover all the 96 equations, i.e., identifying the correcting $\mathbf{\Gamma}$ in Eq.~\eqref{eq:ode1}, under smaller SNR values.}
\begin{subtable}{0.43\textwidth}
  \centering
  \caption{Results under Type 1 noise.}
  \label{tab:lorenz96_Gaussian}
  \begin{tabular}{cc|cc}
      \toprule
             &	 & \multicolumn{2}{c}{The metric ${A}(\mathbf{\Gamma})$} \\[2pt]
  Noise std  &	SNR & PySINDy  &   CS-MIO	\\
   \midrule
%0.01	&	3521465.902	&	96	&	96	\\
%0.1	&	35214.659	&	96	&	96	\\
1	&	352.147	&	\textbf{96}	&	\textbf{96}	\\
10	&	3.521	&	93	&	\textbf{96}	\\
20	&	0.880	&	20	&	\textbf{96}	\\
30	&	0.391	&	0	&	\textbf{96}	\\
40	&	0.220	&	0	&	\textbf{96}	\\
50	&	0.141	&	0	&	\textbf{96}	\\
70	&	0.072	&	0	&	88	\\
%90	&	0.043	&	0	&	69	\\
%110	&	0.029	&	0	&	35	\\
%130	&	0.021	&	0	&	21	\\
150	&	0.016	&	0	&	9	\\
%170	&	0.012	&	0	&	5	\\
%190	&	0.010	&	0	&	1	\\
%210	&	0.008	&	0	&	0	\\
230	&	0.007	&	0	&	0	\\
\bottomrule
\end{tabular}
\end{subtable}
\hspace{0.5cm}
\begin{subtable}{0.48\textwidth}
    \centering
  \caption{Results under Type 2 noise.}
  \label{tab:lorenz96_tvd}
  \begin{tabular}{cc|cc}
      \toprule
             &	 & \multicolumn{2}{c}{The metric ${A}(\mathbf{\Gamma})$} \\[2pt]
  Noise std  &	SNR & PySINDy  &   CS-MIO	\\
   \midrule
0.01	&	132501.470	&	\textbf{96}	&	\textbf{96}	\\
0.05	&	5300.059	&	95	&	\textbf{96}	\\
0.1	&	1325.015	&	91	&	\textbf{96}	\\
0.2	&	331.254	&	45	&	\textbf{96}	\\
0.4	&	82.813	&	14	&	\textbf{96}	\\
0.6	&	36.806	&	0	&	\textbf{96}	\\
0.8	&	20.703	&	0	&	\textbf{96}	\\
1	&	13.250	&	0	&	92	\\
10	&	0.133	&	0	&	0	\\
\bottomrule
\end{tabular}
\end{subtable}
\end{table}

% \begin{figure}
% \centering
% \includegraphics[width=3.5in]{Lorenz96_l2error.pdf}
% \caption{$l_2$ error vs time of the trajectories of the recovered Lorenz 96 system (red dashed) comparing to the ground truth (blue solid) from $t=0$ to $t=20$ under three noise magnitudes: 0.01, 1 and 50. The exact recovery fails when $\sigma$ ls larger than 300.}
% \end{figure}

% \begin{table}
% \footnotesize
%   \centering
%   \caption{Number of exactly recovered equations for Lorenz 96: only $\mathbf{x}$ can be measured; noise is added to $\mathbf{x}$; $\mathbf{\dot{x}}$ is computed by TVD \cite{chartrand2011numerical}.}
%   \label{tab:lorenz96_tvd}
%   \begin{tabular}{rr|rr}
%       \toprule
%   Noise: $\sigma$  &	SNR &   PySINDy &   CS-MIO	\\
%   \midrule
% 0.01	&	132400.400	&	\textbf{96}	&	\textbf{96}	\\
% 0.05	&	5296.016	&	95	&	\textbf{96}	\\
% 0.1	&	1324.004	&	91	&	\textbf{96}	\\
% 0.2	&	331.001	&	45	&	\textbf{96}	\\
% 0.4	&	82.750	&	14	&	\textbf{96}	\\
% 0.6	&	36.778	&	0	&	\textbf{96}	\\
% 0.8	&	20.688	&	0	&	\textbf{96}	\\
% 1	&	13.240	&	0	&	92	\\
% 10	&	0.132	&	0	&	0	\\
% \bottomrule
% \end{tabular}
% \end{table}

\begin{figure}
\begin{subfigure}[c]{0.45\textwidth}
    \centering
    \includegraphics[width=2.9in]{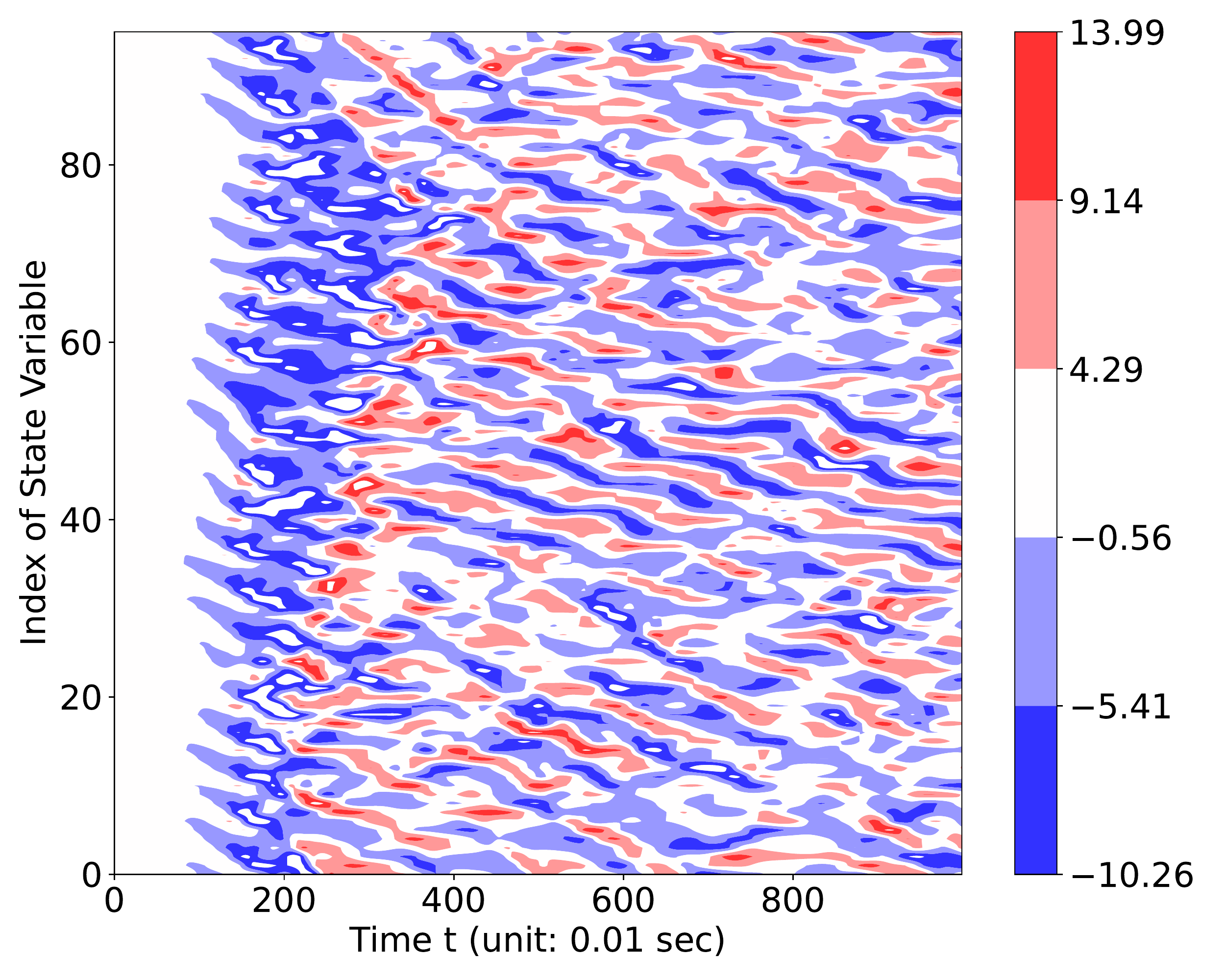}
    \caption{Hoverm{\"o}ller plot under Type 1 noise with $\sigma$=50.}
    \label{fig:lorenz96_Gaussian50_traj}
\end{subfigure}
\hspace{0.8cm}
\begin{subfigure}[c]{0.45\textwidth}
 \centering
    \includegraphics[width=2.9in]{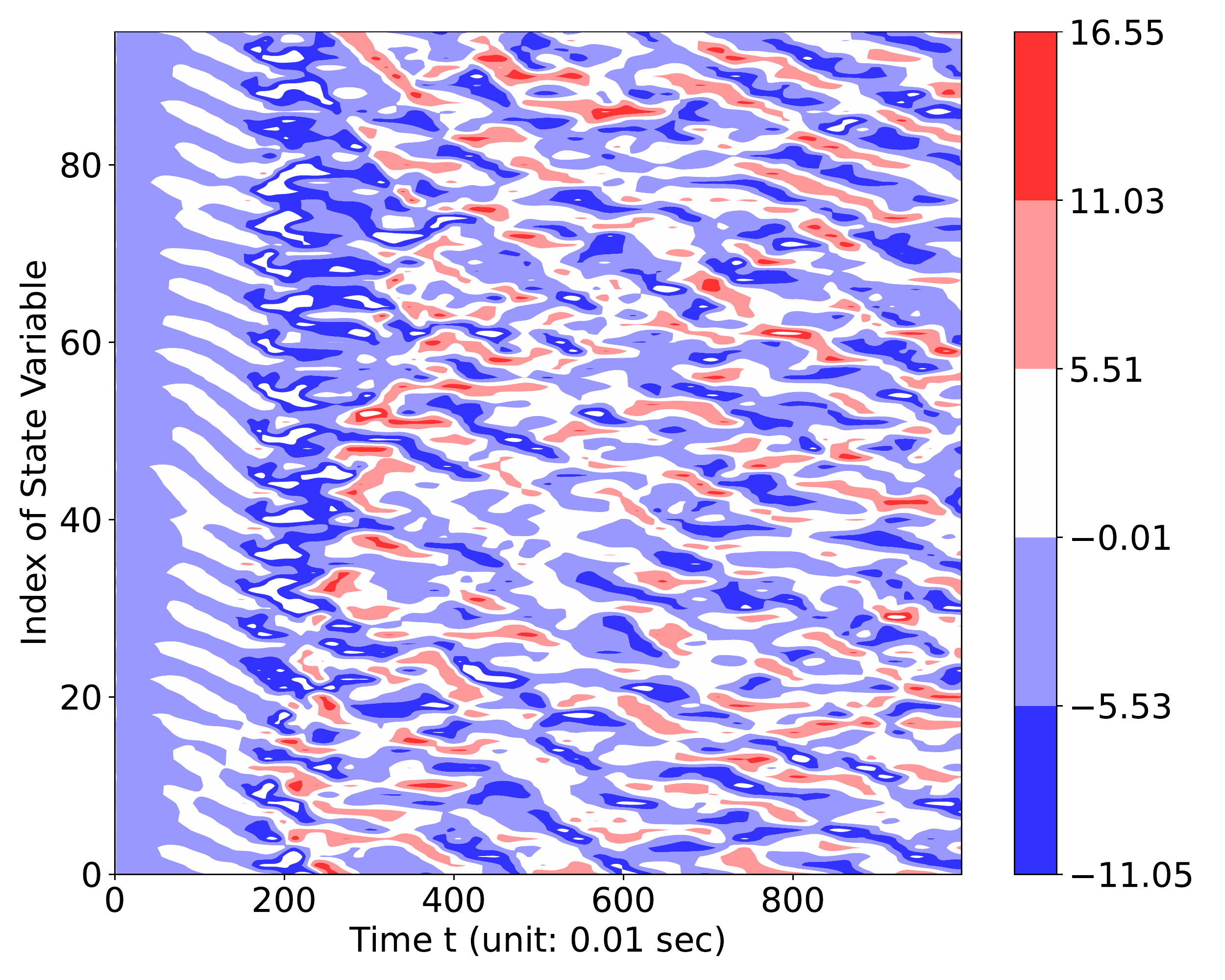}
    \caption{Hoverm{\"o}ller plot under Type 2 noise with $\sigma$=0.8.}
    \label{fig:lorenz96_TVD0.8_traj}
\end{subfigure} 
\caption{Hoverm{\"o}ller plot for difference between the identified system and ground truth of Lorenz 96 system in $t \in[0,10]$. The vertical axis is the index $j$ of the state variable. The values of the colors refer to the difference between the ground truth states $x_j(t)$ and the evolved states $\hat{x}_j(t)$ using the identified equations by CS-MIO. The more white-colored areas indicate the trajectory of the identified system agrees better with the ground truth.}
    \label{fig:lorenz96_50_trajec}
\end{figure}

 Tables \ref{tab:lorenz96_Gaussian} and \ref{tab:lorenz96_tvd} show that that CS-MIO achieves better performance than PySINDy under large noise in terms of the number of exactly recovered equations, i.e., the metric $A(\mathbf{\Gamma})$, under both noisy types. Note from $\sigma$=70 under Type 1 noise and $\sigma$=1 under Type 2 noise, CS-MIO fails to completely discover the 96 equations because the $\ell_1$ regularization fails to include all the ground truth terms within the first 100 significant terms.

Figure \ref{fig:lorenz96_50_trajec} shows in the form of Hoverm{\"o}ller plot the trajectory difference between the identified system and the ground truth. 
In particular, the CS-MIO identified systems at $\sigma$=50 in Table \ref{tab:lorenz96_Gaussian} and $\sigma$=0.8 in Table \ref{tab:lorenz96_tvd}, are used to run simulation in time interval $t\in[0,10]$ with time step $\Delta t=0.01$ sec. In the Hoverm{\"o}ller plot, the horizontal axis is the time and the vertical axis refers to the index of the state variables. The differences between the state values of ground truth $x_j(t)$ and those of the identified system $\hat{x}_j(t)$, $\Delta x_j(t) = x_j(t) - \hat{x}_j(t)$ for $j\in\{1,2,\cdots,96\}$, are shown with different colors. The more white-colored areas indicate the trajectory of the identified system agrees better with the ground truth. For example, it is seen roughly at $t\in[0,1]$ from Figure \ref{fig:lorenz96_Gaussian50_traj}, the identified system trajectory coincides well with the ground truth while the deviation starts to increase after that. We do not show the figure for SINDy since its identified equations result in unstable attractor. More details for the identified models of Lorenz 96 system by CS-MIO are in Appendix \ref{sec:appendix_lorenz96}.

% \subsection{Normal Forms, Bifurcations, and Parameterized Systems}

\subsection{Bifurcations and parameterized systems}

% \begin{table}[t]
% \footnotesize
%   \centering
%   \caption{Number of exactly recovered equations for Hopf Normal form: both $\mathbf{x}$ and $\mathbf{\dot{x}}$ can be measured; noise is added to $\mathbf{\dot{x}}$.}
%   \label{tab:hopf_Gaussian}
%   \begin{tabular}{rr|rr}
%       \toprule
%   Noise: $\sigma$  &	SNR &   PySINDy &   CS-MIO	\\
%   \midrule
% %0.001	&	137726.0217	&	2	&	2	\\
% %0.01	&	1377.2602	&	2	&	2	\\
% 0.1	&	13.773	&	\textbf{2}	&	\textbf{2}	\\
% 0.3	&	1.530	&	\textbf{2}	&	\textbf{2}	\\
% 0.5	&	0.551	&	\textbf{2}	&	\textbf{2}	\\
% 0.7	&	0.281	&	1	&	\textbf{2}	\\
% 1	&	0.138	&	0	&	\textbf{2}	\\
% %2	&	0.0344	&	0	&	2	\\
% 3	&	0.015	&	0	&	\textbf{2}	\\
% %4	&	0.0086	&	0	&	2	\\
% 5	&	0.006	&	0	&	\textbf{2}	\\
% %6	&	0.0038	&	0	&	2	\\
% 7	&	0.003	&	0	&	0	\\
% %8	&	0.0022	&	0	&	0	\\
% %9	&	0.0017	&	0	&	0	\\
% %10	&	0.0014	&	0	&	0	\\
% \bottomrule
% \end{tabular}
% \end{table}

Parameterized systems exhibit rich dynamic behaviors with various parameter values, which is known as bifurcations. We consider two examples of parameterized systems used in \cite{brunton2016discovering}. The first is the 2D Hopf normal form with bifurcation parameter $\mu$,
\begin{align}
    \dot{x} & = \mu x - \omega y + Ax(x^2+y^2),\\
    \dot{y} & = \omega x + \mu y + Ay(x^2+y^2).
\end{align}

\begin{figure}[h!]
\begin{subfigure}[c]{0.32\textwidth}
    \centering
    \hspace{-0.8cm}
    \includegraphics[width=2.2in]{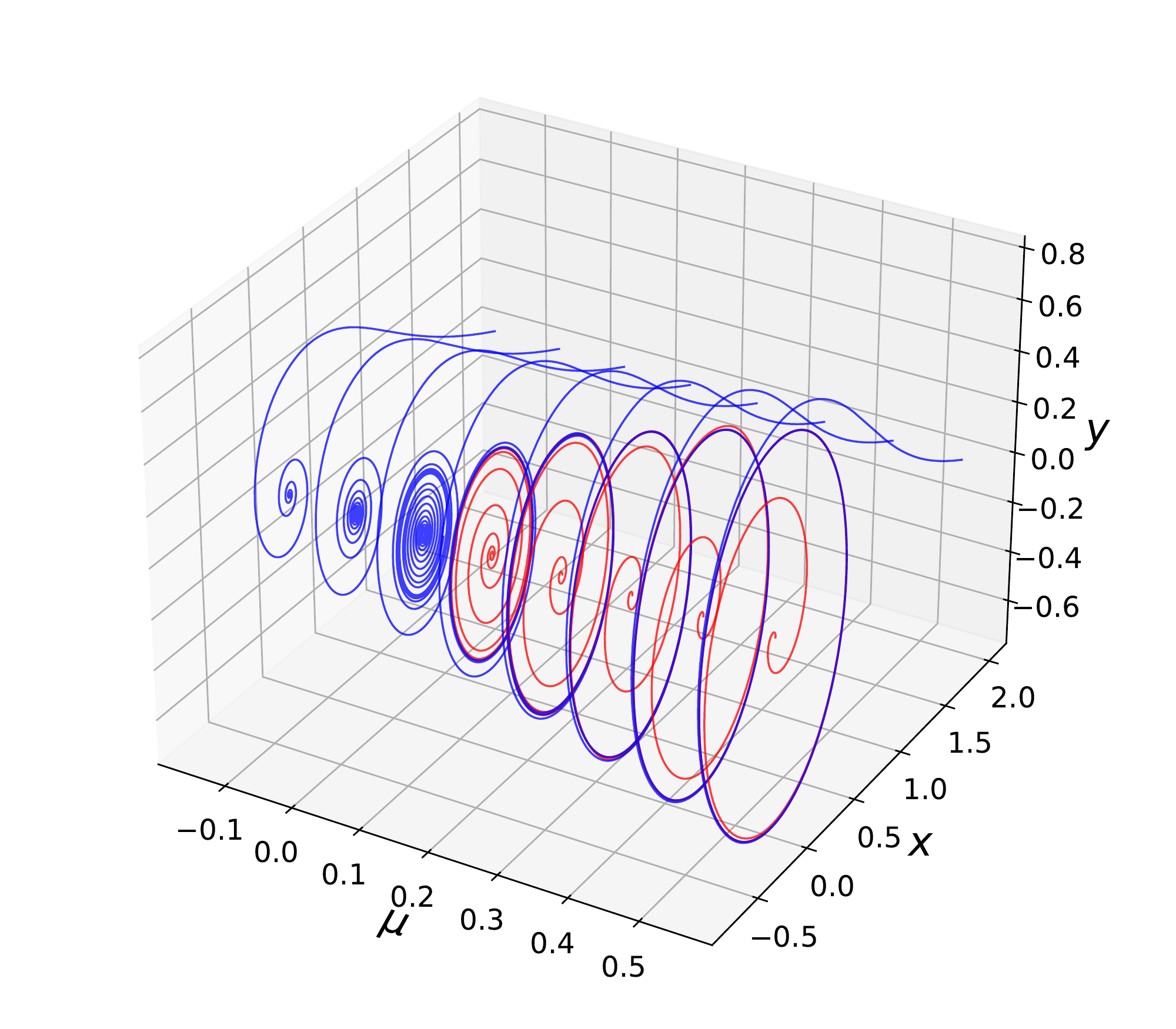}
    \caption{Trajectory of the ground truth.}
    \label{fig:hopf_groundtruth}
\end{subfigure}
\begin{subfigure}[c]{0.33\textwidth}
    \centering
    \hspace{-0.6cm}
    \includegraphics[width=2.2in]{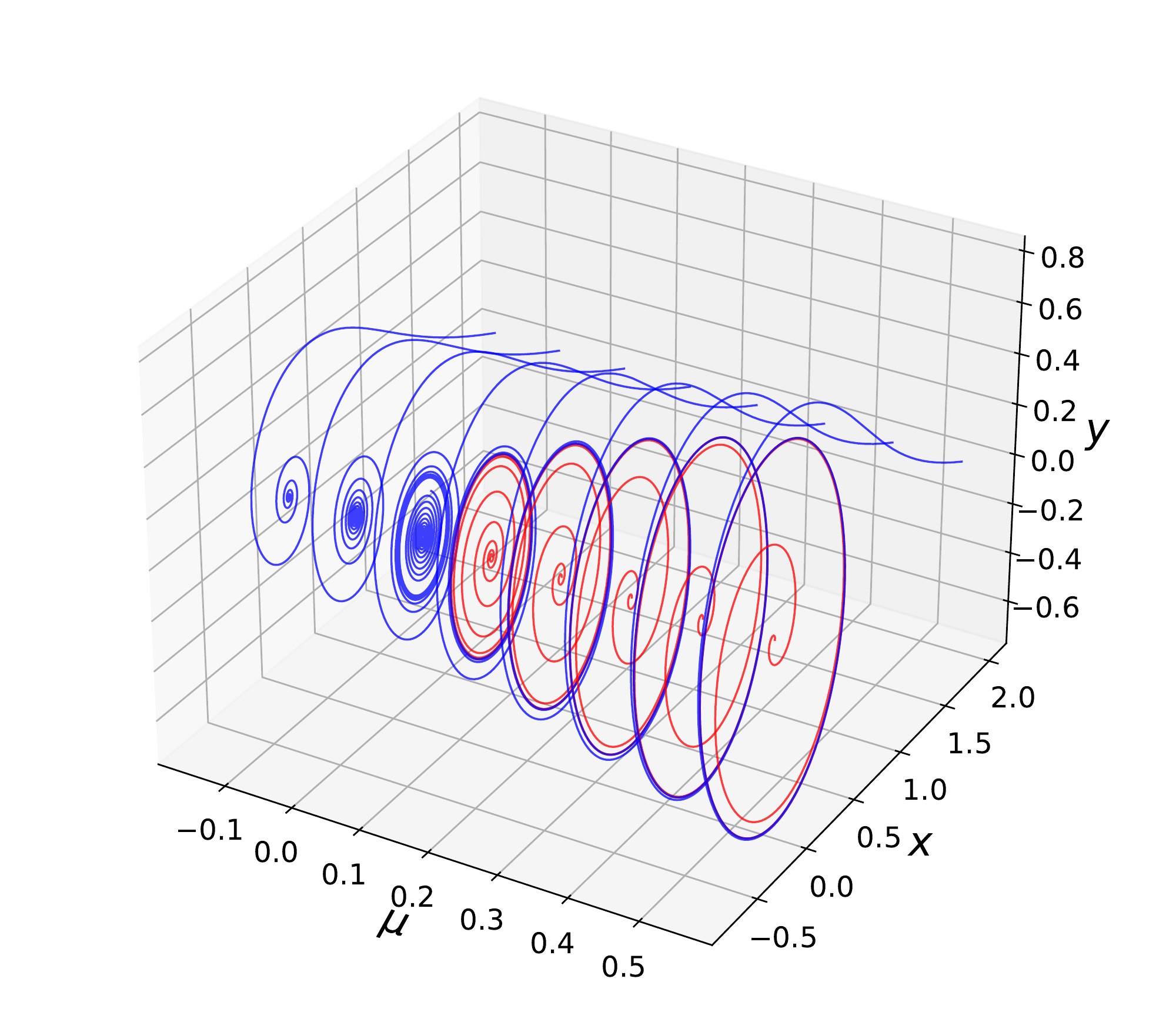}
    \caption{Type 1 noise at $\sigma$=3.}
    \label{fig:hopf_Gaussian6_CS_MIO}
\end{subfigure}    
\begin{subfigure}[c]{0.33\textwidth}
    \centering
    \hspace{-0.6cm}
    \includegraphics[width=2.2in]{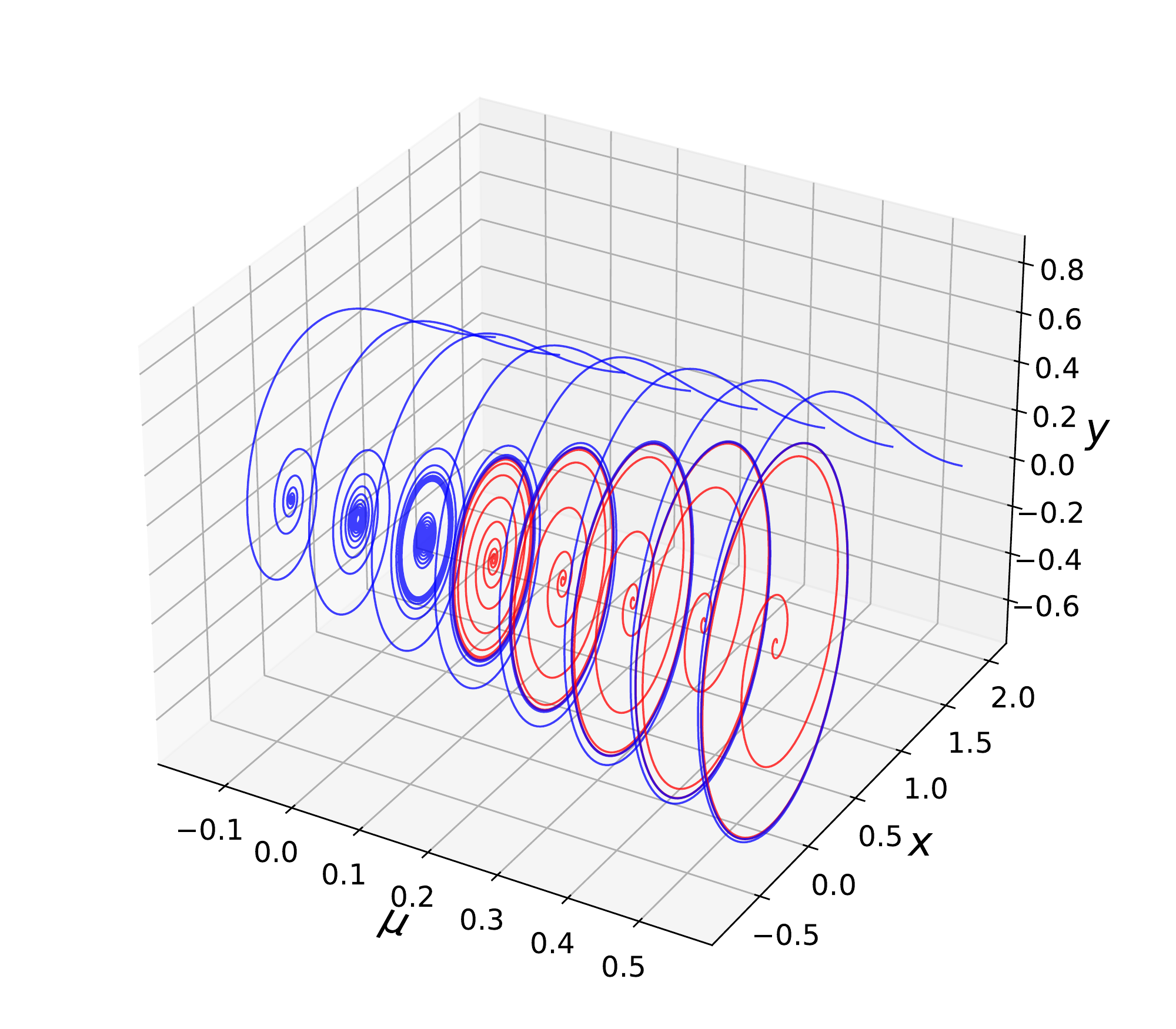}
    \caption{Type 2 noise at $\sigma$=0.015.}
    \label{fig:hopf_TVD0.015_CS_MIO}
\end{subfigure}    
\caption{Trajectory of CS-MIO discovered systems for Hopf normal form. (a) Trajectory of the ground truth full simulation. (b) Trajectory of the CS-MIO identified system under Type 1 noise at at $\sigma$=3. (c) Trajectory of the CS-MIO identified system under Type 2 noise at at $\sigma$=0.015.}
\label{fig:hopf_traj}
\end{figure}

To handle the bifurcation behaviors, the $\mu$ in the Hopf normal form is treated as additional state variables by adding dummy differential equation $\dot{\mu}=0$ to the system \cite{brunton2016discovering}. By adopting this setting, we used 14 values of $\mu$ to generate 14 datasets, with each dataset is collected using $\Delta t=0.0025$ in $t\in[0,75]$. 
We combine these datasets as a single training dataset to identify the governing equation as a function of state $\mathbf{x}$ and bifurcation parameter $\mu$, i.e., $\mathbf{\dot{x}} = \mathbf{f}(\mathbf{x},\mu)$. 
Tables \ref{tab:hopf_Gaussian} and \ref{tab:hopf_NonGaussian} present the number of exact recovered equations of Hopf normal forms under various noise SNRs Type 1 and Type 2 noises, respectively. Note herein we neglect the counting of the dummy differential equation in both examples. Under Type 1 noise, it is seen the lowest SRN can be as low as 0.015 for CS-MIO to exactly recovered all the equations. Figure \ref{fig:hopf_traj} shows the trajectory of the CS-MIO discovered systems under both Type 1 and 2 noise in comparison with the ground truth. More details for the identified models of Hopf normal form by CS-MIO are in Appendix \ref{sec:appendix_hopf}.

\begin{table}[h!]
\footnotesize
\caption{Comparison of the number of exactly recovered equations, i.e., the metric ${A}(\mathbf{\Gamma})$ in Eq.~\eqref{metric}, for Hopf Normal form. Compared with PySINDy, our CS-MIO method can correctly recover all the two equations, i.e., identifying the correcting $\mathbf{\Gamma}$ in Eq.~\eqref{eq:ode1}, under smaller SNR values.}
\begin{subtable}{0.43\textwidth}
   \caption{Results under Type 1 noise.}
  \label{tab:hopf_Gaussian}
  \begin{tabular}{cc|cc}
      \toprule
             &	 & \multicolumn{2}{c}{The metric ${A}(\mathbf{\Gamma})$} \\[2pt]
  Noise std  &	SNR & PySINDy  &   CS-MIO	\\
   \midrule
%0.001	&	137583.905	&	2	&	2	\\
%0.01	&	1375.839	&	2	&	2	\\
0.1	&	13.758	&	\textbf{2}	&	\textbf{2}	\\
0.3	&	1.529	&	\textbf{2}	&	\textbf{2}	\\
0.5	&	0.550	&	1	&	\textbf{2}	\\
0.7	&	0.281	&	1	&	\textbf{2}	\\
1	&	0.138	&	0	&	\textbf{2}	\\
2	&	0.034	&	0	&	\textbf{2}	\\
3	&	0.015	&	0	&	\textbf{2}	\\
4	&	0.009	&	0	&	0	\\
%5	&	0.006	&	0	&	0	\\
%6	&	0.0038	&	0	&	0	\\
%7	&	0.003	&	0	&	0	\\
%8	&	0.0022	&	0	&	0	\\
%9	&	0.0017	&	0	&	0	\\
%10	&	0.0014	&	0	&	0	\\
\bottomrule
\end{tabular}
\end{subtable}
\hspace{0.5cm}
\begin{subtable}{0.48\textwidth}
   \centering
  \caption{Results under Type 2 noise.}
  \label{tab:hopf_NonGaussian}
  \begin{tabular}{cc|cc}
      \toprule
               &	 & \multicolumn{2}{c}{The metric ${A}(\mathbf{\Gamma})$} \\[2pt]
  Noise std  &	SNR & PySINDy  &   CS-MIO	\\
   \midrule
0.001	&	120306.233	&	\textbf{2}	&	\textbf{2}	\\
0.003	&	13367.359	&	\textbf{2}	&	\textbf{2}	\\
0.005	&	4812.249	&	\textbf{2}	&	\textbf{2}	\\
0.007	&	2455.229	&	\textbf{2}	&	\textbf{2}	\\
0.010	&	1203.062	&	0	&	\textbf{2}	\\
0.013	&	711.871	&	0	&	\textbf{2}	\\
0.015	&	534.694	&	0	&	\textbf{2}	\\
0.017	&	416.285	&	0	&	0	\\
%0.02	&	300.766	&	0	&	0	\\
%0.1	&	12.031	&	0	&	0	\\
%1	    &	0.120	&	0	&	0	\\
\bottomrule
\end{tabular}
\end{subtable}
\end{table}

\begin{table}[h!]
\footnotesize
  \centering
  \caption{Comparison of the number of exactly recovered equations, i.e., the metric ${A}(\mathbf{\Gamma})$ in Eq.~\eqref{metric}, for the logistic map. Compared with PySINDy, our CS-MIO method can correctly recover the single equation, i.e., identifying the correcting $\mathbf{\Gamma}$ in Eq.~\eqref{eq:ode1}, under smaller SNR values.}
  \label{tab:logistic_Gaussian}
  \begin{tabular}{cc|cc}
      \toprule
             &	 & \multicolumn{2}{c}{The metric ${A}(\mathbf{\Gamma})$} \\[2pt]
  Noise std  &	SNR & PySINDy  &   CS-MIO	\\
   \midrule
%0.001	&	48377.506	&	1	&	1	\\
%0.01	&	481.877	&	1	&	1	\\
0.1	&   8.985	&	\textbf{1}	&	\textbf{1}	\\
0.2	&	3.619	&	\textbf{1}	&	\textbf{1}	\\
0.3	&	2.146	&	0	&	\textbf{1}	\\
0.4	&	1.455	&	0	&	\textbf{1}	\\
0.5	&	1.098	&	0	&	\textbf{1}	\\
0.6	&	0.874	&	0	&	\textbf{1}	\\
0.7	&	0.738	&	0	&	0	\\
%1	&	0.495	&	0	&	0	\\
\bottomrule
\end{tabular}
\end{table}

\begin{figure}[h!]
\begin{subfigure}{0.48\textwidth}
    \centering
    \includegraphics[width=3.0in]{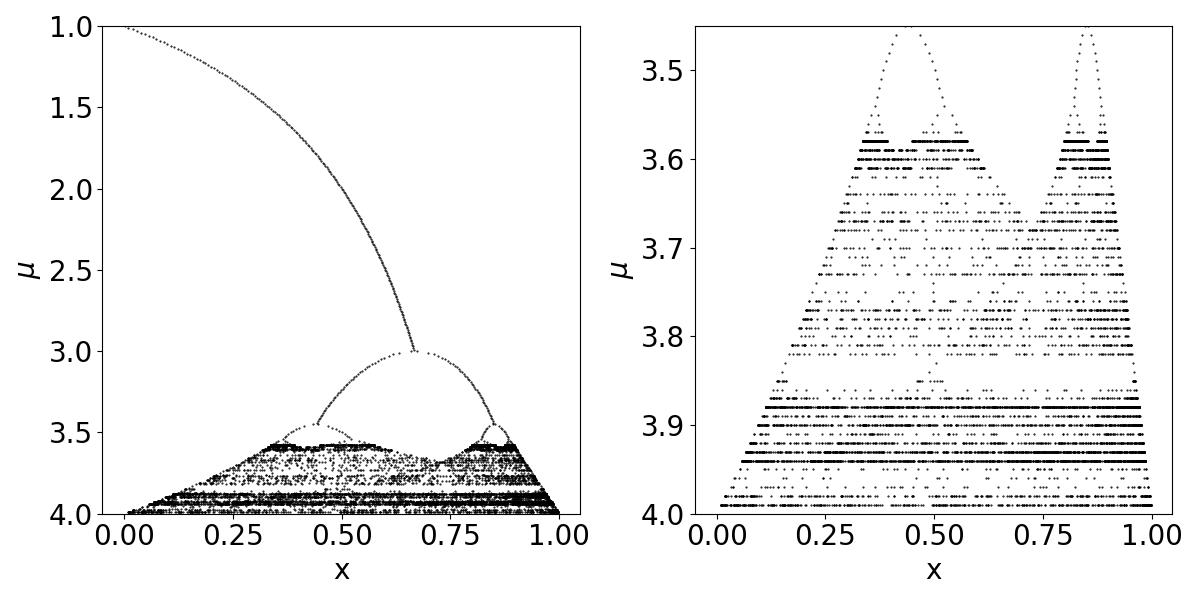}
    \caption{Trajectory of the ground truth.}
    \label{fig:logistic_ground_traj}
\end{subfigure}
\hspace{0.1cm}
\begin{subfigure}{0.48\textwidth}
\centering
    \includegraphics[width=3.0in]{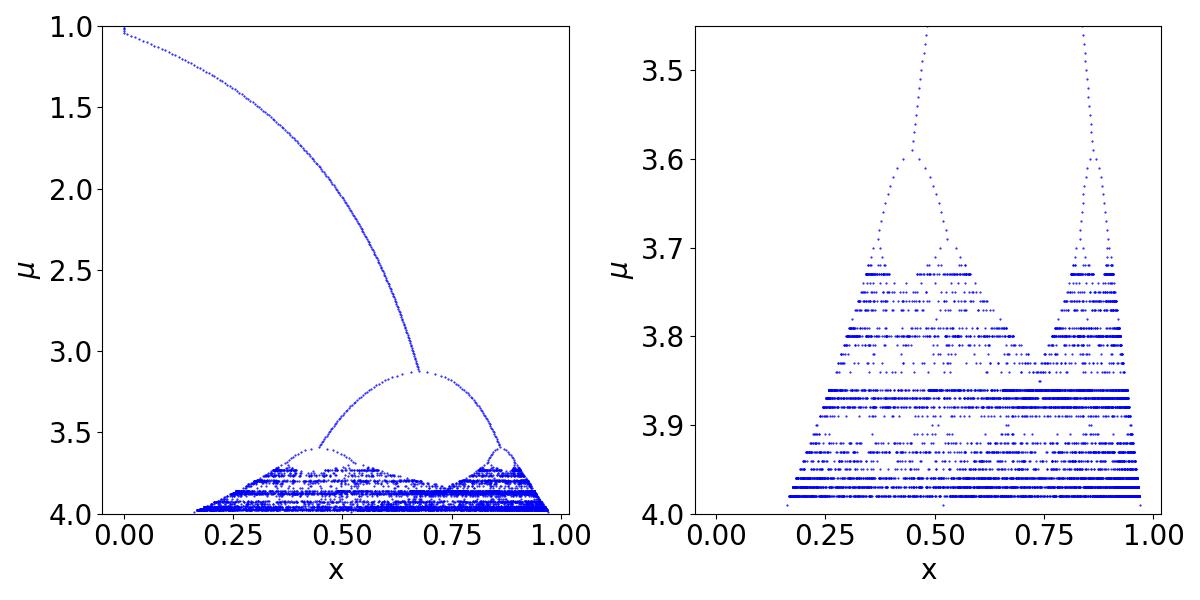}
    \caption{Trajectory of the CS-MIO identified system.}
    \label{fig:logistic_0.2_CS_MIO_traj}
\end{subfigure}
\caption{\footnotesize Trajectories of CS-MIO identified models for logistic map system under noise magnitude 0.2 in (b) and the comparison to the ground truth in (a) for ten values of $r$. }
\label{fig:logistic_simu}
\end{figure}

The second example is the 1D logistic map with stochastic forcing and bifurcation parameter $r$,
\begin{align}
    x_{n+1} = rx_{n}(1-x_n) + \eta_n.
\end{align}
where $\eta_n$ is the stochastic forcing and $r$ is the bifurcation parameter.
Similar manner is imposed to the bifurcation parameter $r$ in the logistic map with dummy equation $r_{n+1} = r_n$. Besides, 10 values of $r$ are used to collect the data. Within each dataset, we evolve the dynamical system for 1000 discrete steps. Note the logistic map is a discrete time dynamical system, so that there is only one manner for adding noise (herein $\eta_n$) to the state variables $\mathbf{x}$. Table \ref{tab:logistic_Gaussian} presents the comparison results under various SNRs of noise. CS-MIO exhibits strong capability of recovering governing equations from large noise. In Figure \ref{fig:logistic_simu}, we compare the trajectories of the CS-MIO identified system with the ground truth. Note we neglect the stochastic forcing $\eta_n$ when evolving the trajectories in both figures, namely $\eta_n=0$.
The right panel of both Figures \ref{fig:logistic_ground_traj} and \ref{fig:logistic_0.2_CS_MIO_traj} limit the $\mu$ in the range of $[3.5, 4]$ for clearer presentation. It can be seen the trajectory of CS-MIO identified system agrees well with the ground truth simulation. More details for the identified models of Logistic map system by CS-MIO are in Appendix \ref{sec:appendix_logistic}.

\subsection{PDE for vortex shedding behind a cylinder}

The last example system is the fluid dynamics for vortex shedding behind a cylinder which are high-dimensional partial differential equations. As discussed in \cite{brunton2016discovering}, the high-dimensional PDEs of cylinder dynamics can evolve on a a low-dimensional attractor governed by ordinary differential equations after dimension reduction using proper orthogonal decomposition (POD). The mean-field model using three POD modes as coordinate system is given as follows.
\begin{align}
    \dot{x} & = \mu x - \omega y + Axz,\label{eq:cylinder1}\\
    \dot{y} & = \omega x + \mu y + Ayz,\label{eq:cylinder2}\\
    \dot{z} & = -\lambda (z-x^2 -y^2).
\end{align}

\begin{table}
\footnotesize
  \centering
  \caption{Identified coefficients of CS-MIO and PySINDy on flow wake behind a cylinder. Quadratic terms are identified. The bold coefficients refer to those in the ground truth mean field model. CS-MIO can identify all the ground truth terms using less nonzeros in $\bf \Gamma$ in comparison with PySINDy.}
  \label{tab:cylinder_order2}
  %\resizebox{\columnwidth}{!}{
  \begin{tabular}{l|rrrrrrrr}
      \toprule
  Term  &	\multicolumn{2}{c}{Equation 1}	&&	\multicolumn{2}{c}{Equation 2}	&&    \multicolumn{2}{c}{Equation 3}	\\
  \cline{2-3} \cline{5-6} \cline{8-9}
  & PySINDy & CS-MIO   && PySINDy & CS-MIO   && PySINDy & CS-MIO\\
  \midrule
Bias	&	-0.1225	&	0	&&	-0.0569	&	0	&&	-21.9002	&	-20.8466	\\
$x$	&	\textbf{-0.0092}	&	\textbf{-0.0092}	&&	\textbf{1.0347}	&	\textbf{1.0346}	&&	-0.0009	&	0	\\
$y$	&	\textbf{-1.0224}	&	\textbf{-1.0225}	&&	\textbf{0.0047}	&	\textbf{0.0046}	&&	0	&	0	\\
$z$	&	-0.0009	&	0	&&	-0.0004	&	0	&&	\textbf{-0.3117}	&	\textbf{-0.2968}	\\
$x^2$	&	0	&	0	&&	0	&	0	&&	\textbf{0.0011}	&	\textbf{0.0011}	\\
$xy$	&	0	&	0	&&	0	&	0	&&	0.0002	&	0	\\
$xz$	&	\textbf{0.0002}	&	\textbf{0.0002}	&&	0.0022	&	0.0022	&&	0	&	0	\\
$y^2$	&	0	&	0	&&	0	&	0	&&	\textbf{0.0009}	&	\textbf{0.0009}	\\
$yz$	&	-0.0019	&	-0.0019	&&	\textbf{-0.0018}	&	\textbf{-0.0018}	&&	0	&	0	\\
$z^2$	&	0	&	0	&&	0	&	0	&&	-0.0011	&	-0.0010\\
\bottomrule
\end{tabular}
%}
\end{table}

\begin{table}
\footnotesize
  \centering
  \caption{Identified coefficients of CS-MIO and PySINDy on flow wake behind a cylinder. Cubic terms are identified for the mean field model. The bold terms indicate the ground truth. CS-MIO can identify equal or more ground truth terms using less nonzeros in comparison with PySINDy.}
  \label{tab:cylinder_3order}
  %\resizebox{\columnwidth}{!}{
  \begin{tabular}{l|rrrrrrrr}
      \toprule
  Term  &	\multicolumn{2}{c}{Equation 1}	&&	\multicolumn{2}{c}{Equation 2}	&&    \multicolumn{2}{c}{Equation 3}	\\
   \cline{2-3} \cline{5-6} \cline{8-9}
   & PySINDy & CS-MIO   && PySINDy & CS-MIO   && PySINDy & CS-MIO\\
   \midrule
Bias	&	0	&	0	&&	0	&	0	&&	0	&	-9.66082	\\
$x$	&	\textbf{0}	&	\textbf{0}	&&	\textbf{0}	&	\textbf{1.02896}	&&	0	&	0	\\
$y$	&	\textbf{-1.04203}	&	\textbf{-0.21545}	&&	\textbf{0.00621}	&	\textbf{0.24547}	&&	0.00025	&	0	\\
$z$	&	0.00002	&	0	&&	-0.00004	&	0	&&	\textbf{0.47502}	&	\textbf{0.19082}	\\
$x^2$	&	0	&	0	&&	0	&	0	&&	\textbf{0.00006}	&	\textbf{0.00047}	\\
$xy$	&	0	&	0	&&	0	&	0	&&	-0.00019	&	0	\\
$xz$	&	0.00138	&	0.00275	&&	-0.00744	&	0.00222	&&	0	&	0	\\
$y^2$	&	0	&	0	&&	0	&	0	&&	\textbf{-0.00006}	&	\textbf{0.00038}	\\
$yz$	&	-0.00367	&	0.00396	&&	-0.00366	&	0	&&	0	&	0	\\
$z^2$	&	0	&	0	&&	0	&	0	&&	0.00532	&	0.00296	\\
$x^3$	&	\textbf{0}	&	\textbf{0}	&&	0.00005	&	0	&&	0	&	0	\\
$x^2y$	&	0	&	-0.00004	&&	\textbf{0}	&	\textbf{-0.00001}	&&	0	&	0	\\
$x^2z$	&	0	&	0	&&	0	&	0	&&	-0.00003	&	-0.00002	\\
$xy^2$	&	\textbf{0}	&	\textbf{0}	&&	0.00005	&	0	&&	0	&	0	\\
$xyz$	&	0	&	0	&&	0	&	0	&&	-0.00002	&	-0.00002	\\
$xz^2$	&	0.00001	&	0.00002	&&	-0.00002	&	0	&&	0	&	0	\\
$y^3$	&	0	&	-0.00004	&&	\textbf{0}	&	\textbf{-0.00001}	&&	0	&	0	\\
$y^2z$	&	0	&	0	&&	0	&	0	&&	-0.00002	&	-0.00002	\\
$yz^2$	&	-0.00002	&	0	&&	-0.00002	&	0	&&	0	&	0	\\
$z^3$	&	0	&	0	&&	0	&	0	&&	0.00001	&	0.00001	\\
%\midrule
%Number Nonzero Terms	&	\textbf{6}	&	\textbf{6}	&&	8	&	\textbf{5}	&&	10	&	\textbf{9}	\\
%Number of Correct Terms/Ground Truth	&	\textbf{1/4}	&	\textbf{1/4}	&&	1/4	&	\textbf{4/4}	&&	\textbf{3/3}	&	\textbf{3/3}	\\
\bottomrule
\end{tabular}
%}
\end{table}

\begin{figure}[t!]
\begin{subfigure}[c]{0.32\textwidth}
    \centering
    \includegraphics[width=2in]{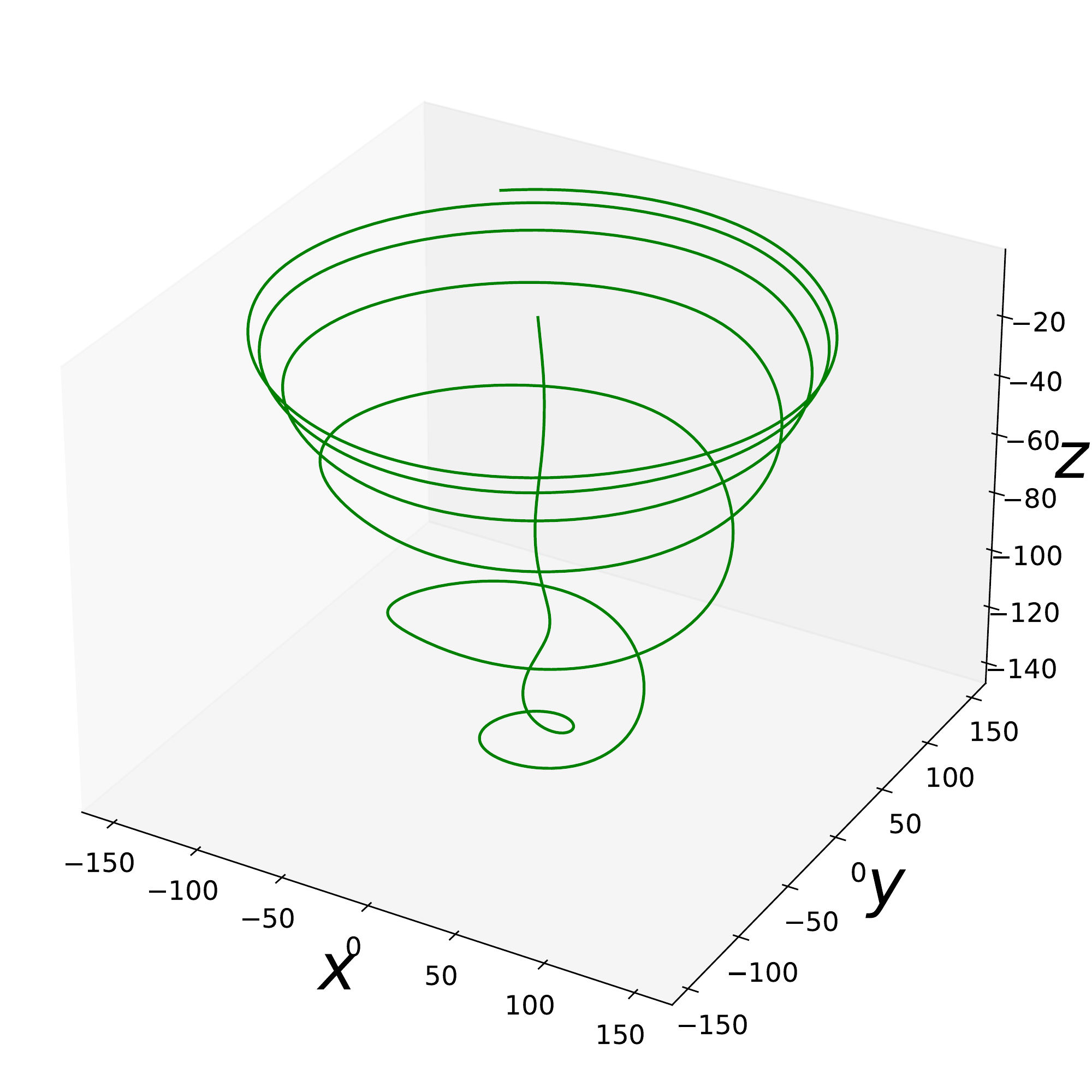}    
    \caption{Ground truth trajectory.}\label{fig:simu_cylinder}
\end{subfigure}
\begin{subfigure}[c]{0.32\textwidth}
    \centering
    \includegraphics[width=2in]{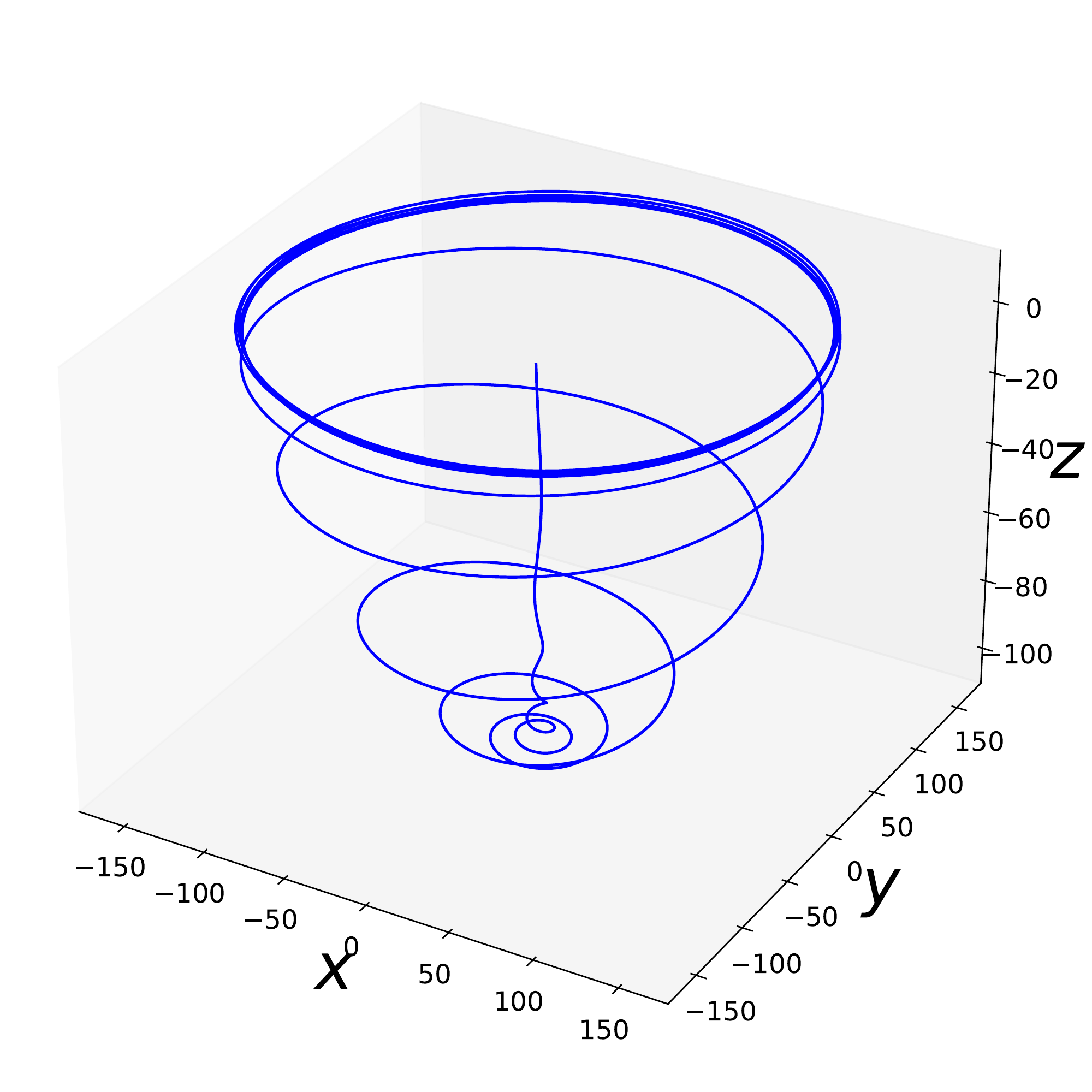}    
    \caption{Discovered 2nd order system}\label{fig:CS_MIO_cylinder_2order}
\end{subfigure}
\begin{subfigure}[c]{0.32\textwidth}
    \centering
    \includegraphics[width=2in]{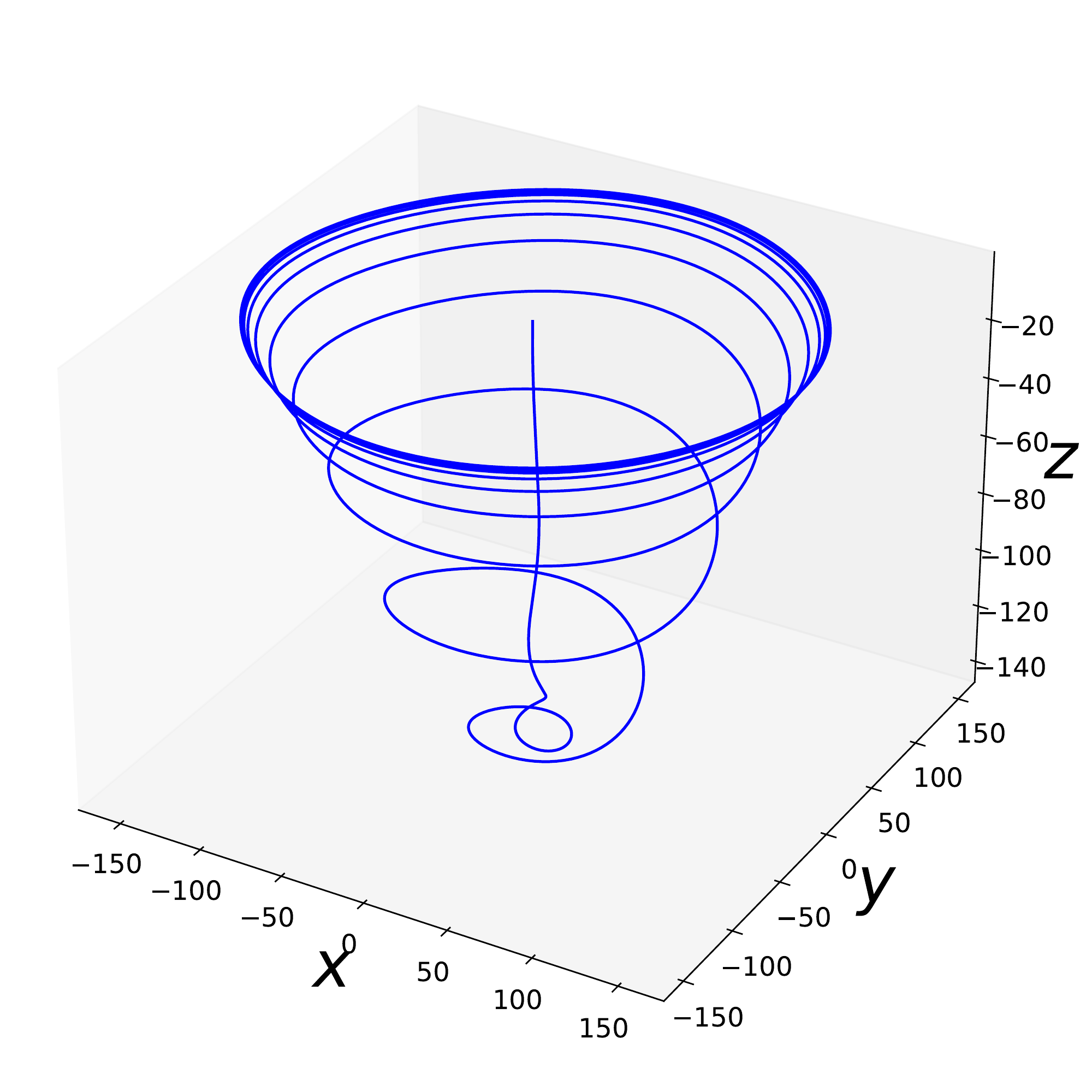}    
    \caption{Discovered 3rd order system}\label{fig:CS_MIO_cylinder_3order}
\end{subfigure}
\caption{\footnotesize Trajectories of the full simulation and the CS-MIO identified system for cylinder dynamics.}
\label{fig:cylinder_traj}
\end{figure}

Herein we use the same dataset used in \cite{brunton2016discovering}, which is originally generated using direct numerical simulations of the 2D Navier-Stokes equations originally by \cite{taira2007immersed,colonius2008fast}. We do not employ either Type 1 or Type 2 noise instead of using this single dataset.
We identify the differential equations using CS-MIO as presented in Table \ref{tab:cylinder_order2}. Note here only second order polynomials are used for CS-MIO and PySINDy. In this case, neither CS-MIO nor PySINDy is able to exactly identify the differential equations. However, it is seen that CS-MIO uses 4, 4 and 5 terms, respectively, for each equation to include all those in the ground truth, while PySINDy uses 6, 6 and 7 terms and includes many false terms. 
In a word, CS-MIO uses much less nonlinear terms to include all those in the ground truth than PySINDy.

If $\lambda$ is large, then the dynamics on $z$ coordinate is fast, resulting in the quick transient dynamics from the mean flow to the parabolic slow manifold, that is $z=x^2+y^2$ given by the amplitude of the vortex shedding. This dynamics are seen in Figure \ref{fig:simu_cylinder} as the sharp decreasing along $z$ coordinate and then correcting to the parabolic slow manifold. If substituting $z=x^2+y^2$ into Equations \ref{eq:cylinder1} and \ref{eq:cylinder2}, we obtain a Hopf normal form system on the slow manifold, which include cubic nonlinearities.
We thus set polynomial order to be three in both CS-MIO and PySINDy for recovery. The results are shown in Table \ref{tab:cylinder_3order}. In this case, both CS-MIO and PySINDy fail to include all the ground truth terms although they all involve many redundant terms. This is reasonable since higher order nonlinearities can express the dynamics of lower order nonlinearities. From Figures \ref{fig:CS_MIO_cylinder_2order} and \ref{fig:CS_MIO_cylinder_3order}, it can be seen the CS-MIO identified system agree almost perfectly with the full simulation using the original dataset.

\section{Conclusion}

We have developed a compressive-sensing-assisted mixed-integer optimization method for recovery of dynamical systems from highly noisy data. As there remain many unknown governing equations across various disciplines in science and engineering, our developed method is critical for uncovering the unknown equations from
the noisy data that is practically observed in such systems. The proposed method is developed grounded on the important foundation, that is, the identification of terms in the governing equations is essentially a discrete optimization problem. Because of this, our method is able to separately control the exact sparsity of the governing equations, and estimate the associated coefficients. This differs significantly from existing research where sparsity is incurred by penalty on the coefficients. We also combine the mixed-integer optimization with compressive sensing and other regularization techniques for enhancing the capability for dealing with highly noisy and high-dimensional problems. Case studies using the classical dynamical system examples demonstrate the powerful capability of the proposed method to uncover the governing equations under large noise, significantly outperforms the state-of-the-art method. This work opens several doors for future directions. First, advanced algorithms could be developed to enhance the efficiency of the method for large-scale instances of the studied problem. In addition, the domain knowledge for specifying the number of active terms can be used to discover new governing equations in specific fields. The construction of candidate terms using rich symbolic expression is further an exciting potential direction.

\section*{Acknowledgement}
This material is based upon work supported in part by the U.S. Department of Energy, Office of Science, Office of Advanced Scientific Computing Research, and by the Laboratory Directed Research and Development program at the Oak Ridge National Laboratory, which is operated by UT-Battelle, LLC, for the U.S.~Department of Energy under Contract DE-AC05-00OR22725. This manuscript is partially supported by the Science Alliance GATE (Graduate Advancement, Training and Education) Award of the University of Tennessee Knoxville.

% One limitation of our method is that, when the number of active terms $k$ is large, the branch-and-cut algorithm used in CS-MIO might get stuck and cannot obtain good solutions with suboptimality guaranteed. This is essentially caused by the $\mathcal{NP}$-hardness of the studied sparse regression problem. While this represents a challenge, approaches for screening can be developed more explicitly to shrink the size of the candidate terms to a relatively small scale. Nevertheless, the parsimonious model assumption holds for almost all the known governing equations as discussed, as such there can be little change that we need to try large $k$.
% 
\bibliographystyle{abbrv}
%\bibliography{pnas-sample,pnas-sample1}

\newpage

\FloatBarrier

\begin{appendices}

\FloatBarrier
\section{Additional results for the Chaotic Lorenz 3 System}\label{sec:appendix_lorenz3}

We provide additional results for the Chaotic Lorenz 3 system. 

\begin{table}[h!]
\centering
\footnotesize
  \caption{Identified coefficients of Lorenz 3 system using CS-MIO under Type 1 noise.}
  \label{tab:lorenz3_Gaussian_coefs}
  \begin{tabular}{lr|rrr}
    \toprule
      \multirow{3}{*}{Noise: $\sigma$}  &	\multirow{3}{*}{SNR} &   $x$ &   $y$ &	\\
        &   &   $x$ &  $y$  &   $xz$    \\
        &   &   $z$ &   $xy$    &   \\
     \midrule
%      \multicolumn{2}{c|}{\multirow{3}{*}{Ground Truth}}	&	-10	&	10	&		\\
% 	&		&	28	&	-1	&	-1	\\
% 	&		&	-8/3	&	1	&		\\
% \midrule
\multirow{3}{*}{0.01}	&	 \multirow{3}{*}{41914317.129}	&	-10.0000	&	10.0000	&		\\
	&		&	28.0000	&	-1.0000	&	-1.0000	\\
	&		&	-2.6667	&	1.0000	&		\\
\midrule
\multirow{3}{*}{0.1}	&	 \multirow{3}{*}{419143.171}	&	-9.9999	&	9.9999	&		\\
	&		&	28.0004	&	-1.0001	&	-1.0000	\\
	&		&	-2.6667	&	1.0000	&		\\
\midrule
\multirow{3}{*}{1}	&	 \multirow{3}{*}{4191.621}	&	-9.9991	&	9.9990	&		\\
	&		&	28.0043	&	-1.0010	&	-1.0001	\\
	&		&	-2.6666	&	1.0000	&		\\
\midrule
\multirow{3}{*}{10}	&	 \multirow{3}{*}{41.916}	&	-9.9908	&	9.9900	&		\\
	&		&	28.0426	&	-1.0103	&	-1.0009	\\
	&		&	-2.6658	&	1.0000	&		\\
\midrule
\multirow{3}{*}{50}	&	 \multirow{3}{*}{1.677}	&	-9.9540	&	9.9499	&		\\
	&		&	28.2131	&	-1.0516	&	-1.0044	\\
	&		&	-2.6624	&	1.0001	&		\\
\midrule
\multirow{3}{*}{100}	&	 \multirow{3}{*}{0.419}	&	-9.9081	&	9.8999	&		\\
	&		&	28.4263	&	-1.1032	&	-1.0089	\\
	&		&	-2.6581	&	1.0001	&		\\
\midrule
\multirow{3}{*}{150}	&	 \multirow{3}{*}{0.186}	&	-9.8621	&	9.8498	&		\\
	&		&	28.6394	&	-1.1548	&	-1.0133	\\
	&		&	-2.6538	&	1.0002	&		\\
\midrule
\multirow{3}{*}{200}	&	 \multirow{3}{*}{0.105}	&	-9.8161	&	9.7998	&		\\
	&		&	28.8525	&	-1.2065	&	-1.0177	\\
	&		&	-2.6495	&	1.0002	&		\\
\midrule
\multirow{3}{*}{250}	&	 \multirow{3}{*}{0.067}	&	-9.7701	&	9.7497	&		\\
	&		&	29.0657	&	-1.2581	&	-1.0222	\\
	&		&	-2.6452	&	1.0003	&		\\
\midrule
\multirow{3}{*}{300}	&	 \multirow{3}{*}{0.047}	&	-9.7242	&	9.6997	&		\\
	&		&	29.2788	&	-1.3097	&	-1.0266	\\
	&		&	-2.6409	&	1.0003	&		\\
\bottomrule
    \end{tabular}
\end{table}

\begin{figure}[h!]
\begin{subfigure}[c]{0.75\textheight}
\centering
    \includegraphics[width=4.7in]{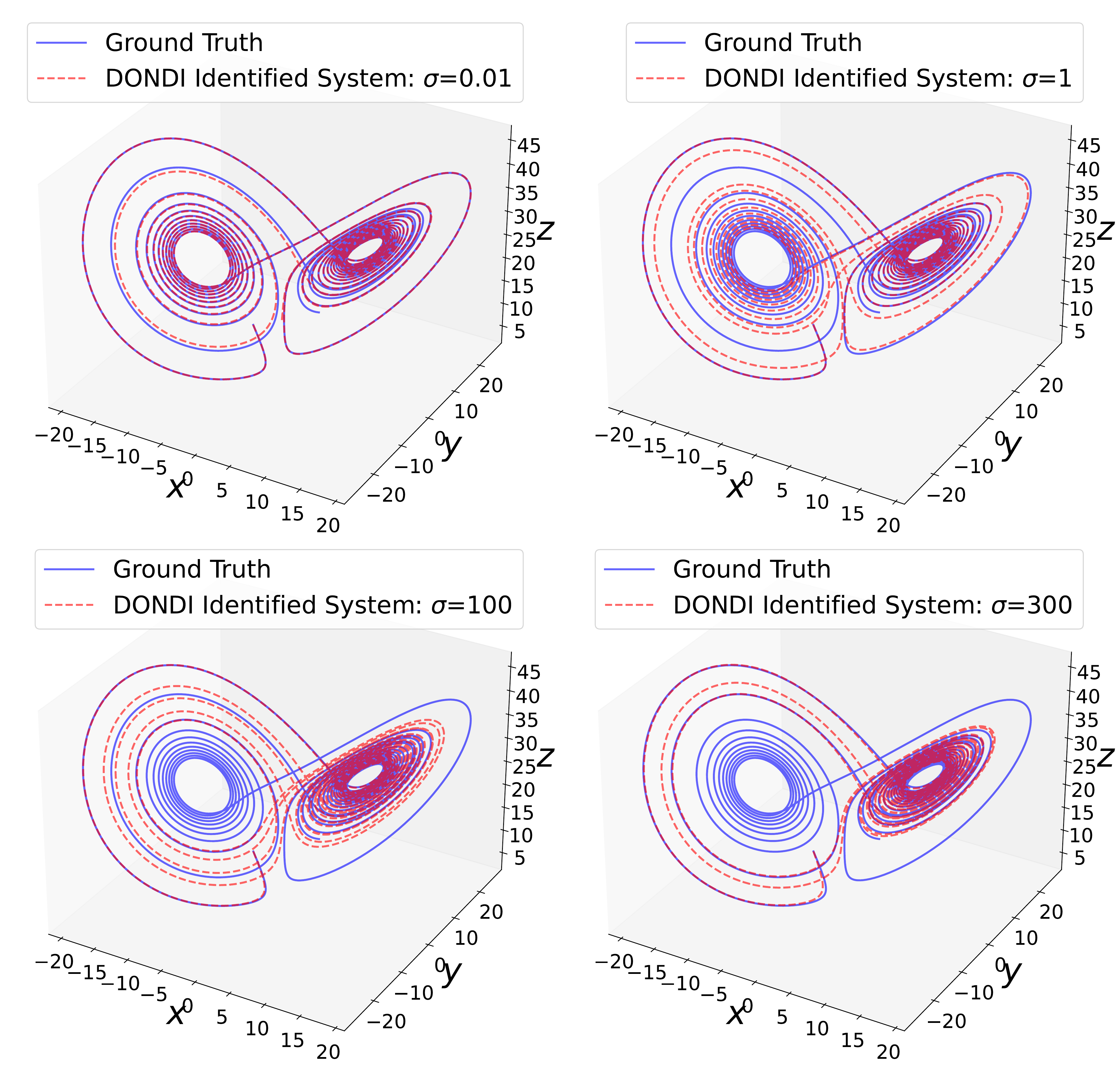}
    \caption{Trajectories of the CS-MIO identified system and the ground truth.}
    \label{fig:lorenz3_Gaussian_simu}
\end{subfigure}
\begin{subfigure}[c]{1\textwidth}
    \centering
    \includegraphics[width=4.7in]{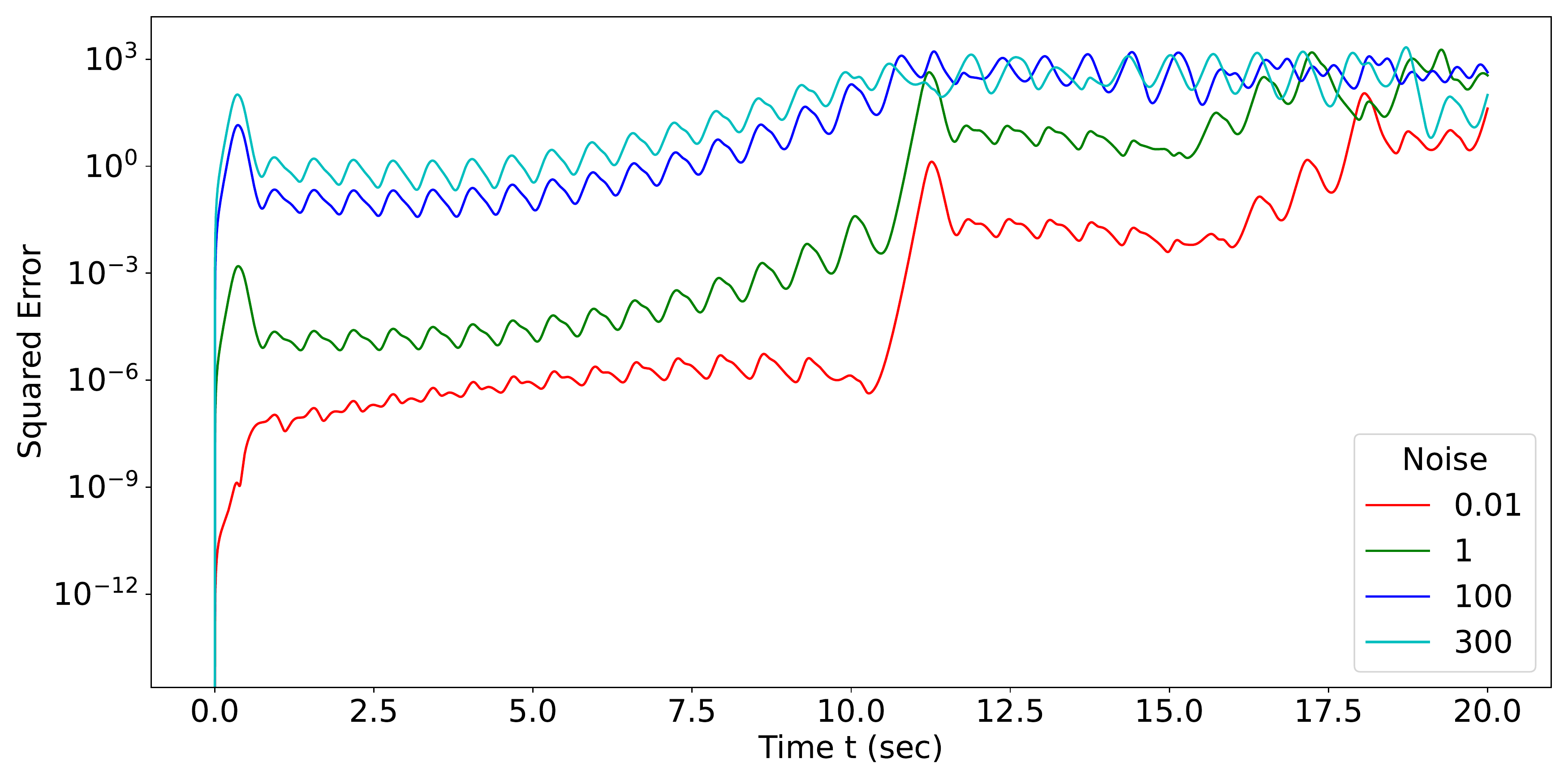}    \caption{$\ell_2$ error between the CS-MIO identified system and the ground truth.}
    \label{fig:lorenz3_Gaussian_l2error}
\end{subfigure}
\caption{Simulation results of the CS-MIO identified Lorenz 3 system comparing to the ground truth from $t=0$ to $t=20$ under Type 1 noise with four noise magnitudes $\sigma$: 0.01, 1, 100, and 300. The exact recovery fails when $\sigma$ is larger than 300.
(a) Trajectories of the CS-MIO identified system (red dashed) and ground truth (blue solid).
(b) $\ell_2$ error vs time of the trajectories of the recovered Lorenz 3 system ($\hat{\mathbf{x}}(t)$) comparing to the ground truth ($\mathbf{x}(t)$), i.e., $||\hat{\mathbf{x}}(t) - \mathbf{x}(t)||_2^2$ as a function of $t$ from $t=0$ to $t=20$.}
\label{fig:lorenz3_Gaussian_traj}
\end{figure}

\begin{table}[h!]
\footnotesize
  \centering
  \caption{Identified coefficients of Lorenz 3 system using CS-MIO under Type 2 noise.}
  \label{tab:lorenz3_tvd_coefs}
  \begin{tabular}{rr|rrr}
    \toprule
        \multirow{3}{*}{Noise: $\sigma$}  &	\multirow{3}{*}{SNR} &   $x$ &   $y$ &	\\
        &   &   $x$ &  $y$  &   $xz$    \\
        &   &   $z$ &   $xy$    &   \\
     \midrule
%      \multicolumn{2}{c|}{\multirow{3}{*}{Ground Truth}}	&	-10	&	10	&		\\
% 	&		&	28	&	-1	&	-1	\\
% 	&		&	-8/3	&	1	&		\\
% \midrule
\multirow{3}{*}{0.01}	&	\multirow{3}{*}{729427.159}	&	-9.9851	&	10.0000	&		\\
	&		&	27.6974	&	-0.8682	&	-0.9939	\\
	&		&	-2.6602	&	0.9997	&		\\
\midrule
\multirow{3}{*}{0.05}	&	\multirow{3}{*}{29178.677}	&	-9.9852	&	9.9999	&		\\
	&		&	27.6968	&	-0.8682	&	-0.9939	\\
	&		&	-2.6602	&	0.9997	&		\\
\midrule
\multirow{3}{*}{0.1}	&	\multirow{3}{*}{7295.616}	&	-9.9851	&	9.9997	&		\\
	&		&	27.6954	&	-0.8680	&	-0.9938	\\
	&		&	-2.6603	&	0.9997	&		\\
\midrule
\multirow{3}{*}{0.5}	&	\multirow{3}{*}{292.848}	&	-9.9791	&	9.9934	&		\\
	&		&	27.6635	&	-0.8596	&	-0.9929	\\
	&		&	-2.6608	&	0.9999	&		\\
\midrule
\multirow{3}{*}{1}	&	\multirow{3}{*}{73.982}	&	-9.9573	&	9.9730	&		\\
	&		&	27.5722	&	-0.8335	&	-0.9906	\\
	&		&	-2.6605	&	0.9997	&		\\
\midrule
\multirow{3}{*}{2}	&	\multirow{3}{*}{19.255}	&	-9.8670	&	9.8916	&		\\
	&		&	27.2258	&	-0.7319	&	-0.9823	\\
	&		&	-2.6571	&	0.9984	&		\\
    \bottomrule
    \end{tabular}
\end{table}

\begin{figure}[h!]
\begin{subfigure}[c]{1\textwidth}
    \centering
    \includegraphics[width=4.7in]{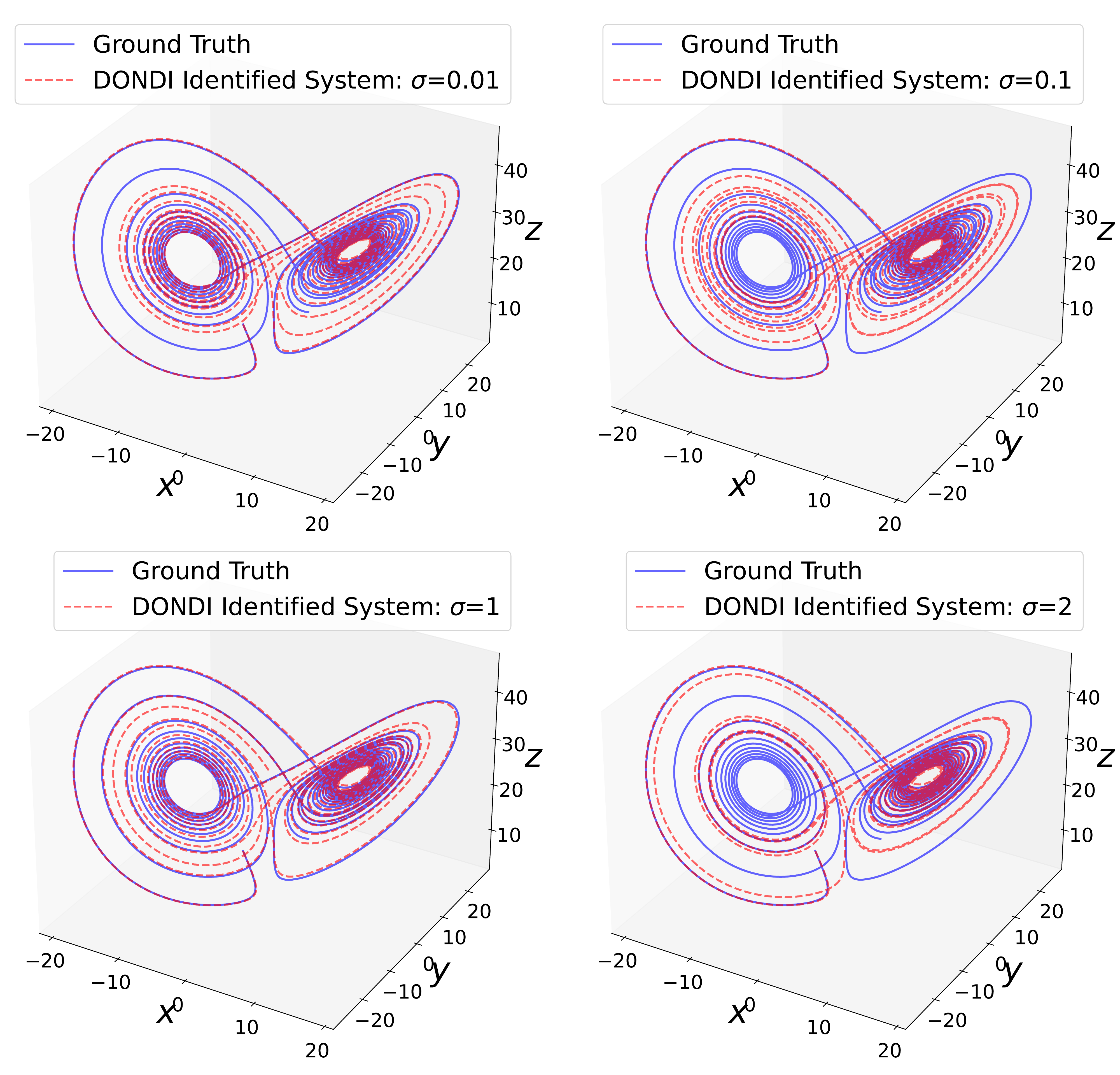}
    \caption{Trajectories of the CS-MIO identified system and the ground truth.}
    \label{fig:lorenz3_TVD_simu}
\end{subfigure}
\begin{subfigure}[c]{1\textwidth}
\centering
 \centering
    \includegraphics[width=4.6in]{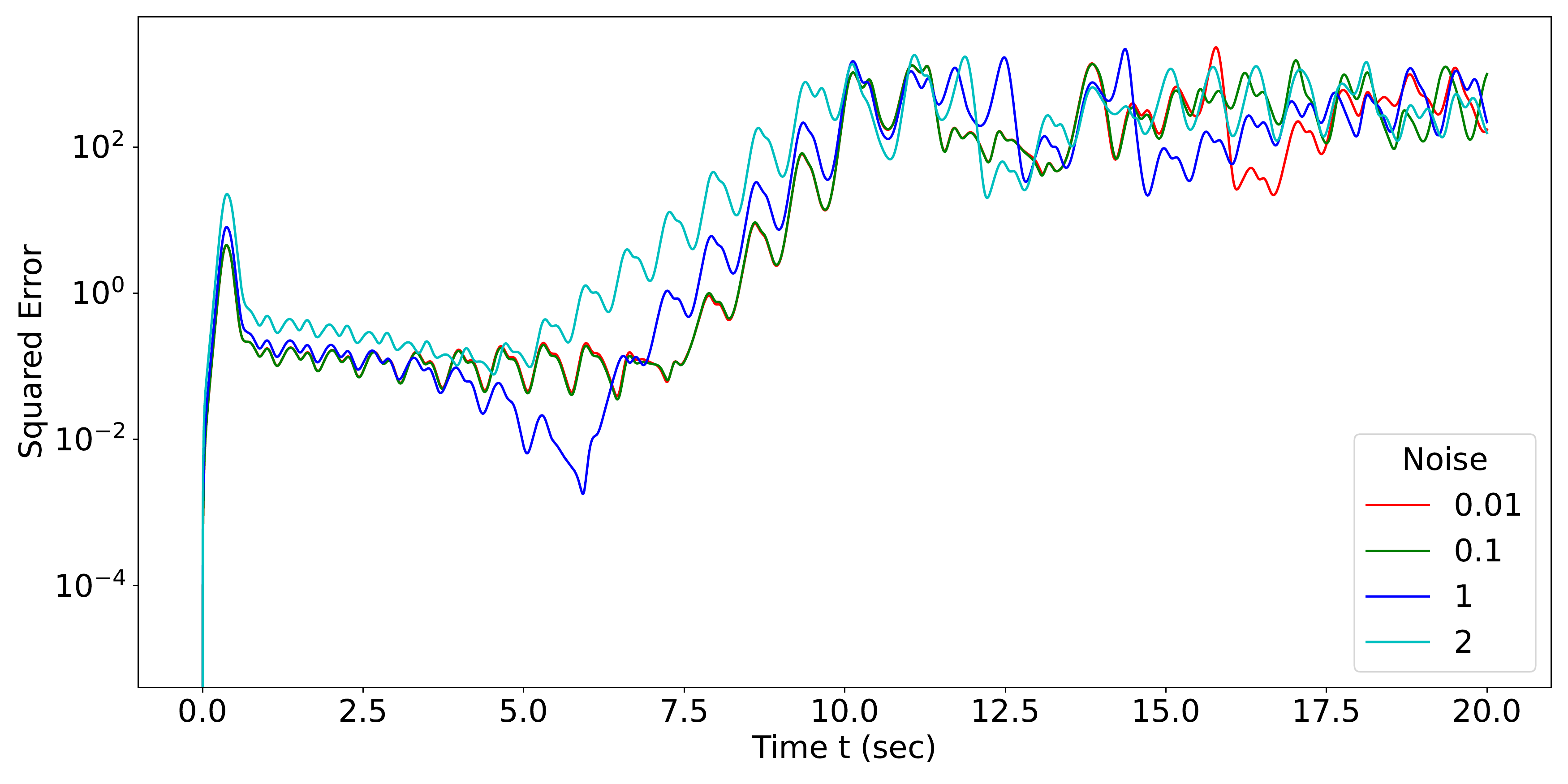}    
    \caption{$\ell_2$ error between the CS-MIO identified system and the ground truth.}
    \label{fig:lorenz3_TVD_l2error}
\end{subfigure}
\caption{Simulation results of the CS-MIO identified Lorenz 3 system comparing to the ground truth from $t=0$ to $t=20$ under Type 2 noise with four noise magnitudes $\sigma$: 0.01, 0.1, 1 and 2. The exact recovery fails when $\sigma$ ls larger than 2. 
(a) Trajectories of the CS-MIO identified system (red dashed) and ground truth (blue solid).
(b) $\ell_2$ error vs time of the trajectories of the recovered Lorenz 3 system ($\hat{\mathbf{x}}(t)$) comparing to the ground truth ($\mathbf{x}(t)$), i.e., $||\hat{\mathbf{x}}(t) - \mathbf{x}(t)||_2^2$ as a function of $t$ from $t=0$ to $t=20$.}
\label{fig:lorenz3_TVD_traj}
\end{figure}

\FloatBarrier
\section{Additional results for the Chaotic Lorenz 96 System}\label{sec:appendix_lorenz96}

\begin{table}[h]
\footnotesize
\centering
  \caption{Identified coefficients of Lorenz 96 system using CS-MIO under Type 1 noise with magnitude 50: Part A.}
  \label{tab:lorenz96_coefs_50_Gaussian_A}
  \begin{tabular}{c|rrrr}
    \toprule
      Equation Index $j$ &    $F$     &   $x_{j+1}x_{j-1}$   & $x_{j-2}x_{j-1} $  &  $x_j$ 	\\
     \midrule
     1	&	7.6636	&	0.9954	&	-0.9919	&	-0.9053	\\
2	&	8.2107	&	1.0137	&	-0.9957	&	-1.0136	\\
3	&	8.0884	&	1.0031	&	-0.9776	&	-1.0558	\\
4	&	7.7594	&	0.9958	&	-0.9925	&	-1.0204	\\
5	&	7.6274	&	0.9997	&	-1.0020	&	-0.9596	\\
6	&	7.6735	&	1.0009	&	-0.9790	&	-0.9568	\\
7	&	8.1661	&	1.0126	&	-1.0122	&	-1.0837	\\
8	&	7.7873	&	1.0004	&	-0.9930	&	-0.9782	\\
9	&	8.1277	&	0.9978	&	-1.0040	&	-1.0211	\\
10	&	8.2015	&	1.0310	&	-1.0276	&	-0.9769	\\
11	&	7.7901	&	0.9904	&	-1.0143	&	-0.8445	\\
12	&	8.0504	&	0.9844	&	-0.9902	&	-1.0707	\\
13	&	7.9518	&	0.9889	&	-1.0144	&	-1.0927	\\
14	&	7.8521	&	1.0055	&	-0.9897	&	-0.9773	\\
15	&	8.4207	&	0.9905	&	-0.9984	&	-1.0906	\\
16	&	8.2164	&	1.0162	&	-1.0171	&	-0.9698	\\
17	&	8.3087	&	0.9924	&	-0.9893	&	-1.0197	\\
18	&	8.4646	&	1.0119	&	-1.0123	&	-1.0714	\\
19	&	8.2501	&	0.9759	&	-0.9886	&	-1.0430	\\
20	&	7.8796	&	0.9919	&	-0.9873	&	-0.9722	\\
21	&	7.8357	&	0.9929	&	-1.0093	&	-0.9649	\\
22	&	8.5176	&	0.9749	&	-1.0268	&	-1.1041	\\
23	&	8.4018	&	1.0369	&	-1.0221	&	-1.0326	\\
24	&	7.9646	&	0.9937	&	-0.9830	&	-1.0083	\\
25	&	7.7756	&	0.9871	&	-0.9960	&	-0.9909	\\
26	&	7.6644	&	0.9907	&	-1.0124	&	-1.0218	\\
27	&	8.1968	&	1.0030	&	-1.0100	&	-1.0355	\\
28	&	8.2314	&	0.9985	&	-0.9925	&	-1.0586	\\
29	&	8.2346	&	1.0033	&	-1.0061	&	-1.0567	\\
30	&	7.6806	&	1.0109	&	-0.9958	&	-0.8968	\\
31	&	8.3632	&	1.0160	&	-1.0046	&	-1.0581	\\
32	&	8.1515	&	1.0083	&	-0.9858	&	-0.9228	\\
33	&	8.0453	&	1.0109	&	-1.0234	&	-0.9918	\\
34	&	7.8578	&	1.0029	&	-0.9876	&	-0.9044	\\
35	&	7.8325	&	1.0070	&	-1.0029	&	-0.9443	\\
36	&	8.0977	&	0.9873	&	-0.9925	&	-0.9986	\\
37	&	8.6882	&	0.9952	&	-1.0037	&	-1.0654	\\
38	&	7.6767	&	0.9983	&	-0.9685	&	-1.0238	\\
39	&	8.0614	&	0.9718	&	-0.9919	&	-0.9915	\\
40	&	7.8614	&	1.0046	&	-1.0125	&	-1.0181	\\
41	&	7.2833	&	0.9881	&	-0.9548	&	-0.9379	\\
42	&	8.0186	&	0.9879	&	-0.9965	&	-1.0273	\\
43	&	8.0823	&	0.9994	&	-1.0185	&	-1.0349	\\
44	&	7.9811	&	0.9758	&	-1.0048	&	-1.0275	\\
45	&	7.9975	&	0.9956	&	-0.9970	&	-0.9966	\\
46	&	7.7917	&	0.9986	&	-0.9901	&	-0.9353	\\
47	&	7.7668	&	0.9971	&	-0.9891	&	-1.0075	\\
48	&	7.5492	&	1.0106	&	-0.9717	&	-1.0031	\\
\bottomrule
    \end{tabular}
\end{table}

\begin{table}[h]
\footnotesize
\centering
  \caption{Identified coefficients of Lorenz 96 system using CS-MIO under Type 1 noise with magnitude 50: Part B.}
  \label{tab:lorenz96_coefs_50_Gaussian_B}
  \begin{tabular}{c|rrrr}
    \toprule
      Equation Index $j$ &    $F$     &   $x_{j+1}x_{j-1}$   & $x_{j-2}x_{j-1}$  &  $x_j$ 	\\
     \midrule
49	&	8.1873	&	1.0016	&	-1.0014	&	-0.9289	\\
50	&	7.8786	&	1.0036	&	-1.0096	&	-1.0100	\\
51	&	7.7901	&	1.0235	&	-0.9655	&	-1.0325	\\
52	&	7.4264	&	0.9828	&	-0.9877	&	-0.8955	\\
53	&	7.7280	&	0.9898	&	-0.9848	&	-0.9488	\\
54	&	8.0840	&	0.9995	&	-0.9977	&	-0.9669	\\
55	&	8.8206	&	0.9805	&	-1.0038	&	-1.0711	\\
56	&	8.5845	&	0.9934	&	-1.0120	&	-1.0968	\\
57	&	8.1714	&	0.9805	&	-0.9903	&	-0.9710	\\
58	&	8.5471	&	1.0031	&	-1.0284	&	-1.0648	\\
59	&	8.2529	&	1.0088	&	-0.9877	&	-0.9947	\\
60	&	8.2181	&	1.0111	&	-1.0064	&	-1.0332	\\
61	&	8.3878	&	1.0310	&	-1.0099	&	-1.0172	\\
62	&	8.4394	&	0.9898	&	-1.0120	&	-1.0398	\\
63	&	8.4044	&	1.0326	&	-0.9997	&	-1.0862	\\
64	&	8.1365	&	1.0058	&	-1.0072	&	-1.0033	\\
65	&	7.9563	&	0.9895	&	-0.9893	&	-1.0268	\\
66	&	8.3239	&	0.9874	&	-0.9911	&	-1.0302	\\
67	&	8.1219	&	1.0129	&	-1.0013	&	-1.0345	\\
68	&	7.8250	&	1.0037	&	-0.9819	&	-0.9692	\\
69	&	7.8334	&	1.0085	&	-1.0088	&	-0.9828	\\
70	&	8.0292	&	1.0082	&	-0.9822	&	-1.1236	\\
71	&	8.1128	&	0.9926	&	-1.0058	&	-1.0314	\\
72	&	8.2170	&	1.0018	&	-0.9909	&	-1.0880	\\
73	&	8.1572	&	1.0020	&	-0.9912	&	-0.9913	\\
74	&	7.9162	&	1.0058	&	-0.9838	&	-0.9800	\\
75	&	8.4242	&	1.0166	&	-1.0046	&	-1.1802	\\
76	&	8.2365	&	1.0021	&	-1.0086	&	-1.0709	\\
77	&	8.4704	&	1.0057	&	-1.0092	&	-1.0961	\\
78	&	8.2634	&	1.0109	&	-1.0043	&	-1.0300	\\
79	&	8.0802	&	0.9659	&	-0.9874	&	-1.0338	\\
80	&	8.2453	&	1.0183	&	-0.9971	&	-1.0583	\\
81	&	8.5512	&	0.9834	&	-1.0118	&	-1.0790	\\
82	&	7.8437	&	0.9997	&	-0.9978	&	-0.9929	\\
83	&	8.1340	&	0.9818	&	-0.9918	&	-0.9468	\\
84	&	8.3957	&	1.0119	&	-1.0095	&	-1.0589	\\
85	&	8.1792	&	1.0053	&	-0.9906	&	-1.0231	\\
86	&	8.1009	&	0.9933	&	-0.9980	&	-0.9385	\\
87	&	8.0357	&	0.9629	&	-1.0263	&	-0.8889	\\
88	&	7.9619	&	0.9825	&	-0.9979	&	-0.9141	\\
89	&	8.2250	&	0.9839	&	-0.9932	&	-1.0779	\\
90	&	7.8981	&	0.9894	&	-0.9881	&	-0.8814	\\
91	&	7.7766	&	0.9819	&	-1.0157	&	-0.9386	\\
92	&	7.9365	&	1.0125	&	-1.0174	&	-0.9876	\\
93	&	8.3197	&	1.0293	&	-1.0184	&	-1.0195	\\
94	&	7.6349	&	0.9901	&	-0.9642	&	-0.9583	\\
95	&	7.8329	&	1.0115	&	-1.0067	&	-0.9771	\\
96	&	7.9165	&	0.9803	&	-0.9942	&	-1.0080	\\
\bottomrule
    \end{tabular}
\end{table}

%coefficients

\begin{figure}[h!]
\begin{subfigure}[c]{0.7\textheight}
    \centering
    \includegraphics[width=5in]{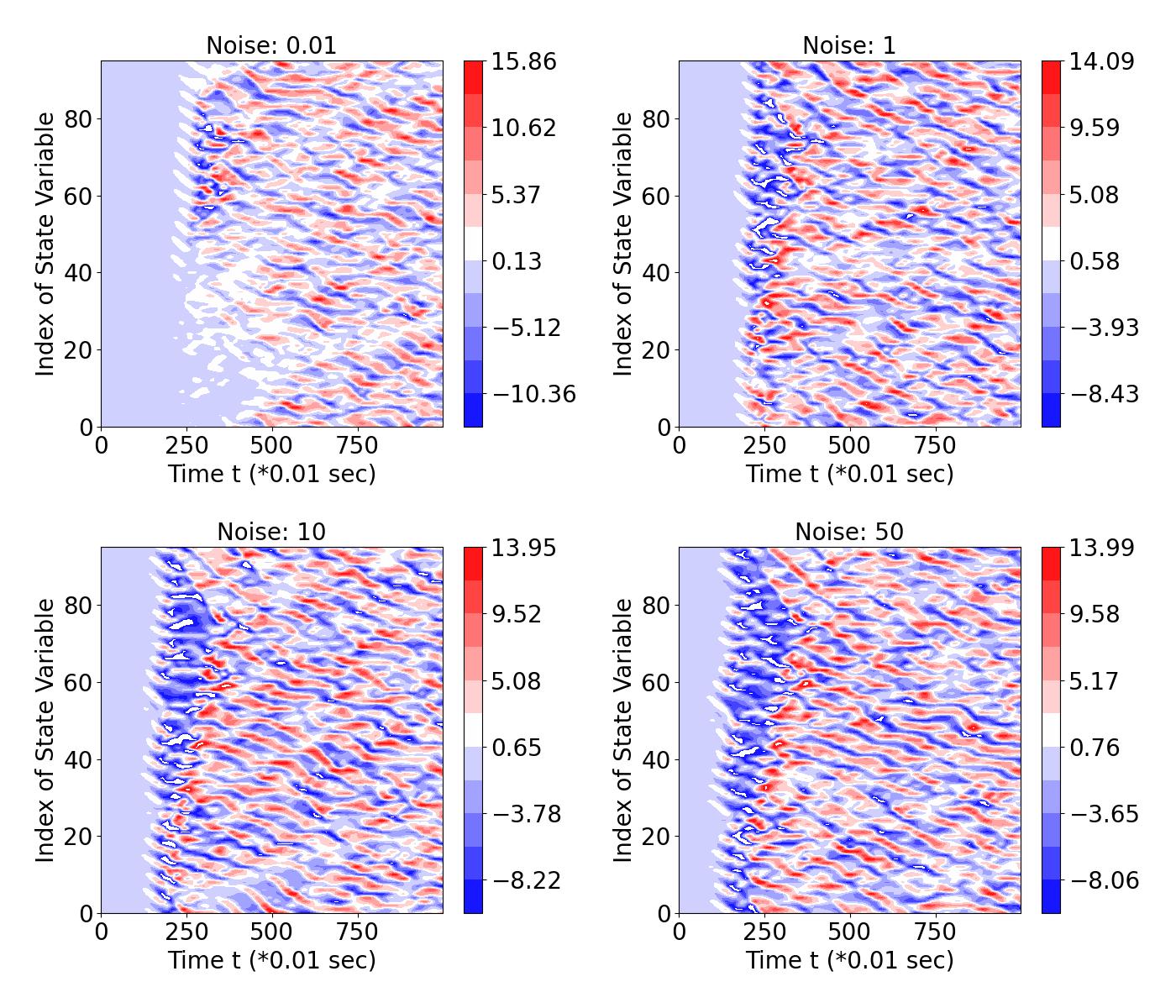}
    \caption{Hoverm{\"o}ller plot for difference between the identified system and ground truth.}
    \label{fig:lorenz96_Gaussian_diff}    
\end{subfigure}
\begin{subfigure}[c]{1\textwidth}
    \centering
    \includegraphics[width=4.3in]{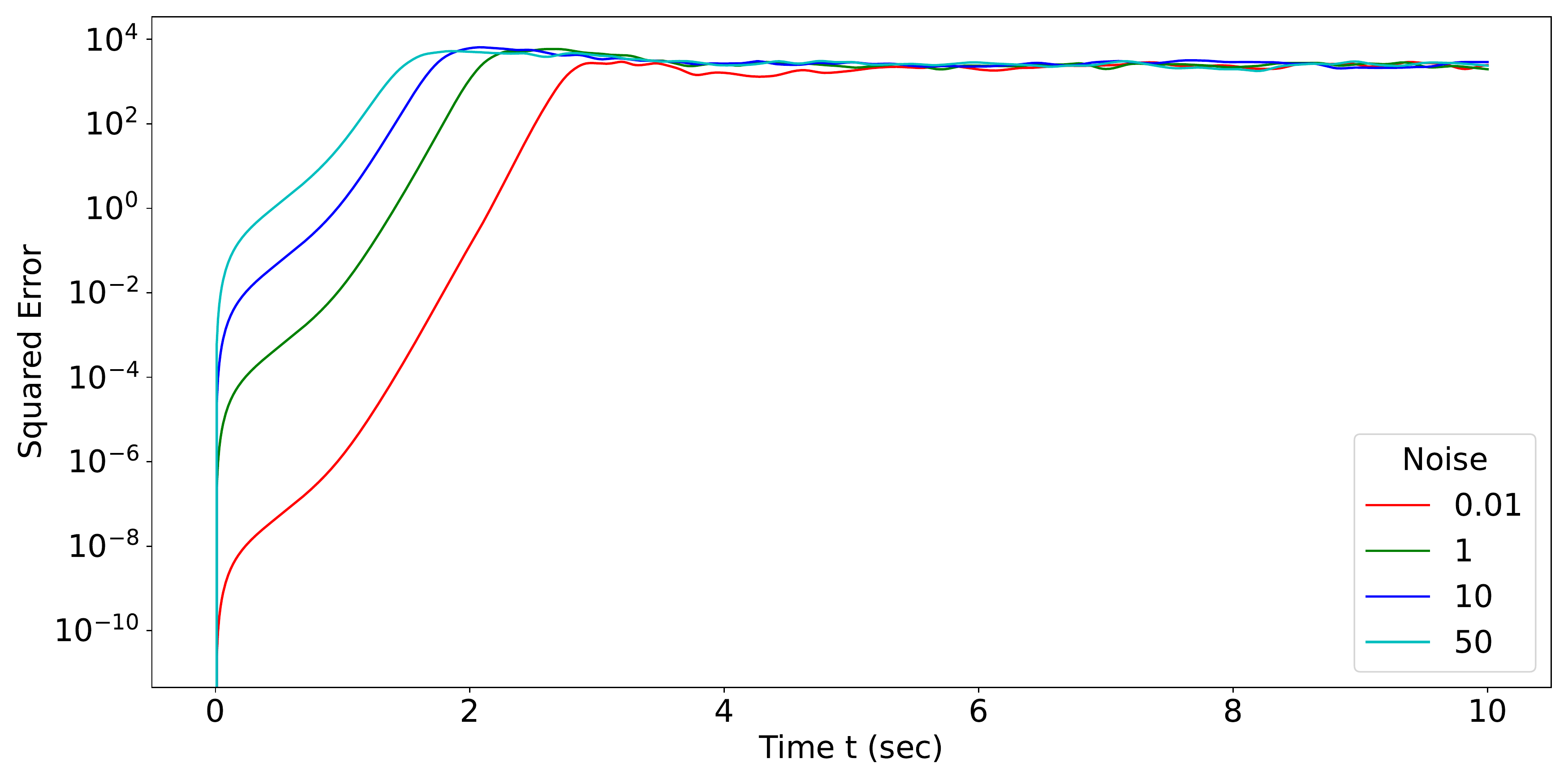}    \caption{$\ell_2$ error between the identified system and ground truth.}
    \label{fig:lorenz96_Gaussian_l2error}
\end{subfigure}
\caption{\footnotesize Simulation results of CS-MIO identified Lorenz 96 system using 60k data and under Type 1 noise with four noise magnitudes, namely 0.01, 1, 10 and 50. (a) Hoverm{\"o}ller plot for difference between the identified system and ground truth of Lorenz 96 system in $t \in[0,10]$. The vertical axis is the index $j$ of the state variable. The values of the colors refer to the difference between the ground truth states $x_j(t)$ and the evolved states $\hat{x}_j(t)$ using the identified equations by CS-MIO, i.e., $\Delta x_j(t) = x_j(t) - \hat{x}_j(t)$ for $j\in[96]$. (b) $l_2$ error vs time of the trajectories of the CS-MIO recovered Lorenz 96 system from $t=0$ to $t=10$. The exact recovery fails when $\sigma$ ls larger than 50.}    
\label{fig:lorenz96_Gaussian_simu}
\end{figure}

\begin{figure}[h!]
\begin{subfigure}[c]{0.75\textheight}
    \centering
    \includegraphics[width=5.4in]{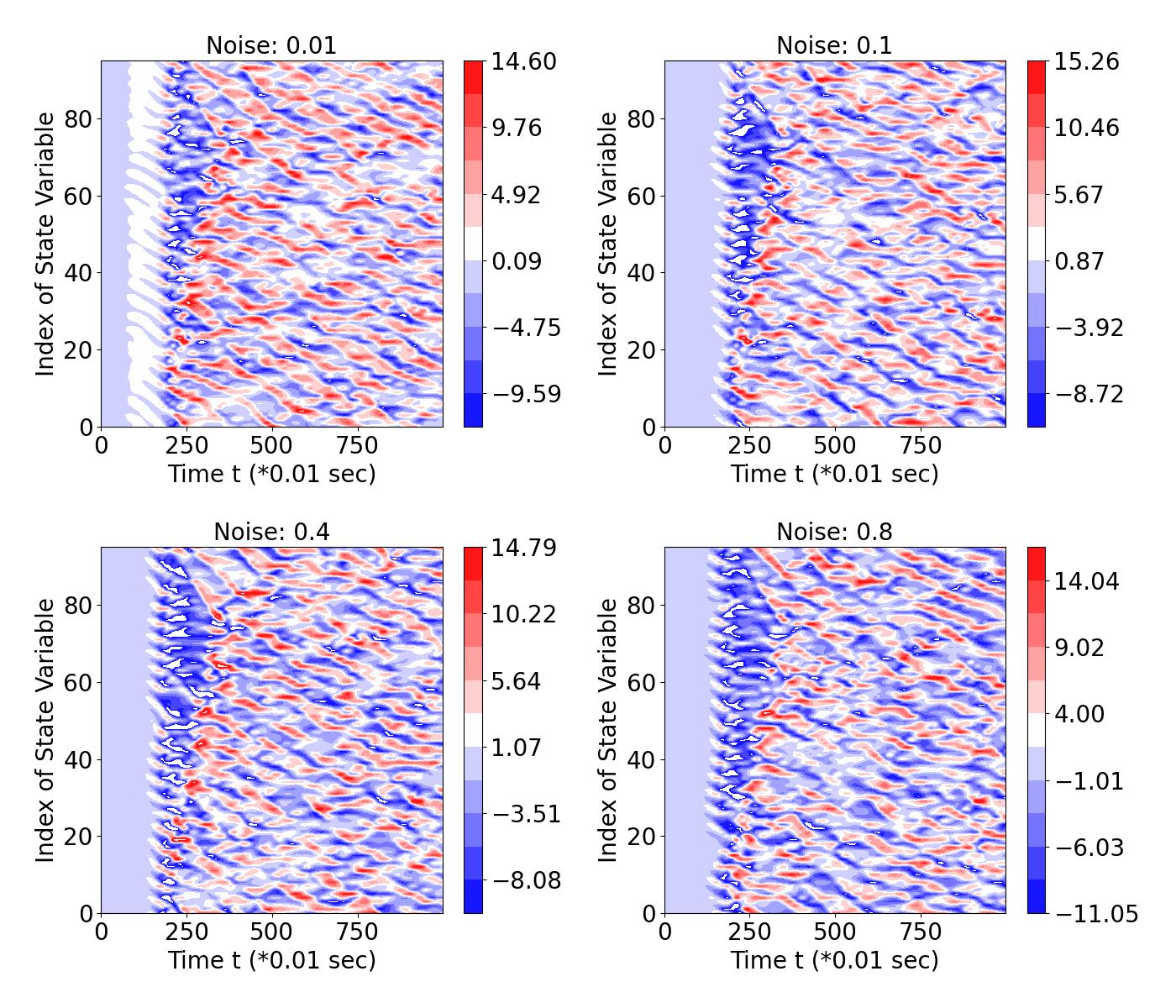}
    \caption{Hoverm{\"o}ller plot for difference between the identified system and ground truth.}
    \label{fig:lorenz96_TVD_diff}
\end{subfigure}
\begin{subfigure}[c]{1\textwidth}
\centering
    \includegraphics[width=4.7in]{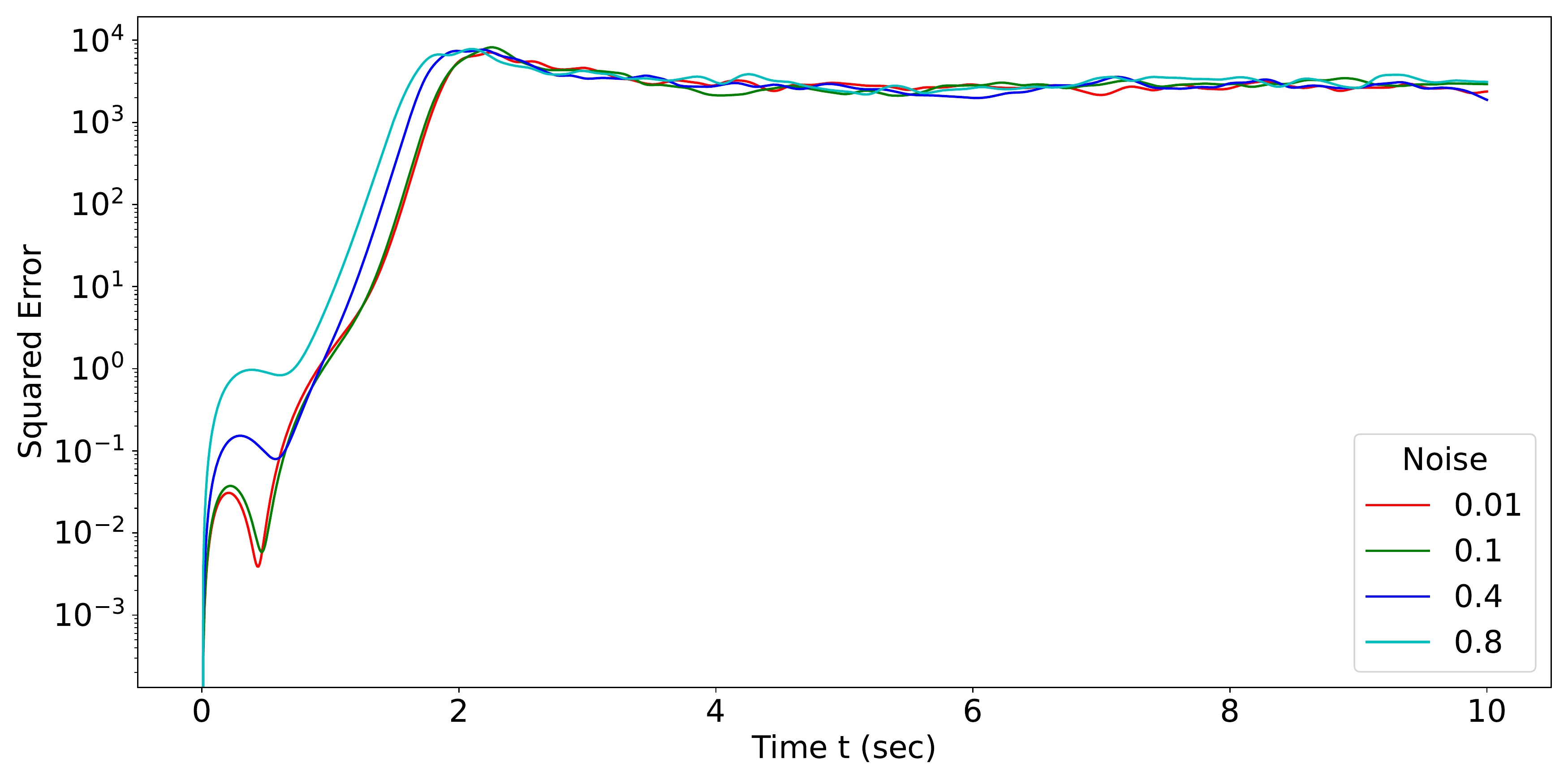}    
    \caption{$\ell_2$ error between the identified system and ground truth.}
    \label{fig:lorenz96_TVD_l2error}
\end{subfigure}
\caption{\footnotesize Simulation results of CS-MIO identified Lorenz 96 system using 60k data and under Type 2 noise with four noise magnitudes, namely 0.01, 0.1, 0.4, and 0.8. (a) Hoverm{\"o}ller plot for difference between the identified system and ground truth of Lorenz 96 system in $t \in[0,10]$. The vertical axis is the index $j$ of the state variable. The values of the colors refer to the difference between the ground truth states $x_j(t)$ and the evolved states $\hat{x}_j(t)$ using the identified equations by CS-MIO, i.e., $\Delta x_j(t) = x_j(t) - \hat{x}_j(t)$ for $j\in[96]$. (b) $l_2$ error vs time of the trajectories of the recovered Lorenz 96 system from $t=0$ to $t=10$. The exact recovery fails when $\sigma$ ls larger than 0.8.}
\label{fig:lorenz96_TVD_simu}
\end{figure}

\begin{table}[h]
\footnotesize
\centering
  \caption{Identified coefficients of Lorenz 96 system using CS-MIO under Type 2 noise with magnitude 0.8: Part A.}
  \label{tab:lorenz96_coefs_0.8_Gaussian_A}
  \begin{tabular}{c|rrrr}
    \toprule
      Equation Index $j$ &    $F$     &   $x_{j+1}x_{j-1}$   & $x_{j-2}x_{j-1} $  &  $x_j$ 	\\
     \midrule
 1	&	7.0568	&	0.9459	&	-0.9658	&	-0.6757	\\
2	&	7.0112	&	0.9483	&	-0.9639	&	-0.6958	\\
3	&	7.0591	&	0.9475	&	-0.9581	&	-0.6963	\\
4	&	7.0493	&	0.9571	&	-0.9597	&	-0.6776	\\
5	&	7.1518	&	0.9481	&	-0.9672	&	-0.7060	\\
6	&	7.0869	&	0.9523	&	-0.9635	&	-0.6799	\\
7	&	7.1948	&	0.9430	&	-0.9682	&	-0.6986	\\
8	&	7.1661	&	0.9491	&	-0.9674	&	-0.7077	\\
9	&	6.9831	&	0.9408	&	-0.9571	&	-0.6745	\\
10	&	7.0809	&	0.9559	&	-0.9623	&	-0.7037	\\
11	&	7.0943	&	0.9516	&	-0.9639	&	-0.6838	\\
12	&	7.2745	&	0.9574	&	-0.9713	&	-0.7263	\\
13	&	7.0001	&	0.9630	&	-0.9728	&	-0.6775	\\
14	&	7.1846	&	0.9497	&	-0.9675	&	-0.7080	\\
15	&	7.1597	&	0.9475	&	-0.9624	&	-0.6962	\\
16	&	7.0932	&	0.9526	&	-0.9670	&	-0.6986	\\
17	&	7.0520	&	0.9541	&	-0.9552	&	-0.7014	\\
18	&	6.8848	&	0.9549	&	-0.9643	&	-0.6501	\\
19	&	7.3256	&	0.9529	&	-0.9710	&	-0.7050	\\
20	&	7.1578	&	0.9634	&	-0.9677	&	-0.6955	\\
21	&	6.9812	&	0.9586	&	-0.9734	&	-0.6890	\\
22	&	6.9374	&	0.9463	&	-0.9665	&	-0.6560	\\
23	&	7.2987	&	0.9527	&	-0.9702	&	-0.7081	\\
24	&	7.0471	&	0.9486	&	-0.9610	&	-0.6813	\\
25	&	7.0503	&	0.9581	&	-0.9643	&	-0.7093	\\
26	&	6.9403	&	0.9514	&	-0.9609	&	-0.6487	\\
27	&	7.2672	&	0.9534	&	-0.9608	&	-0.7060	\\
28	&	6.9671	&	0.9564	&	-0.9652	&	-0.6765	\\
29	&	7.0891	&	0.9571	&	-0.9648	&	-0.6691	\\
30	&	7.0759	&	0.9487	&	-0.9712	&	-0.7038	\\
31	&	7.0496	&	0.9579	&	-0.9665	&	-0.6851	\\
32	&	6.9661	&	0.9421	&	-0.9614	&	-0.6288	\\
33	&	7.2776	&	0.9547	&	-0.9588	&	-0.7411	\\
34	&	6.8697	&	0.9531	&	-0.9710	&	-0.6503	\\
35	&	7.2938	&	0.9561	&	-0.9661	&	-0.7087	\\
36	&	6.9934	&	0.9543	&	-0.9610	&	-0.7026	\\
37	&	7.1151	&	0.9545	&	-0.9641	&	-0.6916	\\
38	&	7.1185	&	0.9554	&	-0.9673	&	-0.6762	\\
39	&	6.9873	&	0.9584	&	-0.9661	&	-0.6900	\\
40	&	7.0552	&	0.9493	&	-0.9635	&	-0.6868	\\
41	&	7.2372	&	0.9544	&	-0.9669	&	-0.7220	\\
42	&	6.9862	&	0.9515	&	-0.9619	&	-0.6670	\\
43	&	7.0901	&	0.9551	&	-0.9586	&	-0.6930	\\
44	&	7.0438	&	0.9531	&	-0.9620	&	-0.7080	\\
45	&	7.0553	&	0.9517	&	-0.9625	&	-0.6809	\\
46	&	7.0447	&	0.9509	&	-0.9646	&	-0.6727	\\
47	&	7.3077	&	0.9685	&	-0.9679	&	-0.7102	\\
48	&	7.0477	&	0.9698	&	-0.9671	&	-0.6916	\\
\bottomrule
    \end{tabular}
\end{table}

\begin{table}[h]
\footnotesize
\centering
  \caption{Identified coefficients of Lorenz 96 system using CS-MIO under Type 2 noise with magnitude 0.8: Part B.}
  \label{tab:lorenz96_coefs_0.8_Gaussian_B}
  \begin{tabular}{c|rrrr}
    \toprule
      Equation Index $j$ &    $F$     &   $x_{j+1}x_{j-1}$   & $x_{j-2}x_{j-1}$  &  $x_j$ 	\\
     \midrule
49	&	7.1351	&	0.9523	&	-0.9683	&	-0.6815	\\
50	&	7.0338	&	0.9526	&	-0.9637	&	-0.6808	\\
51	&	7.0436	&	0.9391	&	-0.9574	&	-0.6857	\\
52	&	7.0394	&	0.9566	&	-0.9659	&	-0.6886	\\
53	&	7.1219	&	0.9587	&	-0.9692	&	-0.6887	\\
54	&	7.1177	&	0.9464	&	-0.9631	&	-0.7365	\\
55	&	7.0162	&	0.9496	&	-0.9678	&	-0.6666	\\
56	&	7.2848	&	0.9554	&	-0.9675	&	-0.7262	\\
57	&	7.0027	&	0.9609	&	-0.9591	&	-0.6802	\\
58	&	7.0277	&	0.9523	&	-0.9640	&	-0.6903	\\
59	&	6.8512	&	0.9535	&	-0.9636	&	-0.6585	\\
60	&	7.2209	&	0.9483	&	-0.9651	&	-0.6831	\\
61	&	7.0277	&	0.9493	&	-0.9647	&	-0.6740	\\
62	&	7.0550	&	0.9545	&	-0.9701	&	-0.6907	\\
63	&	7.0577	&	0.9549	&	-0.9617	&	-0.6875	\\
64	&	7.1290	&	0.9515	&	-0.9595	&	-0.6763	\\
65	&	7.0478	&	0.9457	&	-0.9674	&	-0.6722	\\
66	&	7.1539	&	0.9407	&	-0.9648	&	-0.7014	\\
67	&	6.9451	&	0.9446	&	-0.9607	&	-0.6592	\\
68	&	7.0965	&	0.9588	&	-0.9702	&	-0.6837	\\
69	&	7.0340	&	0.9422	&	-0.9566	&	-0.6943	\\
70	&	7.0540	&	0.9546	&	-0.9654	&	-0.7106	\\
71	&	6.9663	&	0.9495	&	-0.9619	&	-0.6307	\\
72	&	7.2584	&	0.9489	&	-0.9604	&	-0.7156	\\
73	&	7.0076	&	0.9572	&	-0.9680	&	-0.6768	\\
74	&	7.1294	&	0.9450	&	-0.9585	&	-0.7171	\\
75	&	7.1466	&	0.9466	&	-0.9758	&	-0.6850	\\
76	&	7.1545	&	0.9560	&	-0.9680	&	-0.7093	\\
77	&	6.9454	&	0.9473	&	-0.9591	&	-0.6868	\\
78	&	7.0249	&	0.9494	&	-0.9603	&	-0.7016	\\
79	&	7.0989	&	0.9560	&	-0.9571	&	-0.6800	\\
80	&	7.1336	&	0.9661	&	-0.9697	&	-0.7055	\\
81	&	7.1208	&	0.9680	&	-0.9631	&	-0.6860	\\
82	&	7.1551	&	0.9685	&	-0.9713	&	-0.7370	\\
83	&	6.9564	&	0.9583	&	-0.9638	&	-0.6756	\\
84	&	7.1512	&	0.9496	&	-0.9649	&	-0.6969	\\
85	&	7.0909	&	0.9450	&	-0.9696	&	-0.6817	\\
86	&	6.9989	&	0.9461	&	-0.9636	&	-0.6907	\\
87	&	7.0673	&	0.9494	&	-0.9697	&	-0.6880	\\
88	&	7.0741	&	0.9495	&	-0.9611	&	-0.6648	\\
89	&	7.1304	&	0.9542	&	-0.9704	&	-0.7008	\\
90	&	7.0863	&	0.9455	&	-0.9587	&	-0.6852	\\
91	&	7.1075	&	0.9504	&	-0.9620	&	-0.7066	\\
92	&	7.0741	&	0.9534	&	-0.9761	&	-0.6834	\\
93	&	7.2857	&	0.9491	&	-0.9676	&	-0.6904	\\
94	&	6.8830	&	0.9580	&	-0.9587	&	-0.6762	\\
95	&	7.1757	&	0.9499	&	-0.9766	&	-0.7105	\\
96	&	7.1995	&	0.9550	&	-0.9631	&	-0.6822	\\
\bottomrule
    \end{tabular}
\end{table}

\FloatBarrier
\section{Additional results for the Hopf Normal Form}\label{sec:appendix_hopf}

\begin{table}[h]
\footnotesize
\centering
  \caption{Identified coefficients of Hopf normal form system using CS-MIO under Type 1 noise. }
  \label{tab:Hopf_Gaussian_coefs}
  \begin{tabular}{lr|rrrr}
    \toprule
      \multirow{2}{*}{Noise: $\sigma$}  &	\multirow{2}{*}{SNR} &   $y$   & $\mu x$  &  $x^3$  &   $xy^2$	\\
      & &   $x$ &   $\mu y$ &   $x^2y$  &   $y^3$   \\
     \midrule
 \multirow{2}{*}{0.001}	&	 \multirow{2}{*}{137583.905}	&	-1.0000	&	1.0000	&	-1.0000	&	-1.0000	\\
	&		&	1.0000	&	1.0000	&	-1.0000	&	-1.0000	\\
	 \midrule
 \multirow{2}{*}{0.01}	&	 \multirow{2}{*}{1375.839}	&	-1.0000	&	1.0000	&	-0.9999	&	-0.9998	\\
	&		&	1.0000	&	1.0000	&	-0.9995	&	-1.0000	\\
	 \midrule
 \multirow{2}{*}{0.1}	&	 \multirow{2}{*}{13.758}	&	-1.0001	&	1.0003	&	-0.9992	&	-0.9981	\\
	&		&	0.9998	&	0.9996	&	-0.9952	&	-1.0001	\\
	 \midrule
 \multirow{2}{*}{0.3}	&	 \multirow{2}{*}{1.529}	&	-1.0002	&	1.0009	&	-0.9975	&	-0.9942	\\
	&		&	0.9995	&	0.9988	&	-0.9857	&	-1.0003	\\
	 \midrule
 \multirow{2}{*}{0.5}	&	 \multirow{2}{*}{0.550}	&	-1.0004	&	1.0015	&	-0.9959	&	-0.9904	\\
	&		&	0.9991	&	0.9980	&	-0.9762	&	-1.0005	\\
	 \midrule
 \multirow{2}{*}{0.7}	&	 \multirow{2}{*}{0.281}	&	-1.0005	&	1.0021	&	-0.9942	&	-0.9866	\\
	&		&	0.9988	&	0.9972	&	-0.9667	&	-1.0008	\\
	 \midrule
 \multirow{2}{*}{1}	&	 \multirow{2}{*}{0.138}	&	-1.0007	&	1.0031	&	-0.9918	&	-0.9808	\\
	&		&	0.9983	&	0.9961	&	-0.9525	&	-1.0011	\\
	 \midrule
 \multirow{2}{*}{2}	&	 \multirow{2}{*}{0.034}	&	-1.0014	&	1.0061	&	-0.9835	&	-0.9616	\\
	&		&	0.9966	&	0.9921	&	-0.9049	&	-1.0022	\\
	 \midrule
 \multirow{2}{*}{3}	&	 \multirow{2}{*}{0.015}	&	-1.0021	&	1.0092	&	-0.9753	&	-0.9424	\\
	&		&	0.9949	&	0.9882	&	-0.8574	&	-1.0033	\\
\bottomrule
    \end{tabular}
\end{table}

\begin{table}[h!]
\footnotesize
\centering
  \caption{Identified coefficients of Hopf normal form system using CS-MIO under Type 2 noise.}
  \label{tab:Hopf_TVD_coefs}
  \begin{tabular}{lr|rrrr}
    \toprule
      \multirow{2}{*}{Noise: $\sigma$}  &	\multirow{2}{*}{SNR} &   $y$   & $\mu x$  &  $x^3$  &   $xy^2$	\\
      & &   $x$ &   $\mu y$ &   $x^2y$  &   $y^3$   \\
     \midrule
 \multirow{2}{*}{0.001}	&	 \multirow{2}{*}{120306.233}	&	-0.9951	&	0.9680	&	-0.9681	&	-0.9680	\\
	&		&	0.9951	&	0.9681	&	-0.9681	&	-0.9682	\\
	\midrule
 \multirow{2}{*}{0.003}	&	 \multirow{2}{*}{13367.359}	&	-0.9949	&	0.9545	&	-0.9542	&	-0.9543	\\
	&		&	0.9949	&	0.9555	&	-0.9554	&	-0.9552	\\
	\midrule
 \multirow{2}{*}{0.005}	&	 \multirow{2}{*}{4812.249}	&	-0.9948	&	0.9289	&	-0.9282	&	-0.9284	\\
	&		&	0.9947	&	0.9291	&	-0.9287	&	-0.9283	\\
	\midrule
 \multirow{2}{*}{0.007}	&	 \multirow{2}{*}{2455.229}	&	-0.9946	&	0.8925	&	-0.8913	&	-0.8915	\\
	&		&	0.9945	&	0.8913	&	-0.8905	&	-0.8898	\\
	\midrule
 \multirow{2}{*}{0.010}	&	 \multirow{2}{*}{1203.062}	&	-0.9942	&	0.8237	&	-0.8215	&	-0.8218	\\
	&		&	0.9940	&	0.8201	&	-0.8187	&	-0.8175	\\
	\midrule
 \multirow{2}{*}{0.013}	&	 \multirow{2}{*}{711.871}	&	-0.9937	&	0.7459	&	-0.7425	&	-0.7429	\\
	&		&	0.9934	&	0.7404	&	-0.7381	&	-0.7364	\\
	\midrule
 \multirow{2}{*}{0.015}	&	 \multirow{2}{*}{534.694}	&	-0.9932	&	0.6928	&	-0.6887	&	-0.6891	\\
	&		&	0.9929	&	0.6865	&	-0.6836	&	-0.6816	\\
\bottomrule
    \end{tabular}
\end{table}

\FloatBarrier
\section{Additional results for the logistic Map}\label{sec:appendix_logistic}

\begin{table}[h]
\centering
\footnotesize
  \caption{Identified coefficients of the logistic map system using CS-MIO.}
  \label{tab:logistic_Gaussian_coef}
  \begin{tabular}{rr|rr}
      \toprule
  Noise: $\sigma$  &	SNR &   $rx_n$  &   $rx_n^2$	\\
   \midrule
0.001	&	48377.506	&	1.0000	&	-1.0000	\\
0.01	&	481.877	&	0.9999	&	-0.9999	\\
0.1	&	8.985	&	0.9902	&	-0.9862	\\
0.2	&	3.619	&	0.9543	&	-0.9386	\\
0.3	&	2.146	&	0.9212	&	-0.8956	\\
0.4	&	1.455	&	0.8759	&	-0.8382	\\
0.5	&	1.098	&	0.8406	&	-0.7907	\\
0.6	&	0.874	&	0.8031	&	-0.7443	\\
\bottomrule
\end{tabular}
\end{table}

\end{appendices}

%\showmatmethods{} % Display the Materials and Methods section

% \acknow{Please include your acknowledgments here, set in a single paragraph. Please do not include any acknowledgments in the Supporting Information, or anywhere else in the manuscript.}
%
%\showacknow{} % Display the acknowledgments section

% Bibliography

\end{document}